\documentclass{article}

\usepackage{arxiv}

\usepackage[utf8]{inputenc} 
\usepackage[T1]{fontenc}    
\usepackage{hyperref}       
\usepackage{url}            
\usepackage{booktabs}       
\usepackage{amsfonts}       
\usepackage{nicefrac}       
\usepackage{microtype}      
\usepackage{doi}
\usepackage{style}
 \usepackage{enumitem}
 \usepackage[math]{cellspace}
\title{{\postprocess}: REcover Polynomial completeness ALGorithm for arbitrary polytopal basis functions}

\author{ \href{https://orcid.org/0000-0001-8642-4258}{\includegraphics[scale=0.06]{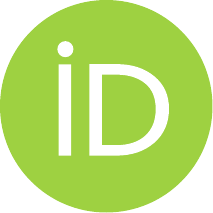}\hspace{1mm}Stefano~Berrone} \\
	Dipartimento di Scienze Matematiche\\
	``G. L. Lagrange'',\\
	Politecnico di Torino,  \\
    Corso Duca degli Abruzzi 24, 10129, Torino, Italy,\\
	\texttt{stefano.berrone@polito.it} \\
	\And
    \href{https://orcid.org/0000-0002-8481-0312}{\includegraphics[scale=0.06]{orcid.pdf}\hspace{1mm}Claudio~Canuto} \\
	Dipartimento di Scienze Matematiche\\
	``G. L. Lagrange''\\
	Politecnico di Torino,  \\
    Corso Duca degli Abruzzi 24, 10129, Torino, Italy,\\
	\texttt{claudio.canuto@polito.it} \\
	\And
	\href{https://orcid.org/0000-0002-2837-2792}{\includegraphics[scale=0.06]{orcid.pdf}\hspace{1mm}Moreno~Pintore} \\
	Laboratoire Jacques-Louis Lions, \\
	Sorbonne Universit\'e, \\
	MEGAVOLT Team, Inria, \\
	4 place Jussieu, 75005 Paris, France \\
	\texttt{moreno.pintore@sorbonne-universite.fr} \\
	\And
	\href{https://orcid.org/0000-0002-8540-3639}{\includegraphics[scale=0.06]{orcid.pdf}\hspace{1mm}Gioana~Teora} \\
	Dipartimento di Scienze Matematiche\\
	``G. L. Lagrange''\\
	Politecnico di Torino,  \\
    Corso Duca degli Abruzzi 24, 10129, Torino, Italy,\\
	\texttt{gioana.teora@polito.it} \\
}

\renewcommand{\shorttitle}{{\postprocess}: An Algorithm to Recover Polynomial Completeness}

\hypersetup{
pdftitle={RePAlg: An Algorithm to Recover Polynomial Completeness},
pdfsubject={q-bio.NC, q-bio.QM},
pdfauthor={Stefano~Berrone, Claudio~Canuto, Moreno~Pintore, Gioana~Teora},
pdfkeywords={Virtual Element Method, Neural Network, Basis Functions, Polygonal Meshes},
}

\begin{document}
\maketitle

\begin{abstract}
In this work, we introduce a general algebraic algorithm, called {\postprocess} (Recover Polynomial Algorithm), that transforms arbitrary closed-form non-polynomial shape functions over generic polygons into polynomial-complete bases. The proposed approach relies only on standard linear algebra operations, namely LU factorizations, and provides an efficient way to enhance existing approximation spaces without altering their essential structure. To prove the effectiveness of the algorithm, we introduce the high-order extension of the Boundary-Neural Approximated Virtual Element Method, where neural networks are used to learn virtual element basis functions over polygons. For this formulation, we derive convergence estimates before and after the application of the algorithm, proving that the proposed correction eliminates the limiting effect of the neural-network approximation error and restores optimal polynomial convergence rates. Numerical experiments confirm the theoretical results and assess the impact of quadrature errors associated with these non-polynomial basis functions.
\end{abstract}


\section{Introduction}

The simulation of physical phenomena plays a crucial role in many engineering applications, where the numerical solution of partial differential equations is often required. Among the available numerical techniques, the Finite Element Method (FEM) \cite{ciarlet2002finite, brennerscott} is one of the most widely used, owing to its solid theoretical foundations and computational flexibility. Standard FEM formulations, however, are constrained by the requirements imposed on the domain discretization: in two-dimensional problems, such as those considered in this work, the mesh is typically restricted to triangular or convex quadrilateral elements.

Over the last few decades, several strategies have been proposed to broaden the class of admissible elements. Among them, Discontinuous Galerkin (DG) methods \cite{DiPietro2011, lymph}, Hybrid High-Order (HHO) methods \cite{basicHHO}, and Virtual Element Methods (VEM) \cite{LBe13, LBe14} all circumvent the difficulty of explicitly constructing high-order shape functions on generic polygons: DG methods resort to fully discontinuous polynomial spaces coupled with suitable numerical fluxes, whereas HHO and VEM rely on local projection operators possibly supplemented by a stabilization term.

A complementary line of research has instead pursued the explicit construction of basis functions directly on polygonal elements. Since the pioneering work of Wachspress \cite{Wachspress1975}, who first defined rational nodal shape functions on convex polygons that form a partition of unity and reproduce linear polynomials exactly (i.e. are linearly complete), numerous studies have sought to generalize this idea \cite{SukumarTabarraei2004}. A non-exhaustive list of representative works, limited to low-order basis functions, includes Mean Value Coordinates (MVC) \cite{FLOATER2003, FLOATER2005}, Maximum Entropy coordinates (MAXENT) \cite{Sukumar2004}, Natural Neighbor-based coordinates \cite{Sukumar1998, Sukumar2000, talischi2010polygonal}, and the Smoothed Finite Element Method (SFEM) \cite{liu2007, dai2007, Liu2018}. More recently, a related class of methods has made the otherwise implicit basis functions of the VEM explicitly computable, granting access to their point-wise values: this is the case of the Lightning Virtual Element Method (LVEM) \cite{TrezziZerbinati2024}, which approximates the virtual basis functions through the Laplace Solver \cite{Gopal2019}, the reduced-basis Virtual Element Method (rbVEM) \cite{Credali2024}, that exploits reduced order models and the Neural Approximated Virtual Element Method (NAVEM) \cite{PintoreTeora2025, navemElasticity, PintoreTeora2026}.

Extending these explicit constructions to higher order, possibly in the presence of hanging nodes, is considerably more involved. Nevertheless, the superior accuracy and convergence properties of high-order approximations have long motivated research efforts in this direction \cite{Duczek2016}. A non-exhaustive list includes the Spectral Element Method (SEM) \cite{semBook, semReview},  $hp$-methods \cite{Schwab1998}, Serendipity Finite Element Methods \cite{arbogast2022direct, arbogast2023direct}, transfinite elements \cite{EISENTRAGER2021, provatidis2025Macro, provatidis2025}, and the Zipped Finite Element Methods \cite{Bunge2020, TeoraNeva2025}.

In the present manuscript, we propose {\postprocess}, a general algebraic algorithm that enables the construction of polynomial-complete basis functions starting from arbitrary shape functions defined over generic polygonal elements and available in \textit{closed-form}, i.e., explicitly known and directly evaluable at any point of the element without requiring any additional approximation or numerical reconstruction. The resulting basis preserves the key structural properties of the original functions, including their Lagrange interpolation property and the minimal continuity required by the underlying discretization, while guaranteeing that the polynomial space of the prescribed degree is contained in their span. Since the proposed construction relies exclusively on standard linear algebra operations, namely the LU factorization, it is easy to implement and readily applicable to any family of shape functions.

To assess the generality and effectiveness of the proposed approach, we apply it to two different families of non-polynomial basis functions. The first one consists of a high-order lifting of Wachspress coordinates enriched with bubble functions \cite{Sukumar2026}. The second one is a novel high-order extension of the Boundary-Neural Approximated Virtual Element Method (B-NAVEM) \cite{PintoreTeora2026}, in which neural networks are employed to learn Virtual Element basis functions \cite{LBe13}. The latter represents a particularly compelling application of the proposed framework, as the approximation properties of the numerical method are directly influenced by the accuracy of the neural-network approximation. Indeed, although the learned basis functions satisfy the required conformity conditions, the approximation errors introduced by the neural networks generally prevent the resulting space from exactly reproducing polynomials, thereby limiting the asymptotic convergence rates. The proposed polynomial enrichment algorithm overcomes this limitation. Rather than replacing the neural-network approximation, it acts as an algebraic post-processing step that preserves the structural properties of the learned basis, while enforcing polynomial completeness. Consequently, the approximation space inherits the optimal polynomial approximation properties independently of the approximation error introduced by the neural network. For the proposed high-order B-NAVEM formulation, we derive rigorous a priori error estimates both before and after the application of the enrichment algorithm. The analysis explicitly characterizes the effect of the neural-network approximation error on the convergence of the method and proves that, once the proposed algebraic correction is applied, this limiting term is eliminated.

The outline of the paper is as follows. In Section \ref{sec:problem_formulation}, we introduce the model problem and formulate a general discretization framework over polygonal elements based on closed-form representations of the underlying basis functions. In Section \ref{sec:discrete_space}, we present the proposed algebraic algorithm, based on the LU factorization, and describe its application to the high-order extension of the B-NAVEM method. Section \ref{sec:apriori} provides the a priori error estimates for the space spanned by the basis functions obtained through the proposed post-processing {\postprocess} procedure. In Section \ref{sec:boundary_funcs}, we introduce a procedure to construct functions satisfying prescribed boundary conditions exactly, and we show how this construction can be exploited to impose the continuity of B-NAVEM basis functions and to generate additional families of approximation bases. Finally, Section \ref{sec:numeri_experiements} presents a wealth of numerical experiments validating the theoretical results and investigating the influence of quadrature errors on the accuracy of the method, whereas \ref{app:qr} shows how alternative factorizations can be used to restore polynomials.

Throughout the paper, $C$ will denote a positive quantity independent of the relevant parameters, which may be different in different occurrences.

\section{The model problem and its discretization} \label{sec:problem_formulation}

In the following, we adopt standard notation for Sobolev spaces. Given a generic Sobolev space $\mathcal{H}$, we use the symbols $\norm[\mathcal{H}]{}$ and $\seminorm[\mathcal{H}]{}$ to indicate the norm and the seminorm in $\mathcal{H}$, resp., and $\scal[\mathcal{H}]{}{}$ to denote its scalar product. When $\mathcal{H}=\leb{2}{\genericset}$, with $\genericset$ being a generic open subset $\genericset \subset \R^2$, we adopt a lighter notation, denoting by $\scal[\genericset]{}{}$ the scalar product in $\leb{2}{\genericset}$; precisely, given two scalar functions $p,q \in \leb{2}{\genericset}$ and two vector fields $\bm{a},\ \vv \in \vleb{2}{\genericset}{2}$, we define
\begin{align*}
\scal[\genericset]{p}{q} = \int_{\genericset} p q,\qquad\scal[\genericset]{\bm{a}}{\vv} = \int_{\omega} \bm{a}\cdot \vv.
\end{align*}
Moreover, for a vector $\xx \in \R^n$ and matrix $\mathbf{A} \in \R^{n \times m}$, for given positive integer numbers $n, m$, the symbols $\norm[2]{\xx}$ and $\norm[2]{\mathbf{A}}$ denote the Euclidean and the spectral norms, resp., whereas $\norm[\infty]{\mathbf{A}}$ denotes the infinity norm of matrix $\mathbf{A}$.

Given an integer $k \geq 1$ and a polygon $E$, we denote by $\Poly{k}{E}$ the space of polynomials of total degree up to $k$ defined on $E$, with $n_k = \dim \Poly{k}{E} = \frac{(k+1)(k+2)}{2}$. We further adopt the standard convention $ \Poly{-1}{E} = \{0\}$ and $n_{-1} = 0$. Moreover, we denote by $m_{\alpha}$ the scaled monomial 
\begin{equation}
    m_{\alpha}(\xx) = \left(\frac{\xx - \xx_E}{h_E}\right)^{\bm{\alpha}},
    \label{eq:monomial}
\end{equation}
where $\xx_E$ and $h_E$ are the centroid and the diameter of $E$, resp., $\bm{\alpha} \in \N^2$ is a multi-index s.t. $\vert \bm{\alpha} \vert \leq k$ and $\alpha = \xi(\bm{\alpha}) \in \{1,\dots, n_{k}\}$ with $\xi$ denoting the classical Cantor diagonal counting of the elements of $\N^2$
\begin{equation*}
    0 \leftrightarrow (0, 0),\quad 1 \leftrightarrow (1, 0),\quad 2 \leftrightarrow (0, 1),\quad 3 \leftrightarrow (2, 0),\dots.
\end{equation*}

Let us consider a convex open bounded polygonal domain $\Omega\subset\R^2$ with boundary $\Gamma=\partial\Omega$ and the following boundary-value problem:
\begin{equation}
    \begin{cases}
        -\Delta u = f & \text{in }\Omega,\\
        u  = 0&\text{on }\Gamma,
    \end{cases}
    \label{eq:continuous_poisson}
\end{equation}
where $f\in \leb{2}{\Omega}$ is the problem forcing term. 
Introducing the Hilbert space $\V \coloneqq \sob[0]{1}{\Omega}$, the variational formulation of Problem \eqref{eq:continuous_poisson} reads as: \textit{Find $u\in \V$ such that}
\begin{equation}
    \dbilin{u}{v} = \rightlin{v} \quad \forall v\in \V,
    \label{eq:exact_var_pb}
\end{equation}
where the bilinear form $\dbilin{}{}:\V\times \V\to\R$ and the linear form $\rightlin{}:\V\to\R$ are defined as
\begin{equation*}
    \dbilin{v}{z} = \int_\Omega \nabla v\cdot\nabla z, \quad \forall v, z\in \V,\qquad \quad
\rightlin{v} = \int_\Omega f \,v, \quad \forall v\in \V.
\end{equation*}

\begin{remark}
    While we consider Problem \eqref{eq:continuous_poisson} for the sake of simplicity, our theory applies to any variational problem of the form \eqref{eq:exact_var_pb}, where $\V$ is a Hilbert space satisfying $\sob[0]{1}{\Omega} \subseteq \V \subseteq \sob{1}{\Omega}$, $a : \V \times \V \to \R$ is a continuous and coercive bilinear form, and $\mathcal{F}: \V \to \R$ is a linear form.
\end{remark}

Let $\Th$ be a polygonal mesh of $\overline{\Omega}$ made up of polygonal elements $E$ satisfying suitable assumptions. In the following, $\Nv[E]$ and $\Eh[E]$, resp., will denote the number of vertices (and edges) of $E$, and its set of edges, resp.. Moreover, we set $h = \max_{E \in \Th} h_E$. Let us consider the following standard mesh assumptions \cite{LBe17}.
\begin{assumption}[Mesh assumptions]
\label{ass:meshassumption}
There exists a positive real number $\rho \in (0, 1)$, independent of $h$, such that
    \begin{itemize}
        \item for each element $E\in\Th$ and each edge $e \in \Eh[E]$, it holds $\vert e \vert \geq \rho h_E$;
        \item each element $E\in\Th$ is star-shaped with respect to a ball of radius $\geq \rho h_E$;
        \item $\Th$ is a quasi-uniform mesh, i.e. $h_E \approx h$ for each $E \in \Th$.
    \end{itemize}
\end{assumption}

Given an integer $k \geq 1$, our goal is to introduce a local discrete space $\pVh[E]{k}$ for each $E\in \Th$, such that its global counterpart
\begin{equation}
    \pVh{k} \coloneq \{v \in \V \cap \con{0}{\overline{\Omega}}: v_{|E} \in \pVh[E]{k}\ \forall E \in \Th\},
\end{equation}
contains functions that are available in closed-form and satisfies
\begin{equation}
   \Poly{k}{\Omega; \Th} \subseteq \pVh{k},
   \label{eq:polynomial_reproducibility}
\end{equation}
where
\begin{equation*}
    \Poly{k}{\Omega; \Th} \coloneq \{p \in V : p_{|E} \in \Poly{k}{E}\  \forall E \in \Th\}.
\end{equation*}

Assuming that such a space $\pVh{k}$ is available, Problem \eqref{eq:exact_var_pb} can be discretized as in a standard finite element solver, that is: \textit{Find $u_h^{\NAPEM}\in \pVh{k}$ such that}
\begin{equation}
    \dbilin{u_h^{\NAPEM}}{v_h^{\NAPEM}} = \rightlin{v_h^{\NAPEM}} \quad \forall v_h^{\NAPEM}\in \pVh{k}.
    \label{eq:discrete_var_pb_p}
\end{equation}
\begin{remark}[The role of quadrature formula]\label{rem:quadrature}
In Equation \eqref{eq:discrete_var_pb_p}, we assume that all integrals are computed exactly. In practical computations, however, these integrals are generally approximated by means of numerical quadrature, thereby introducing an additional source of error. Throughout the theoretical analysis, we neglect this contribution and assume exact integration. The effect of quadrature errors is instead investigated in the numerical experiments reported in Section \ref{sec:numeri_experiements}.
\end{remark}

\section{The discrete local spaces} \label{sec:discrete_space}

To build a space $\pVh[E]{k}$ for an element $E \in \Th$, we can adopt, for instance, a B-NAVEM strategy \cite{PintoreTeora2026}, which defines basis functions that are known in closed-form and continuous over $\overline{\Omega}$ as the approximation of standard virtual basis functions, and then enforce a polynomial reproducibility property through {\postprocess} (see Section~\ref{sec:pol_reproducibility}).

In the same spirit of the finite element construction, we can exploit the existence of an \emph{inertia element} $\widehat{E}$ and an \emph{inertia map} $F_{\mathcal{I}}^{-1}: E \to \widehat{E}$ to build a local space, whose properties depend as weakly as possible on the geometry of the underlying element.  By definition, an inertia element is a polygon with unit diameter, is centered at the origin, and has an inertia tensor equal to a constant times the identity tensor \cite{Teora2024, PintoreTeora2024}. The inertia map $F_{\mathcal{I}}^{-1}$, introduced in \cite{Teora2024}, is a smooth bijection between the two elements: the original element $E$ and the associated inertia element $\widehat{E}$. In the following, the symbol $\ \widehat{q}\ $ will denote a quantity $q$ related to the inertia element. Nonetheless, when no confusion arises, we omit this symbol for the ease of notation.

\subsection{The local virtual element space}\label{sec:vem}

Given an integer $k \geq 1$ and an element $E \in \Th$, the local virtual element space is defined as \cite{LBe16}
\begin{align}
     \label{eq:vem:harmonicity}\Vh[E]{k} \coloneq \Big\{ v \in \sob{1}{E}:\quad &(i)\ \Delta v \in \Poly{k-2}{E}, \\ 
     \label{eq:vem:linearity}&(ii) \ v \in \con{0}{\partial E} : v_{|e}\in \Poly{k}{e},\ \forall e \subset \partial E\Big\}, 
\end{align}
where $\partial E$ is the boundary of $E$. We observe that $\Poly{k}{E} \subseteq \Vh[E]{k}$.

We remark that the VEM space employed herein is the first VEM space introduced in \cite{LBe13}, and not the \textit{enhanced} space that is commonly adopted within the virtual element method framework \cite{LBe16}. This choice is motivated by the fact that B-NAVEM does not require the computation of polynomial projections (since basis functions are explicitly known), making the enhancement property unnecessary.

The local degrees of freedom (DOFs in short) for a function $v \in \Vh[E]{k}$ are defined as
\begin{enumerate}[label=D.\arabic*]
	\item\label{dof:vertex} the values of $v$ at the vertices of $E$;
	\item\label{dof:edge} if $k \geq 2$, the values of $v$ at the $k-1$ Gauss-Lobatto quadrature points internal to each edge of $E$;
	\item\label{dof:internal} if $k \geq 2$, the $n_{k-2}$ internal scaled moments defined over the inertia element $\widehat{E} = F_{\mathcal{I}}^{-1}(E)$: $\frac{1}{\vert \widehat{E} \vert}\int_{\widehat{E}} \widehat{v}\  \widehat{m}_{\alpha}$ for all $\alpha = 1,\dots,n_{k-2}$.
\end{enumerate}
The number of local degrees of freedom is
\begin{equation*}
    \Ndof[E] = k \Nv[E] + n_{k-2}.
\end{equation*}

\begin{remark}\label{rem:inftynorm}
We highlight that the local degrees of freedom \ref{dof:vertex}-\ref{dof:internal} are chosen in such a way that, for any function $v \in \Vh[E]{k}$, it holds true
\begin{equation*}
    \norm[\infty]{\bm{v}} \coloneq \max_{j} \vert v_j \vert \leq \norm[\leb{\infty}{E}]{v},
\end{equation*}
where $\bm{v} \in \R^{\Ndof[E]}$ is the vector collecting the degrees of freedom $\{v_j\}_{j=1}^{\Ndof[E]}$ of the function $v$.
\end{remark}

For a generic element $E$, let $\dof_j : \Vh[E]{k} \to \R$ be the linear operator that associates to each function $v \in \Vh[E]{k}$ its $j$-th degree of freedom, where $j$ first counts the vertex DOFS \ref{dof:vertex}, then the edge DOFs \ref{dof:edge} and finally the internal DOFs \ref{dof:internal}. In the following, this ordering will be referred to as the \textit{natural ordering} of DOFs. Let $\{\vvarphi_{j,E}\}_{j=1}^{\Ndof_E}$ be the $\Ndof_E$ Lagrange basis functions of $\Vh[E]{k}$ associated with the degrees of freedom \ref{dof:vertex}--\ref{dof:internal}. These functions satisfy
\begin{equation*}
    \dof_i(\vvarphi_{j,E}) = \delta_{ij}\quad \forall i,j =1,\dots,\Ndof[E]
\end{equation*}
and are the solutions of the following elemental problem:
\begin{equation}
\begin{cases}
\Delta \vvarphi_{j,E} = g^{\Delta}_{j,E} & \text{in } E,\\
\vvarphi_{j,E} = g^{\partial}_{j,E} & \text{on } \partial E,
\end{cases}
\label{eq:dfproblem_vem_basis}
\end{equation}
where $g^{\Delta}_{j,E}$ are suitable known polynomials of degree $k-2$ in $E$ and $g^{\partial}_{j,E}$ are known continuous functions on $\partial E$, which restrict to a polynomial of degree $k$ over each edge belonging to the boundary of $E$.

In the virtual element framework, the elemental problem \eqref{eq:dfproblem_vem_basis} is never solved, and proper polynomial projectors and stability operators are introduced to mimic the point-wise evaluations of virtual functions and of their derivatives and to recover the coercivity of the discrete problem, causing different issues when dealing with strongly anisotropic problems or non-linear problems \cite{PintoreTeora2025, MarconTeora2023}, or in the post-processing phase of the quantities of interest \cite{Credali2024}.

Finally, the global space is given by
\begin{equation*}
    \Vh{k} \coloneq \{v \in \V : v_{|E} \in \Vh[E]{k} \ \forall E \in \Th\}.
\end{equation*}
The dimension of this space is given by
\begin{equation*}
    \Ndof \coloneq \dim \Vh{k} = \Nv + (k - 1) \Ne + \Nc n_{k-2},
\end{equation*}
where $\Nv$ and $\Ne$ are the number of vertices and edges of the tessellation not belonging to the Dirichlet boundary $\Gamma$, whereas $\Nc$ denotes the number of polygons belonging to $\Th$. 

Under Assumption \ref{ass:meshassumption}, we note that in a quasi-uniform mesh, the number of vertices, edges, and elements scales as
\begin{equation*}
    \Nv \approx \Ne \approx \Nc \approx \frac{1}{h^2}.
\end{equation*}
As a consequence, the number of total degrees of freedom scales as
\begin{equation*}
    \Ndof \approx \frac{1}{h^2} + \frac{1}{h^2} (k -1) + \frac{1}{h^2} \frac{(k-1)k}{2} \approx \frac{k^2}{h^2}.
\end{equation*}

\subsection{The B-NAVEM space} \label{sec:navem_space}

Now, let us describe the B-NAVEM strategy, introduced in \cite{PintoreTeora2026}, to approximate virtual basis functions with functions that are known in closed-form.

For each class of polygons $\mathcal{C}$ identified by a given number of vertices $\Nv[\mathcal{C}] \geq 4$, we build a dataset $\mathcal{D}_{\mathcal{C}}$ composed of a copious number of inertia elements $E$ with $\Nv[\mathcal{C}]$ vertices. Then, consider a fully-connected feed-forward neural network $\xx_0 \mapsto {\mathcal N}(\xx_0)$ given by
\begin{equation} 
    \begin{aligned}
        &\xx_\ell = \rho(A_\ell \xx_{\ell-1} + b_\ell), \hspace{2cm} \ell = 1,...,L-1, \\
        &{\mathcal N}(\xx_0) = A^{L} \xx_{L-1}  + b_L,
    \end{aligned}
    \label{eq:nn_formula}
\end{equation}
where  $L\in\N$ is the number of layers, $A_\ell\in\R^{N_\ell\times N_{\ell-1}}$ and $b_\ell\in\R^{N_\ell}$ are the matrices and vectors containing the trainable weights, and $\rho$ is a nonlinear activation function; here, we always set $\rho \equiv \tanh$. The vector $\xx_0$ encodes the triplet $(\xx,j,E)$, where $E$ is an inertia element, $j$ is the index of the degree of freedom related to the basis functions on $E$ we aim to predict, and $\xx$ is an evaluation point inside $E$. Thus, the input dimension, which must store the 2 coordinates of the evaluation point, the $2\Nv[E]$ coordinates of the vertices, and the index $j$, will be $2 + 2 \Nv[E] + 1$, whereas the output dimension is always chosen to be $1$, i.e. $N_L = 1$.

Assume for a moment that the neural network has been suitably trained. Then, for any inertia element $E$ and every index $j=1,\dots,\Ndof[E]$, consider the \emph{sliced neural network} defined on $E$ by
\begin{equation*}
{\mathcal N}_{j,E}(\xx) := {\mathcal N}\big( (\xx,j,E)\big), \qquad \forall \xx \in E,
\end{equation*}
and define the $j$-th basis function $\bvarphi_{j,E}$ on  $E$ (designed to approximate $\varphi_{j,E}^V$) by
\begin{equation}  
        \bvarphi_{j,E} :=\mathcal{B}_{j,E}({\mathcal N}_{j,E}),
        \label{eq:varphi_j_n}
\end{equation}        
where the operator $\mathcal{B}_{j,E}: \con{0}{E}\to \con{0}{E}$ is defined such that
the function $\mathcal{B}_{j,E}(v)$ satisfies
\begin{equation}
    \mathcal{B}_{j,E}(v) = \varphi^{\VEM}_{j,E} = {g}_{j,E}^{\partial}\quad \text{on } \partial {E}.
    \label{eq:boundary_prop_boperator}
\end{equation}
Moreover, noting that $\mathcal{N}_{j,E} \in\con{\infty}{E}$ due to the choice of the activation function, we assume that the application of $\mathcal{B}_{j,E}$ to the neural network is regular enough to guarantee that the B-NAVEM basis functions are at least as regular as the virtual basis functions. More specifically, we assume the following hypothesis holds true.
\begin{assumption}\label{ass:navem-regularity}
    For each inertial element $E$, $\bvarphi_{j,E} \in \con{0}{E} \cap \sob{1}{E}$ and $\Delta \bvarphi_{j,E} \in \leb{2}{E}$ with $\norm[\leb{2}{E}]{\Delta \bvarphi_{j,E}} \leq C$ uniformly in $E$. 
\end{assumption}

Different choices for constructing such an operator will be introduced in Section \ref{sec:boundary_funcs}.

At this point, we have to explain how the neural network \eqref{eq:nn_formula} is trained. To this end, for any element $E$ in the dataset $\mathcal{D}_{\mathcal{C}}$ and any index $j =1,\dots,\Ndof[E]$, let us introduce the local loss function
\begin{equation}
    {\epsilon}_{j,E}^{\BNAVEM} \coloneq \norm[\leb{2}{{E}}]{{\Delta} {\varphi}^{\BNAVEM}_{j,E} - {g}_{j,E}^{{\Delta}}} .
    \label{eq:bnavem:loss-jE}
\end{equation}
For each class of polygons ${\cal C}$, the neural network is trained to minimize the quadratic mean of these quantities, namely the global loss function
\begin{equation}
    \mathcal{L}^{\BNAVEM}_{\mathbf{C}} = \sqrt{\frac{1}{\Ndof[{\cal C}] \# \mathcal{D}_{\mathcal{C}}} \sum_{E \in \mathcal{D}_{\mathcal{C}}} \sum_{j=1}^{\Ndof[{\cal C}]} \left({\epsilon}_{j,E}^{\BNAVEM} \right)^2},
\label{eq:bnavem:loss}    
\end{equation}
where $\Ndof[{\cal C}]$ denotes the number of degrees of freedom common to all polygons in the class ${\cal C}$.

At last, we deal with the case of a polygon $E$ in the class ${\cal C}$, which is not an inertia element. 
In the same spirit of the finite element construction, we rely on the previous construction on inertia elements, by defining the $j$-th basis function $\bvarphi_{j,E}$ on the original element $E$ as follows:
\begin{equation}  
    \bvarphi_{j,E}(\xx) \coloneq \widehat{\varphi}_{j,\widehat{E}}^{\BNAVEM}(\widehat{\xx}), \qquad \widehat{\xx}= F_{\mathcal{I}}^{-1}(\xx).
    \label{eq:varphi_j_n}
\end{equation}
Under Assumptions \ref{ass:meshassumption} and the properties of the inertia map, we have
\begin{gather}
    \nonumber C h_E \norm[\leb{2}{\widehat{E}}]{\widehat{\varphi}^{\BNAVEM}_{j,\widehat{E}}}\leq \norm[\leb{2}{E}]{\bvarphi_{j,E}} \leq C h_E \norm[\leb{2}{\widehat{E}}]{\widehat{\varphi}^{\BNAVEM}_{j,\widehat{E}}},\\
     \nonumber   C  \norm[\leb{2}{\widehat{E}}]{\widehat{\nabla} \widehat{\varphi}^{\BNAVEM}_{j,\widehat{E}}}\leq \norm[\leb{2}{E}]{\nabla \bvarphi_{j,E}} \leq C  \norm[\leb{2}{\widehat{E}}]{\widehat{\nabla}\widehat{\varphi}^{\BNAVEM}_{j,\widehat{E}}},\\
      \label{eq:scaling-lapl}  C  h_E^{-1}\norm[\leb{2}{\widehat{E}}]{\widehat{\Delta} \widehat{\varphi}^{\BNAVEM}_{j,\widehat{E}}}\leq \norm[\leb{2}{E}]{\Delta \bvarphi_{j,E}} \leq C h_E^{-1}  \norm[\leb{2}{\widehat{E}}]{\widehat{\Delta}\widehat{\varphi}^{\BNAVEM}_{j,\widehat{E}}}.
\end{gather}

As a result, training and prediction are always carried out on the inertia element $\widehat{E}$, and, consequently, the basis functions and their gradients are predicted on $\widehat{E}$ and then mapped back to the original element $E$.

\begin{remark}
    For the class of triangles, i.e. $\Nv[\mathcal{C}] = 3$, we do not train a neural network, since we can adopt the set of finite element basis functions that are known in a closed-form for each order $k \geq 1$.
\end{remark}
\begin{remark}
    We recall that operating on inertia-mapped elements reduces the geometric variability of the inputs presented to the neural networks, thereby improving their approximation accuracy \cite{PintoreTeora2024}.
\end{remark}

Finally, the B-NAVEM local space is defined as 
\begin{equation}
    \bVh[E]{k} \coloneq \myspan\{\bvarphi_{j,E}:\, j=1,\dots,\Ndof[E]\}.
    \label{eq:localspacebnavem}
\end{equation}

\subsubsection{Properties of B-NAVEM functions}

The B-NAVEM strategy produces functions that satisfy Property \eqref{eq:vem:linearity} and are continuous across $\overline{\Omega}$. Property \eqref{eq:vem:harmonicity} is instead only weakly satisfied through the neural network. As a consequence, generally, we have
\begin{equation*}
    \Poly{k}{E} \not\subseteq \bVh[E]{k}.
\end{equation*}

Next theorem shows that B-NAVEM functions are close to VEM basis functions if the neural network is sufficiently accurate. In particular, this property implies that B-NAVEM functions are linearly independent (see Theorem \ref{theor:linearlyind}).

\begin{theorem}\label{theor:bnavem:local_basis}
For a given pair $(j,E) \in \{1,\dots,\Ndof[E]\} \times \Th$, let $\varphi_{j,E}^{\VEM}$ and $\varphi_{j,E}^{\BNAVEM}$ be the VEM basis function and its B-NAVEM approximation, respectively. Then
    \begin{equation*}
        \norm[\leb{2}{E}]{\nabla \varphi_{j,E}^{\BNAVEM} - \nabla \varphi_{j,E}^{\VEM}} \leq C h_E   \epsilon_{j,E}^{\BNAVEM},
    \end{equation*}
for a given constant $C$ independent of $h_E$ and dependent only on the mesh shape parameter $\rho$, where $\epsilon_{j,E}^{\BNAVEM}$ denotes the test loss
\begin{equation}
 \epsilon_{j,E}^{\BNAVEM} \coloneq \norm[\leb{2}{E}]{\Delta \bvarphi_{j,E} - g_{j,E}^{\Delta}}.
    \label{eq:testloss}
\end{equation}
\end{theorem}
\begin{proof}
    Since $\varphi_{j,E}$ and $\varphi_{j,E}^{\BNAVEM}$ coincide on $\partial E$, the difference
    \begin{equation*}
        w_{j,E} \coloneq \varphi_{j,E}^{\BNAVEM} - \varphi_{j,E}^{\VEM} \in \sob[0]{1}{E},
    \end{equation*}
    satisfies the following problem:
    \begin{equation*}
    \begin{cases}
        -\Delta w_{j,E} = - (\Delta \varphi_{j,E}^{\BNAVEM} - g^{\Delta}_{j,E}) =: f^{\BNAVEM}_{j,E} & \text{in } E,\\
        w_{j,E} = 0 & \text{on } \partial E.
    \end{cases}
    \end{equation*}
    A standard integration-by-parts argument yields
    \begin{equation*}
        \int_E f^{\BNAVEM}_{j,E}\, w_{j,E} = \norm[\leb{2}{E}]{\nabla w_{j,E}}^2.
    \end{equation*}
    As a consequence, by using the first Poincaré inequality \cite{CanutoNochetto2024} with constant $C$, we obtain
    \begin{align*}
        \norm[\leb{2}{E}]{\nabla w_{j,E}}^2 = \int_E f^{\BNAVEM}_{j,E}\, w_{j,E} &\leq \norm[\leb{2}{E}]{f^{\BNAVEM}_{j,E}} \norm[\leb{2}{E}]{w_{j,E}}
         \leq C h_E  \norm[\leb{2}{E}]{f^{\BNAVEM}_{j,E}} \norm[\leb{2}{E}]{\nabla w_{j,E}}\\
        &\leq C h_E  \norm[\leb{2}{E}]{\Delta \varphi_{j,E}^{\BNAVEM} - g^{\Delta}_{j,E}} \norm[\leb{2}{E}]{\nabla w_{j,E}} 
        = C h_E  \epsilon_{j,E}^{\BNAVEM} \norm[\leb{2}{E}]{\nabla w_{j,E}}.
    \end{align*}
    This concludes the proof.
\end{proof}

Here, we report some properties of VEM Lagrangian basis functions that we will use later.

\begin{lemma}[Lemma 1 in \cite{CangianiManzini2016}]\label{lemma:vem-basis}
    Under Assumption \ref{ass:meshassumption}, let $E \in \Th$ and $\{\vvarphi_{j,E}\}_{j=1}^{\Ndof[E]}$ be the Lagrangian basis functions of the VEM conforming space $\Vh[E]{k}$ with respect to the degrees of freedom \ref{dof:vertex}--\ref{dof:internal}. It holds true
    \begin{equation*}
        C^{-1} h_E \leq \norm[\leb{2}{E}]{\vvarphi_{j,E}} \leq C h_E, \quad C^{-1} \leq \norm[\leb{2}{E}]{\nabla \vvarphi_{j,E}} \leq C,
    \end{equation*}
    for a positive constant $C$ that depends only on the polynomial degree $k$ and the mesh regularity parameter $\rho$.
\end{lemma}

\begin{lemma}\label{lemma:vem-laplacian}
    Under Assumption \ref{ass:meshassumption}, it holds true
    \begin{equation*}
        \norm[\leb{2}{E}]{g^{\Delta}_{j,E}} \leq C h_E^{-1}, \quad \forall E \in \Th,
    \end{equation*}
    for a positive constant $C$ that depends only on the polynomial degree $k$ and the mesh regularity parameter $\rho$.
\end{lemma}
\begin{proof}
    Let $E$ be star-shaped with respect to a disk $\mathcal{D}_E$ of radius $\rho h_E$; we denote by $\widehat{\mathcal{D}}_E$ the disk concentric with $\mathcal{D}_E$ with radius $h_E$. It is clear that $\mathcal{D}_E \subset E \subset \widehat{\mathcal{D}}_E$. Now, let $b_E$ a bubble function defined over $\widehat{\mathcal{D}}_E$, since $g^{\Delta}_{j,E} \in \Poly{k-2}{E}$, by using the equivalence of norms on finite-dimensional vector
spaces, integration by parts, and the polynomial inverse inequality, it holds true \cite{DaVeigaMascotto2022}
    \begin{align*}
        \norm[\leb{2}{E}]{g^{\Delta}_{j,E}}^2 &\leq \norm[\leb{2}{\widehat{\mathcal{D}}_E}]{g^{\Delta}_{j,E}}^2 \leq C  \int_{\widehat{\mathcal{D}}_E} b_E (g^{\Delta}_{j,E})^2 = - C  \int_{\widehat{\mathcal{D}}_E} \nabla (b_E g^{\Delta}_{j,E}) \cdot \nabla \vvarphi_{j,E} \\&\leq C \norm[\leb{2}{\widehat{\mathcal{D}}_E}]{\nabla \vvarphi_{j,E}} \norm[\leb{2}{\widehat{\mathcal{D}}_E}]{\nabla( b_E g^{\Delta}_{j,E})} \leq C h_E^{-1} \norm[\leb{2}{\widehat{\mathcal{D}}_E}]{\nabla \vvarphi_{j,E}} \norm[\leb{2}{\widehat{\mathcal{D}}_E}]{b_E g^{\Delta}_{j,E}}
        \\&\leq C h_E^{-1} \norm[\leb{2}{\widehat{\mathcal{D}}_E}]{\nabla \vvarphi_{j,E}} \norm[\leb{2}{\widehat{\mathcal{D}}_E}]{g^{\Delta}_{j,E}}
        \\&\leq C h_E^{-1} \norm[\leb{2}{E}]{g^{\Delta}_{j,E}},
    \end{align*}
    where in the last step, we use Lemma~\ref{lemma:vem-basis} and results in \cite[Section 2.3]{Brenner2018}. This concludes the proof.
\end{proof}

\begin{lemma}\label{lemma:min-eigen}
Let $\mathbf{M}^{\VEM} \in \R^{\Ndof[E] \times \Ndof[E]}$ be the VEM mass matrix. The minimum eigenvalue of $\mathbf{M}^{\VEM}$ satisfies
\begin{equation*}
    \lambda_{\rm{min}}^{\VEM} \geq C h_E^2,
\end{equation*}
for a positive constant $C$ depending only on the shape-parameter $\rho$ but not on $h_E$.
\end{lemma}
\begin{proof}
     Since $\mathbf{M}^{\VEM}$ is symmetric, by the Courant-Fischer Formula, its minimum eigenvalue is defined as
\begin{equation*}
    \lambda_{\rm{min}}^{\VEM} = \min_{\vv \neq \bm{0}} \frac{\vv^T \mathbf{M}^{\VEM} \vv}{\vv^T \vv}.
\end{equation*}
Let $v_h \in \Vh[E]{k}$ be a VEM function and let $\vv$ be the related vector of degrees of freedom. We have
\begin{equation*}
    \vv^T \mathbf{M}^{\VEM} \vv = \scal[E]{v_h}{v_h} = \norm[\leb{2}{E}]{v_h}^2 \geq C h_E^2 \norm[2]{\vv}^2,
\end{equation*}
where we use Theorem 4.5 in \cite{ChenHuang2018}. Thus, we conclude that
\begin{equation*}
    \lambda_{\rm{min}}^{\VEM} \geq C h_E^2.
\end{equation*}
\end{proof}

Given a polygon $E$, let us define
\begin{equation*}
     \epsilon^{\BNAVEM}_E \coloneq \max_{j=1,\dots,\Ndof[E]} \epsilon^{\BNAVEM}_{j,E}.
\end{equation*}
The following result holds true.
\begin{proposition}\label{prop:networkaccuracy}
 Under Assumptions~\ref{ass:meshassumption} and ~\ref{ass:navem-regularity}, there exists a positive constant $C^{\BNAVEM}$ independent of $h$ such that
\begin{equation}\label{eq:bound-h-eps}
         \max_{E \in \Th} h_E \epsilon^{\BNAVEM}_{E}  \leq C^{\BNAVEM}.
\end{equation}
\end{proposition}
\begin{proof}
    The result follows from the triangle inequality, using Lemma~\ref{lemma:vem-laplacian} and Assumption~\ref{ass:navem-regularity} combined with \eqref{eq:scaling-lapl}. 
\end{proof}

As a direct consequence of Theorem \ref{theor:bnavem:local_basis} and Proposition \ref{prop:networkaccuracy}, we can derive the following corollary.
\begin{corollary}\label{cor:diff-mass-matrix}
Let $\mathbf{M}^{\VEM} \in \R^{\Ndof[E] \times \Ndof[E]}$ and $\mathbf{M}^{\BNAVEM} \in \R^{\Ndof[E] \times \Ndof[E]}$ be the VEM and B-NAVEM mass matrices, resp.. 
Under Assumption \ref{ass:navem-regularity}, there exists a positive constant $C$ that depends only on the polynomial degree $k$ and the mesh regularity parameter $\rho$ such that, we have
\begin{equation*}
    \norm[\infty]{\mathbf{M}^{\BNAVEM} - \mathbf{M}^{\VEM}} \le C h_E^3 \epsilon^{\BNAVEM}_E.
\end{equation*}  
\end{corollary}
\begin{proof}
   By Theorem \ref{theor:bnavem:local_basis}, we deduce that there exists a function $\varepsilon_{j,E}^{\BNAVEM}: E \to \R$ such that
\begin{equation*}
    \bvarphi_{j,E} = \vvarphi_{j,E} + \varepsilon_{j,E}^{\BNAVEM} \quad \text{with } \varepsilon_{j,E}^{\BNAVEM}  \equiv 0 \text{ on } \partial E,
\end{equation*}
and 
\begin{equation*}
    \norm[\leb{2}{E}]{ \varepsilon_{j,E}^{\BNAVEM}} \leq h_E \norm[\leb{2}{E}]{\nabla \varepsilon_{j,E}^{\BNAVEM}} \leq h_E^2 \epsilon^{\BNAVEM}_{j,E}.
\end{equation*}

For all $i,j=1,\dots,\Ndof[E]$, we have
\begin{align*}
    \vert \mathbf{M}^{\VEM}_{ij} - \mathbf{M}^{\BNAVEM}_{ij} \vert &= \vert \int_E \vvarphi_{i,E} \vvarphi_{j,E} - \int_E \bvarphi_{i,E} \bvarphi_{j,E} \vert  = \vert \int_E  \vvarphi_{i,E}\vvarphi_{j,E} - \int_E (\vvarphi_{i,E} + \varepsilon_{i,E}^{\BNAVEM}) (\vvarphi_{j,E} + \varepsilon_{j,E}^{\BNAVEM}) \vert 
    \\ &= \vert  \int_E \vvarphi_{i,E} \varepsilon_{j,E}^{\BNAVEM} + \int_E \vvarphi_{j,E} \,\varepsilon_{i,E}^{\BNAVEM} + \int_E  \varepsilon_{i,E}^{\BNAVEM}\, \varepsilon_{j,E}^{\BNAVEM} \vert \\
    &\leq \norm[\leb{2}{E}]{\vvarphi_{i,E}} \norm[\leb{2}{E}]{\varepsilon_{j,E}^{\BNAVEM}} + \norm[\leb{2}{E}]{\vvarphi_{j;E}} \norm[\leb{2}{E}]{\varepsilon_{i,E}^{\BNAVEM}} + \norm[\leb{2}{E}]{\varepsilon_{i,E}^{\BNAVEM}} \norm[\leb{2}{E}]{\varepsilon_{j,E}^{\BNAVEM}}\\
    &\leq C (h_E^3 \epsilon^{\BNAVEM}_E +  h_E^4 (\epsilon^{\BNAVEM}_E)^2)
\end{align*}
by virtue of Lemma \ref{lemma:vem-basis} for a constant $C$ independent of $h_E$. 
Using \eqref{eq:bound-h-eps}, we obtain the result.
\end{proof}

We observe that, for $k = 1$, the linear independence of the basis functions follows directly from their boundary values. For $k \geq 2$, however, it is necessary to prove that the B-NAVEM basis functions $\{\bvarphi_{j,E}\}_{j=1}^{\Ndof[E]}$ are linearly independent. The previous corollary provides a useful tool for this purpose: if the B-NAVEM mass matrix $\mathbf{M}^{\BNAVEM}$ is non-singular, as we prove in the following theorem, then the linear independence of the set $\{\bvarphi_{j,E}\}_{j=1}^{\Ndof[E]}$ follows immediately.

\begin{theorem}\label{theor:linearlyind}
    If $h_E \epsilon_E^{\BNAVEM}$ is sufficiently small, the B-NAVEM mass matrix $\mathbf{M}^{\BNAVEM}$ is symmetric positive definite and there exists a positive constant $C$ independent of $h_E$ such that
    \begin{equation*}
        \lambda_{\rm{min}}^{\BNAVEM} \geq C h_E^2,
    \end{equation*}
where $\lambda_{\rm{min}}^{\BNAVEM}$ denotes the smallest eigenvalues of $\mathbf{M}^{\BNAVEM}$.
\end{theorem}
\begin{proof}
We first note that the B-NAVEM mass matrix $\mathbf{M}^{\BNAVEM}$ is symmetric positive semi-definite by construction; thus, we need only to prove its non-singularity. 
We observe that, since both mass matrices $\mathbf{M}^{\BNAVEM}$ and $\mathbf{M}^{\VEM}$ are symmetric, due to Weyl's inequality \cite{franklin2000matrix}, Corollary~\ref{cor:diff-mass-matrix}, and Lemma~\ref{lemma:min-eigen}, we have
\begin{align*}
    \lambda_{\rm{min}}^{\BNAVEM} &= \lambda_{\rm{min}}^{\VEM} + \lambda_{\rm{min}}^{\BNAVEM} - \lambda_{\rm{min}}^{\VEM} 
    \geq \lambda_{\rm{min}}^{\VEM} - \vert \lambda_{\rm{min}}^{\BNAVEM} - \lambda_{\rm{min}}^{\VEM} \vert 
    \geq \lambda_{\rm{min}}^{\VEM} -
    \norm[2]{\mathbf{M}^{\BNAVEM} - \mathbf{M}^{\VEM}}\\
    &\geq \lambda_{\rm{min}}^{\VEM} - C\norm[\infty]{\mathbf{M}^{\BNAVEM} - \mathbf{M}^{\VEM}}
    \geq C h_E^2 \left(1 - C h_E \epsilon_E^{\BNAVEM} \right). 
\end{align*}
By assuming $C h_E \epsilon_E^{\BNAVEM} \leq 1/2$, we conclude the proof.

\end{proof}

Finally, we prove the following stability bound.
\begin{proposition}\label{prop:function_regularity}
    Under Assumptions \ref{ass:meshassumption} and \ref{ass:navem-regularity}, the set of B-NAVEM basis function $\{\bvarphi_{j,E}\}_{j=1}^{\Ndof[E]} \subset \con{0}{E} \cap \sob{1}{E}$ satisfies the following stability bound 
    \begin{equation}
        \sum_{j=1}^{\Ndof[E]} \norm[\leb{2}{E}]{\nabla \varphi_{j,E}^{\BNAVEM}} \le C,
        \label{eq:function_regularity}
    \end{equation}
    holds true for a positive constant $C$ that depends only on the mesh shape parameter $\rho$ and the polynomial degree $k$.
\end{proposition}
\begin{proof}
    We observe that, by using the triangle inequality,
    \begin{align*}
        \norm[\leb{2}{E}]{\nabla \varphi_{j,E}^{\BNAVEM}} &\leq \norm[\leb{2}{E}]{\nabla \varphi_{j,E}^{\BNAVEM} - \nabla \varphi_{j,E}^{\VEM}} + \norm[\leb{2}{E}]{\nabla \varphi_{j,E}^{\VEM}}
        \leq C(h_E\epsilon_{j,E}^{\BNAVEM} + 1) \leq C,
    \end{align*}
    due to Theorem \ref{theor:bnavem:local_basis} and Lemma \ref{lemma:vem-basis}. This completes the proof.
\end{proof}

\subsubsection{The B-NAVEM discrete problem}
Let us introduce the global space 
\begin{equation*}
    \bVh{k} = \{v \in \V: v_{|E} \in \bVh[E]{k}\},
\end{equation*}
where $\bVh[E]{k}$ is defined in \eqref{eq:localspacebnavem}.

Then, the B-NAVEM discrete counterpart of Problem \eqref{eq:exact_var_pb} reads as: 
\textit{Find $u_h^{\BNAVEM}\in \bVh{k}$ such that}
\begin{equation}
    \dbilin{u_h^{\BNAVEM}}{v_h^{\BNAVEM}} = \rightlin{v_h^{\BNAVEM}} \quad \forall v_h^{\BNAVEM}\in \bVh{k}.
    \label{eq:discrete_var_pb}
\end{equation}
or, in algebraic form, as
\begin{equation}
    \mathbf{A}^{\BNAVEM} \uu^{\BNAVEM} = \bm{f}^{\BNAVEM},
    \label{eq:bnavem_problem}
\end{equation}
where
\begin{align*}
    &\mathbf{A}^{\BNAVEM} \in \R^{\Ndof \times \Ndof}: && \mathbf{A}^{\BNAVEM}_{ij} = \sum_{E \in \Th} \int_E \nabla \bvarphi_{j} \cdot \nabla \bvarphi_{i},\\
    &\bm{f}^{\BNAVEM} \in \R^{\Ndof }: && \bm{f}^{\BNAVEM}_{j} = \sum_{E \in \Th} \int_E f \bvarphi_{i},
\end{align*}
and $\uu^{\BNAVEM}$ collects the values $\{u_i^{\BNAVEM}\}_{i=1}^{\Ndof}$ that express the B-NAVEM discrete solution as $u^{\BNAVEM}_h = \sum_{i=1}^{\Ndof} u_i^{\BNAVEM} \bvarphi_i$, where $\bvarphi_i$ represent the global counterpart of the local B-NAVEM basis functions.

Let us introduce the global B-NAVEM test loss 
\begin{equation}
    \epsilon^{\BNAVEM}_h \coloneq \max_{E \in \Th} \epsilon_{E}^{\BNAVEM},
    \label{eq:globalloss}
\end{equation}
the following a priori error estimate can be proved for the B-NAVEM method.

\begin{theorem}\label{theorem:bnavem-estimate}
Let $u \in \sob{k+1}{\Omega}$ be the exact solution of Problem \eqref{eq:continuous_poisson} and $u_h^{\BNAVEM} \in \bVh{k}$ be the solution of the B-NAVEM discrete problem \eqref{eq:bnavem_problem}. The following a priori error estimate holds
\begin{equation*}
    \norm[\leb{2}{\Omega}]{\nabla u - \nabla u_h^{\BNAVEM}} \le C (h^k \seminorm[\sob{k+1}{\Omega}]{u} + h^{-1} \epsilon^{\BNAVEM}_h \norm[\leb{2}{\Omega}]{u}),
\end{equation*}
for a given constant $C$ independent of $h$.
\end{theorem}

\begin{proof}

Let us introduce the VEM and B-NAVEM interpolation operators $I^{\VEM}$ and $I^{\BNAVEM}$ defined as usual as:
\begin{equation*}
    I^{\VEM}u = \sum_{i=1}^{\Ndof} u_i  \,\vvarphi_{i}, \qquad I^{\BNAVEM}u = \sum_{i=1}^{\Ndof} u_i  \,\varphi_{i}^{\BNAVEM},
\end{equation*}
where $\uu = (u_i)_{i=1}^{\Ndof}$ is the vector that collects the application of the operators $\dof_i$ to the exact solution $u$, and $\vvarphi_{i}$ are the global VEM basis functions.

Thanks to the coercivity and continuity of the involved linear and bilinear forms, the error $\norm[\leb{2}{\Omega}]{\nabla u - \nabla u_h^{\BNAVEM}}$ can be initially bounded by C{\'e}a's lemma \cite{ciarlet2002finite}:
\begin{equation}
    \norm[\leb{2}{\Omega}]{\nabla u - \nabla u_h^{\BNAVEM}}  \leq C  \inf_{v_h^{\BNAVEM} \in \bVh{k}} \norm[\leb{2}{\Omega}]{\nabla u - \nabla v_h^{\BNAVEM}}. 
    \label{eq:cea_bnavem}
\end{equation}
Now, let us bound the best approximation error as \cite{LBe13}
\begin{align}
     \nonumber   \inf_{v_h^{\BNAVEM} \in \bVh{k}} \norm[\leb{2}{\Omega}]{\nabla u - \nabla v_h^{\BNAVEM}} &\leq \norm[\leb{2}{\Omega}]{\nabla u - \nabla I^{\BNAVEM} u}\leq \norm[\leb{2}{\Omega}]{\nabla u - \nabla I^{\VEM} u} + \norm[\leb{2}{\Omega}]{\nabla I^{\VEM}u - \nabla I^{\BNAVEM} u}\\
      \nonumber  & \le Ch^k \seminorm[\sob{k+1}{\Omega}]{u} +  \norm[\leb{2}{\Omega}]{\sum_{i = 1}^{\Ndof}u_i \left(\nabla \vvarphi_i -\nabla \varphi_i^{\BNAVEM}\right)}\\
      \nonumber  & \le Ch^{k}\seminorm[\sob{k+1}{\Omega}]{u} + \sum_{E \in \Th}  \sum_{j = 1}^{\Ndof[E]} \vert u_j \vert \norm[\leb{2}{E}]{\nabla \vvarphi_{j,E} - \nabla \varphi_{j,E}^{\BNAVEM}}\\
       \label{eq:bnavem_proof_main_ineq} & \le Ch^k\seminorm[\sob{k+1}{\Omega}]{u} + C \sum_{E \in \Th} \sum_{j = 1}^{\Ndof[E]}  \vert u_j \vert h_E   \epsilon_{j,E}^{\BNAVEM},
\end{align}
where we have exploited Theorem \ref{theor:bnavem:local_basis}. The following bound holds true \cite{ChenHuang2018}
\begin{equation}
   \norm[2]{ \uu }  \leq C  h^{-1} \norm[\leb{2}{\Omega}]{I^{\VEM} u} \leq C h^{-1} (\norm[L^2(\Omega)]{u-I^V u} + \norm[L^2(\Omega)]{u}).
   \label{eq:normcoefficients}
\end{equation}
As a consequence, the last term in \eqref{eq:bnavem_proof_main_ineq} can be bounded as:
\begin{equation}
    \begin{aligned}
        \sum_{E \in \Th} \sum_{j = 1}^{\Ndof[E]}  \vert u_j \vert h_E   \epsilon_{j,E}^{\BNAVEM} \le C h \,  \epsilon^{\BNAVEM}_h \sqrt{\Ndof} \Vert \uu \Vert_2  \le C \epsilon^{\BNAVEM}_h \Vert \uu \Vert_2 \leq C \epsilon^{\BNAVEM}_h h^{-1} 
        (h^{k+1} |u|_{H^{k+1}(\Omega)} +\norm[L^2(\Omega)]{u}).
    \end{aligned}
    \label{eq:bnavem_sums_bound}
\end{equation}
Combining \eqref{eq:cea_bnavem}, \eqref{eq:bnavem_proof_main_ineq}, and \eqref{eq:bnavem_sums_bound} leads to the desired estimate.
\end{proof}

\subsection{Polynomial completeness}\label{sec:pol_reproducibility}

To guarantee polynomial reproducibility, in this section we describe
{\postprocess}, a post-processing algorithm that locally modifies the
basis functions $\{\bvarphi_{j,E}\}_{j=1}^{\Ndof[E]}$ of $\bVh[E]{k}$.
It yields a new set of basis functions
$\{\pvarphi_{j,E}\}_{j=1}^{\Ndof[E]}$ spanning the space
$\pVh[E]{k} \coloneq \myspan\{\pvarphi_{j,E}: j=1,\dots,\Ndof[E]\}$,
which contains $\Poly{k}{E}$.

Let us consider the matrix $\mathbf{D} \in \mathbb{R}^{n_k \times \Ndof[E]}$, whose entries are defined by
\begin{equation}
D_{\alpha i} = \dof_i(m_{\alpha}), \quad \forall \alpha=1,\dots,n_k \quad \forall i =1,\dots,\Ndof[E],
\label{eq:d_matrix}
\end{equation}
where $\{m_{\alpha}\}_{\alpha=1}^{n_k}$ denotes the set of standard scaled monomials over $E$, defined in \eqref{eq:monomial}.
Since $\Poly{k}{E} \subseteq \Vh[E]{k}$, the matrix $\mathbf{D}$ has full rank $n_k$. We assume that the columns of $\mathbf{D}$, and thus the degrees of freedom and the corresponding basis functions are ordered as established by a QR factorization with column pivoting \cite{Engler1997}, rather than by the natural ordering of the degrees of freedom introduced in Section \ref{sec:vem}. In this way, $\mathbf{D}$ is characterized by a non-singular square block $\mathbf{D}(:, 1:n_k)$, where $\mathbf{D}(:, 1:n_k)$ denotes the sub-matrix of $\mathbf{D}$ formed by all its rows and its first $n_k$ columns.

We observe that $\Poly{k}{E} \subset \pVh[E]{k}$ if and only if 
\begin{equation}
    m_{\alpha}(\xx) = \sum_{j =1}^{\Ndof[E]} D_{\alpha j} \pvarphi_{j,E}(\xx) \quad \forall \alpha=1,\dots,n_k \quad \text{and}\quad \forall \xx \in E.
    \label{eq:polynomial_reproducibility}
\end{equation}

By starting from $\{\bvarphi_{j,E}(\xx): j =1,\dots, \Ndof[E]\}$, our goal is to define the functions $\{\pvarphi_{j,E}(\xx): j=1,\dots,\Ndof[E]\}$ such that the condition \eqref{eq:polynomial_reproducibility} holds true for each $\xx \in E$.

Let us define 
\begin{equation*}
    \pvarphi_{j,E}(\xx) \coloneqq \bvarphi_{j,E}(\xx) \quad \forall j = n_{k} + 1,\dots,\Ndof[E];
\end{equation*}
then, system \eqref{eq:polynomial_reproducibility} can be rewritten as 
\begin{equation}\label{eq:postproc_initial_system}
    \begin{bmatrix}
        D_{11}  & D_{1 2} &D_{1 3} &\ldots &D_{1 n_k} \\
        D_{21}  & D_{2 2} &D_{2 3} &\ldots &D_{2 n_k}\\
        D_{31}  & D_{3 2} &D_{3 3} &\ldots &D_{3 n_k}\\
        \vdots & \vdots& \vdots &\ddots & \vdots\\
        D_{n_k 1}  & D_{n_k 2} &D_{n_k 3} &\ldots &D_{n_k n_k}\\
    \end{bmatrix}\begin{bmatrix}
        \pvarphi_{1,E}(\xx)\\
        \pvarphi_{2,E}(\xx)\\
        \pvarphi_{3,E}(\xx)\\
        \vdots\\
        \pvarphi_{n_k,E}(\xx)
    \end{bmatrix} = \begin{bmatrix}
         m_1(\xx) - \sum_{i=n_k+1}^{\Ndof[E]} D_{1 i} \bvarphi_{i,E}(\xx)\\
         m_2(\xx) - \sum_{i=n_k+1}^{\Ndof[E]} D_{2 i} \bvarphi_{i,E}(\xx)\\
         m_3(\xx) - \sum_{i=n_k+1}^{\Ndof[E]} D_{3 i} \bvarphi_{i,E}(\xx)\\
        \vdots \\
         m_{n_k}(\xx) - \sum_{i=n_k+1}^{\Ndof[E]} D_{n_k i} \bvarphi_{i,E}(\xx)\\
    \end{bmatrix},
\end{equation}
i.e. a system in the $n_k$ unknowns $\{\pvarphi_{j,E}(\xx)\}_{i=1}^{n_k}$ that can be solved by a standard LU factorization (with full pivoting). 

More precisely, assuming $D_{11} \neq 0$, we can remove $\pvarphi_{1,E}(\xx)$ from the last $n_{k}-1$ equations, obtaining:
\begin{equation*}
    \begin{bmatrix}
        D_{11}  & D_{1 2} &D_{1 3} &\ldots &D_{1 n_k} \\
        0  & D_{2 2}^{(2)} &D_{2 3}^{(2)} &\ldots &D_{2 n_k}^{(2)}\\
        0  & D_{3 2}^{(2)} &D_{3 3}^{(2)} &\ldots &D_{3 n_k}^{(2)}\\
        \vdots & \vdots& \vdots &\ddots & \vdots\\
        0  & D_{n_k 2}^{(2)} &D_{n_k 3}^{(2)} &\ldots &D_{n_k n_k}^{(2)}\\
    \end{bmatrix}\begin{bmatrix}
        \pvarphi_{1,E}(\xx)\\
        \pvarphi_{2,E}(\xx)\\
        \pvarphi_{3,E}(\xx)\\
        \vdots\\
        \pvarphi_{n_k,E}(\xx)
    \end{bmatrix} = \begin{bmatrix}
         m_1(\xx) - \sum_{i=n_k+1}^{\Ndof[E]} D_{1 i}\bvarphi_{i,E}\\
         m_2^{(2)}(\xx) - \sum_{i=n_k+1}^{\Ndof[E]} D_{2 i}^{(2)} \bvarphi_{i,E}\\
         m_3^{(2)}(\xx) - \sum_{i=n_k+1}^{\Ndof[E]} D_{3 i}^{(2)}\bvarphi_{i,E}\\
        \vdots \\
         m_{n_k}^{(2)}(\xx) - \sum_{i=n_k+1}^{\Ndof[E]} D_{n_k i}^{(2)}\bvarphi_{i,E}\\
    \end{bmatrix},
\end{equation*}
where
\begin{align*}
    &D_{\alpha i}^{(2)} \coloneqq  D_{\alpha i} -\frac{D_{\alpha 1}}{D_{11}}D_{1 i}  &&\hspace{-30pt}\forall \alpha =2,\dots,n_k\quad \forall i=2,\dots, \Ndof[E],\\
    &m_{\alpha}^{(2)} \coloneqq m_{\alpha} -\frac{D_{\alpha 1}}{D_{11}} m_1&&\hspace{-30pt}\forall \alpha =2,\dots,n_k.
\end{align*}
We can immediately observe that
\begin{equation}
    \dof_1(m_{\alpha}^{(2)}) = 0 \quad \dof_{\gamma}(m_{\alpha}^{(2)}) = D_{\alpha \gamma}^{(2)} \quad \forall \gamma = 2,\dots,\Ndof[E]. 
\end{equation}
Now, let us assume that $D_{22}^{(2)} \neq 0$. We can remove $\pvarphi_{2,E}(\xx)$ from the last $n_{k}-2$ equations and obtain
\begin{equation*}
    \begin{bmatrix}
        D_{11}  & D_{1 2} &D_{1 3} &\ldots &D_{1 n_k} \\
        0  & D_{2 2}^{(2)} &D_{2 3}^{(2)} &\ldots &D_{2 n_k}^{(2)}\\
        0  & 0 &D_{3 3}^{(3)} &\ldots &D_{3 n_k}^{(3)}\\
        \vdots & \vdots& \vdots &\ddots & \vdots\\
        0  & 0 &D_{n_k 3}^{(3)} &\ldots &D_{n_k n_k}^{(3)}\\
    \end{bmatrix}\begin{bmatrix}
        \pvarphi_{1,E}(\xx)\\
        \pvarphi_{2,E}(\xx)\\
        \pvarphi_{3,E}(\xx)\\
        \vdots\\
        \pvarphi_{n_k,E}(\xx)
    \end{bmatrix} = \begin{bmatrix}
         m_1(\xx) - \sum_{i=n_k+1}^{\Ndof[E]} D_{1 i}\bvarphi_{i,E}(\xx)\\
         m_2^{(2)}(\xx) - \sum_{i=n_k+1}^{\Ndof[E]} D_{2 i}^{(2)}  \bvarphi_{i,E}(\xx)\\
         m_3^{(3)}(\xx) - \sum_{i=n_k+1}^{\Ndof[E]} D_{3 i}^{(3)} \bvarphi_{i,E}(\xx)\\
        \vdots \\
         m_{n_k}^{(3)}(\xx) - \sum_{i=n_k+1}^{\Ndof[E]} D_{n_k i}^{(3)} \bvarphi_{i,E}(\xx)\\
    \end{bmatrix},
\end{equation*}
where
\begin{align*}
    &D_{\alpha i}^{(3)} \coloneqq D_{\alpha i}^{(2)} -\frac{D_{\alpha 2}^{(2)}}{D_{22}^{(2)}}D_{2 i}^{(2)}&&\hspace{-30pt}\forall \alpha =3,\dots,n_k\quad \forall i=3,\dots, \Ndof[E],\\
    &m_{\alpha}^{(3)} \coloneqq m_{\alpha}^{(2)} -\frac{D_{\alpha 2}^{(2)}}{D_{22}^{(2)}} m_2^{(2)}&&\hspace{-30pt} \forall \alpha =3,\dots,n_k.
\end{align*}
Again, we highlight that it holds
\begin{equation*}
    \dof_i(m_{\alpha}^{(3)}) = 0 \quad \forall i =1,2\quad \dof_{\gamma}(m_{\alpha}^{(3)}) = D_{\alpha \gamma}^{(3)} \quad \forall \gamma = 3,\dots,\Ndof[E]. 
\end{equation*}

By iterating the whole procedure, we obtain the following final upper triangular system
\begin{equation}
\label{eq:last_postproc_system}
    \begin{bmatrix}
        D_{11}  & D_{1 2} &D_{1 3} &\ldots &D_{1 n_k} \\
        0  & D_{2 2}^{(2)} &D_{2 3}^{(2)} &\ldots &D_{2 n_k}^{(2)}\\
        0  & 0 &D_{3 3}^{(3)} &\ldots &D_{3 n_k}^{(3)}\\
        \vdots & \vdots& \vdots &\ddots & \vdots\\
        0  & 0 &0 &\ldots &D_{n_k n_k}^{(n_k)}\\
    \end{bmatrix}\begin{bmatrix}
        \pvarphi_{1,E}(\xx)\\
        \pvarphi_{2,E}(\xx)\\
        \pvarphi_{3,E}(\xx)\\
        \vdots\\
        \pvarphi_{n_k, E}(\xx)
    \end{bmatrix} = \begin{bmatrix}
         m_1(\xx) - \sum_{i=n_k+1}^{\Ndof[E]} D_{1 i}\bvarphi_{i,E}(\xx)\\
         m_2^{(2)}(\xx) - \sum_{i=n_k+1}^{\Ndof[E]} D_{2 i}^{(2)} \bvarphi_{i,E}(\xx)\\
         m_3^{(3)}(\xx) - \sum_{i=n_k+1}^{\Ndof[E]} D_{3 i}^{(3)} \bvarphi_{i,E}(\xx)\\
        \vdots \\
         m_{n_k}^{(n_k)}(\xx) - \sum_{i=n_k+1}^{\Ndof[E]} D_{n_k i}^{(n_k)} \bvarphi_{i,E}(\xx)\\
    \end{bmatrix},
\end{equation}
where
\begin{align}
    \label{eq:def_Dmatrix_term}&D_{\alpha i}^{(\beta + 1)} \coloneqq D_{\alpha i}^{(\beta)} -\frac{D_{\alpha \beta}^{(\beta)}}{D_{\beta \beta}^{(\beta)}}D_{\beta i}^{(\beta)}&&\hspace{-30pt} \forall \alpha =\beta + 1,\dots,n_k\quad \forall i=\beta + 1,\dots, \Ndof[E],\\
\label{eq:def_monomial_term}&m_{\alpha}^{(\beta + 1)} \coloneqq m_{\alpha}^{(\beta)} -\frac{D_{\alpha \beta}^{(\beta)}}{D_{\beta \beta}^{(\beta)}} m_{\beta}^{(\beta)}&&\hspace{-30pt} \forall \alpha =\beta + 1,\dots,n_k,
\end{align}
which satisfy
\begin{equation}
    \dof_i(m_{\alpha}^{(\beta + 1)}) = 0 \quad \forall i =1,\dots,\beta, \qquad \dof_{\gamma}(m_{\alpha}^{(\beta + 1)}) = D_{\alpha \gamma}^{(\beta  + 1)} \quad \forall \gamma = \beta + 1,\dots,\Ndof[E]. 
    \label{dof:monomi_lu}
\end{equation}
We observe that if $\{\bvarphi_{j,E}\}_{j=1}^{\Ndof[E]}$ are in $\con{0}{E} \cap \sob{1}{E}$, then $\pvarphi_{j,E} \in \con{0}{E} \cap \sob{1}{E}$.

\begin{proposition}\label{prop:lagrange}
If the set of basis functions $\{\bvarphi_{j,E}\}_{j=1}^{\Ndof[E]}$ is a set of Lagrange basis functions w.r.t. the degrees of freedom \eqref{dof:vertex}--\eqref{dof:internal}, then also the set of basis functions $\{\pvarphi_{j,E}\}_{j=1}^{\Ndof[E]}$ is a set of Lagrange basis functions w.r.t. the same set of degrees of freedom, i.e. they satisfy
\begin{equation}
    \dof_i(\pvarphi_{j,E}) = \delta_{ij} \quad \forall i,j=1,\dots, \Ndof[E].
\end{equation}
\end{proposition}
\begin{proof}

Under the assumptions, by construction, the last $\Ndof[E] - n_k$ basis functions satisfy
\begin{equation*}
    \dof_i(\pvarphi_{j,E}) = \delta_{ij} \quad \forall i=1,\dots, \Ndof[E], \quad \forall j = n_k +1, \dots, \Ndof[E]
\end{equation*}

The first computed basis function $\varphi_{n_k}^P$ depends solely on the last $\Ndof[E] - n_k$ ones. Thus,
\begin{equation*}
    D_{n_k n_k}^{(n_k)} \pvarphi_{n_k,E} = m_{n_k}^{(n_k)} - \sum_{\beta = n_{k} + 1}^{\Ndof[E]} D_{n_k \beta}^{(n_k)} \pvarphi_{\beta,E}.
\end{equation*}
By applying the operator $\dof_{\gamma}$, for $\gamma=1,\dots,n_k$, to both terms of the previous equation, we obtain:
\begin{equation*}
     D_{n_k n_k}^{(n_k)} \dof_{\gamma} (\pvarphi_{n_k,E}) = \dof_{\gamma}(m_{n_k}^{(n_k)}) - \sum_{\beta = n_{k} + 1}^{\Ndof[E]} D_{n_k \beta}^{(n_k)} \dof_{\gamma}(\pvarphi_{\beta,E}).
\end{equation*}
From Equation \eqref{dof:monomi_lu}, we get
\begin{equation*}
    D_{n_k n_k}^{(n_k)} \dof_{\gamma} (\pvarphi_{n_k,E}) = D_{n_k \gamma}^{(n_k)} - \sum_{\beta = n_{k} + 1}^{\Ndof[E]} D_{n_k \beta}^{(n_k)} \dof_{\gamma}(\pvarphi_{\beta,E}).
\end{equation*}
Now, we must consider three different cases:
\begin{itemize}
    \item $\gamma = n_k$:
    \begin{equation*}
        D_{n_k n_k}^{(n_k)} \dof_{n_k} (\pvarphi_{n_k,E}) = D_{n_k n_k}^{(n_k)}\quad \Rightarrow \quad  \dof_{n_k} (\pvarphi_{n_k,E}) =1,
    \end{equation*}
    \item $\gamma < n_k$:
    \begin{equation*}
        D_{n_k n_k}^{(n_k)} \dof_{\gamma} (\pvarphi_{n_k,E}) = D_{n_k \gamma}^{(n_k)} = 0 \quad \Rightarrow \quad  \dof_{\gamma} (\pvarphi_{n_k,E}) = 0,
    \end{equation*}
    \item $\gamma > n_k$:
    \begin{equation*}
        D_{n_k n_k}^{(n_k)} \dof_{\gamma} (\pvarphi_{n_k,E}) =  D_{n_k \gamma}^{(n_k)} -  D_{n_k \gamma}^{(n_k)} = 0 \quad \Rightarrow \quad  \dof_{\gamma} (\pvarphi_{n_k,E}) = 0.
    \end{equation*}
\end{itemize}
Thus, the $n_k$-function is a Lagrange function w.r.t the degrees of freedom. By applying this procedure iteratively, we can easily prove that the first $n_{k}-1$-functions are Lagrange basis functions w.r.t. the set of degrees of freedom \eqref{dof:vertex}--\eqref{dof:internal}.
\end{proof}

\begin{remark}
    We observe that the set of B-NAVEM functions satisfies
\begin{equation}\label{eq:lagrangian_on_border}
    \dof_i(\bvarphi_{j,E}) = \delta_{ij} \quad \forall i =1,\dots,\Nv[E]k \quad \forall j=1,\dots,\Ndof[E],
\end{equation}
where we are assuming a natural ordering for DOFs. On the other hand, the extent to which the B-NAVEM functions satisfy
    \begin{equation*}
    \dof_i(\bvarphi_{j,E}) = \delta_{ij} \quad \forall i =\Nv[E]k + 1,\dots, \Ndof[E] \quad \forall j=1,\dots,\Ndof[E],
\end{equation*}
depends on the neural network accuracy. As a consequence, the set of basis functions $\{\pvarphi_{j,E}\}_{j=1}^{\Ndof[E]}$ satisfies \eqref{eq:lagrangian_on_border}
but may not be a set of Lagrange basis functions w.r.t. the set of internal DOFs \ref{dof:internal}. 
\end{remark}

\begin{proposition}\label{prop:continuity}
If the set of basis functions $\{\bvarphi_{j,E}\}_{j=1}^{\Ndof[E]}$ satisfies Equation~\eqref{eq:vem:linearity}, namely $\bvarphi_{j,E} \in \Poly{k}{e}$ for all $e \in \Eh[E]$, $\bvarphi_{j,E} \in \con{0}{\partial E}$, and
\begin{equation*}
    \dof_i(\bvarphi_{j,E}) = \delta_{ij} \quad \forall i =1,\dots,\Nv[E]k \quad \forall j=1,\dots,\Ndof[E],
\end{equation*}
where a natural ordering of the degrees of freedom is assumed, so that the corresponding basis functions are continuous on $\overline{\Omega}$, then the basis functions $\{\pvarphi_{j,E}\}_{j=1}^{\Ndof[E]}$ also satisfy Property~\eqref{eq:vem:linearity} and are continuous on $\overline{\Omega}$.
\end{proposition}
\begin{proof}
    The proof follows from the assumptions and Proposition \ref{prop:lagrange}.
\end{proof}

Next Proposition shows that the stability property proved in Proposition~\ref{prop:function_regularity} can be extended to the post-processed functions $\{\pvarphi_{j,E}\}_{j=1}^{\Ndof[E]}$. For this purpose, we first need to prove the following lemma.

\begin{lemma}\label{lemma:ratio-lu}
    Given an integer $\beta=1,\dots,n_k-1$, let us define $\pivot{1}{\alpha \beta} \coloneq \frac{D^{(\beta + 1)}_{\alpha \beta}}{D^{(\beta + 1)}_{\beta \beta}}$ and
    \begin{equation}
        \pivot{\gamma}{\alpha \beta} \coloneq \pivot{1}{\alpha,\beta + 1 - \gamma} - \sum_{j=1}^{\gamma - 1} \pivot{j}{\alpha \beta} \pivot{1}{\beta + 1 -j, \beta + 1- \gamma} \quad \forall \gamma = 2,\dots,\beta,
    \end{equation}
    for any $\alpha = \beta + 1,\dots,n_k$. It holds true 
    \begin{equation}
        1 + \sum_{\gamma=1}^{\beta} \vert \pivot{\gamma}{\alpha \beta} \vert \leq 2^{\beta}.
        \label{eq:lemma3:bound}
    \end{equation}
\end{lemma}
\begin{proof}
We proceed by recursion. First, we observe that $\vert \pivot{1}{\alpha \beta} \vert \leq 1$ due to the pivoting strategy, hence \eqref{eq:lemma3:bound} holds true for $\beta = 1$. Furthermore, we have, for each $\gamma = 2,\dots,\beta$,
\begin{equation}
    \vert \pivot{\gamma}{\alpha \beta} \vert \leq 1 + \sum_{j=1}^{\gamma - 1} \vert \pivot{j}{\alpha \beta} \vert.
    \label{eq:pivot-1}
\end{equation}
We split the sum up to $\gamma-1$ into the sum up to $\gamma-2$ and the $(\gamma-1)$-th term, and then bound the latter using \eqref{eq:pivot-1}. Thus, we have
\begin{align*}
    \vert \pivot{\gamma}{\alpha \beta} \vert &\leq 1 + \sum_{j=1}^{\gamma - 2} \vert \pivot{j}{\alpha \beta} \vert + \vert \pivot{\gamma - 1}{\alpha \beta} \vert\leq 2 \left(1 + \sum_{j=1}^{\gamma - 2} \vert \pivot{j}{\alpha \beta} \vert\right).
\end{align*}
Since $\gamma - 2 < \beta$, the recursive assumption leads to 
    \begin{align*}
        \vert \pivot{\gamma}{\alpha \beta} \vert &\leq 2^{\gamma - 1}.
    \end{align*}
    This concludes the proof by observing that
    \begin{equation*}
        1 + \sum_{\gamma=1}^{\beta} \vert \pivot{\gamma}{\alpha \beta} \vert \leq 1 + \sum_{\gamma=1}^{\beta} 2^{\gamma - 1} = 2^{\beta}.
    \end{equation*}
\end{proof}

\begin{proposition}
\label{prop:function_regularity_pnavem}
    The stability bound 
    \begin{equation}
    \label{eq:function_regularity_post}
        \sum_{j=1}^{\Ndof[E]} \norm[\leb{2}{E}]{\nabla \pvarphi_{j,E}} \le C
    \end{equation}
    holds true for a positive constant $C$ that depends only on the polynomial degree $k$ and on the mesh shape parameter $\rho$.
\end{proposition}
\begin{proof}
Since we perform the post-processing on the inertia element, we prove the proposition on an inertia element $E$, and then the result follows by a standard scaling argument.

The post-processed basis functions are the solution of the final upper triangular system~\eqref{eq:last_postproc_system}, i.e. their gradients can be written as
\begin{align}
    \label{eq:first-p-gradient}&\nabla \pvarphi_{n_k,E} = \frac{1}{{D^{(n_k)}_{n_k\, n_k}}}\left(\nabla m_{n_k}^{(n_k)} - \sum_{i=n_{k} + 1}^{\Ndof[E]} D^{(n_k)}_{n_k i} \nabla \bvarphi_{i,E}\right),\\
    \label{eq:other-p-gradient}&\nabla \pvarphi_{j,E} = \frac{1}{D_{jj}^{(j)}}\left(\nabla m_{j}^{(j)} - \sum_{i=n_{k} + 1}^{\Ndof[E]} D^{(j)}_{j i} \nabla \bvarphi_{i,E} - \sum_{i=j + 1}^{n_k} D^{(j)}_{j i} \nabla \pvarphi_{i,E} \right)\qquad \forall j = n_{k} - 1,\dots,1.
\end{align}

To bound these gradients, we first observe that,  given an integer $\beta=1,\dots,n_k - 1$, by using the quantities $\pivot{\gamma}{\alpha \beta}$, for any $\alpha = \beta + 1,\dots,n_k$ and $\gamma=1,\dots,\beta$ introduced in Lemma \ref{lemma:ratio-lu}, we can rewrite Equations~\eqref{eq:def_Dmatrix_term} and \eqref{eq:def_monomial_term} as
\begin{align*}
    &D_{\alpha i}^{(\beta + 1)} = D_{\alpha i} - \sum_{\gamma=1}^{\beta}\pivot{\gamma}{\alpha\beta}D_{\beta+1-\gamma,i},\\
    &m_{\alpha}^{(\beta + 1)} = m_{\alpha} - \sum_{\gamma=1}^{\beta}\pivot{\gamma}{\alpha\beta}m_{\beta+1-\gamma}.
\end{align*}

Since we are using scaled monomials, for each $\alpha = 1,\dots, n_k$ and $i =1,\dots,\Ndof[E]$, we have $\vert D_{\alpha i} \vert \leq 1$ and $\norm[\leb{2}{E}]{\nabla m_{\alpha}} \leq C$, for a constant $C$ depending on the polynomial degree. 
Therefore, we obtain
\begin{align}
     \label{eq:bound_on_generic_D_entry}&\vert D_{\alpha i}^{(\beta + 1)}\vert \le \vert D_{\alpha i} \vert + \sum_{\gamma=1}^{\beta} \vert \pivot{\gamma}{\alpha,\beta}\vert \vert D_{\beta+1-\gamma,i} \vert \le 1 + \sum_{\gamma=1}^{\beta} \vert \pivot{\gamma}{\alpha,\beta} \vert \leq 2^{\beta},\\
    \label{eq:bound_on_generic_mon_entry}&\norm[\leb{2}{E}]{\nabla m_{\alpha}^{(\beta + 1)}} \le C\left(1 + \sum_{\gamma=1}^{\beta} \vert \pivot{\gamma}{\alpha,\beta}\right) \leq C 2^{\beta},
\end{align}
as a consequence of Lemma \ref{lemma:ratio-lu}.

Second, we note that
\begin{equation*}
    \frac{1}{\vert D_{\beta \beta}^{(\beta)} \vert} \leq \norm[2]{\mathbf{U}^{-1}} \leq \norm[2]{\mathbf{P}_2^T} \norm[2]{\mathbf{A}^{-1}} \norm[2]{\mathbf{P}_1^T} \norm[2]{\mathbf{L}}\quad \forall \beta=1,\dots,n_k,
\end{equation*}
where $\mathbf{A} = \mathbf{D}(:, 1:n_k)$ is the coefficient matrix of the final system~\eqref{eq:postproc_initial_system} and $\mathbf{U}$, $\mathbf{L}$, $\mathbf{P}_1$ and $\mathbf{P}_2$ denotes the LU factors of $\mathbf{A}$, i.e. $\mathbf{P}_1 \mathbf{A} \mathbf{P_2} = \mathbf{L} \mathbf{U}$.

These three norms can be bounded as
\begin{equation*}
    \norm[2]{\mathbf{A}^{-1}} = \frac{1}{\sigma_{\rm{min}}(\mathbf{A})}, \qquad
\norm[2]{\mathbf{P}_1^T} = \norm[2]{\mathbf{P}_2^T} = 1, \qquad
\norm[2]{\mathbf{L}} \le \norm[F]{\mathbf{L}} \le \sqrt{\frac{n_k(n_k+1)}{2}} = C,
\end{equation*}
where $\sigma_{\rm{min}}(\mathbf{A})$ is the smallest singular value of $\mathbf{A}$, and $\norm[F]{\mathbf{L}}$ denotes the Frobenius norm of $\mathbf{L}$. The latter can be bounded by counting the non-zero elements in $\mathbf{L}$, which have always absolute values smaller than $1$. As a consequence, it holds true
\begin{equation}
\label{eq:diagonal_term_bound}
    \left|\frac{1}{D_{\beta,\beta}^{(\beta)}}\right| \le \frac{C}{\sigma_{\rm{min}}(\mathbf{A})}, \qquad \beta=1,\dots,n_k.
\end{equation}

Thus, the gradient in \eqref{eq:first-p-gradient} can be bounded as
\begin{align*}
    \norm[\leb{2}{E}]{\nabla\pvarphi_{n_k,E}} &\le \left|\frac{1}{D_{n_k,n_k}^{(n_k)}}\right| \left(\norm[\leb{2}{E}]{\nabla m_{n_k}^{(n_k)}} + \sum_{i=n_k+1}^{\Ndof[E]} \left|D_{n_k i}^{(n_k)}\right| \norm[\leb{2}{E}]{\nabla \bvarphi_{i,E}} \right) \\
    &\le \frac{C}{\sigma_{\rm{min}}(\mathbf{A})} \left( 2^{n_k-1} + \sum_{i=1}^{\Ndof} 2^{n_k-1}\right) \le  \frac{C}{\sigma_{\rm{min}}(\mathbf{A})} ,
\end{align*}
where we use Proposition \ref{prop:function_regularity} and the bounds in \eqref{eq:bound_on_generic_D_entry} and \eqref{eq:bound_on_generic_mon_entry}.
Such a bound can be used to bound $\norm[\leb{2}{E}]{\nabla\pvarphi_{n_k-1,E}}$ using the $(n_k-1)$-th equation as
\begin{align*}
    \norm[\leb{2}{E}]{\nabla\pvarphi_{n_k-1,E}} &\le \left|\frac{1}{D_{n_k-1,n_k-1}^{(n_k-1)}}\right| \left(\norm[\leb{2}{E}]{\nabla m_{n_k-1}^{(n_k-1)}} + \sum_{i=n_k+1}^{\Ndof[E]} \left|D_{n_k-1, i}^{(n_k-1)}\right| \norm[\leb{2}{E}]{\nabla \bvarphi_{i,E}} + \left|D_{n_k-1, n_k}^{(n_k)}\right|\norm[\leb{2}{E}]{\nabla\pvarphi_{n_k,E}}\right) \\
    &\le  \frac{C}{\sigma_{\rm{min}}(\mathbf{A})}  \left(2^{n_k - 2} + \sum_{i=1}^{\Ndof[E]}2^{n_k - 2} +  2^{n_k - 1} \frac{C}{\sigma_{\rm{min}}(\mathbf{A})}\right) \le  \frac{C}{\sigma_{\rm{min}}(\mathbf{A})^2}.
\end{align*}
This process can be repeated by recursion to bound all the terms in \eqref{eq:other-p-gradient} and prove the Proposition. 
\end{proof}


It is interesting to note that, since we define the ordering of degrees of freedom by using a QR-factorization with column pivoting, we are actually maximizing as much as possible the lower bound for the smallest singular value $\sigma_{\rm{min}}(\mathbf{A})$ of the sub-matrix $\mathbf{A} = \mathbf{D}(:, 1:n_k)$ that appears in the proof of the previous Proposition. Indeed, the following result of linear algebra is valid.
\begin{lemma}
    Let $\mathbf{M}$ be a $n \times m$ matrix with columns $\bm{C}_1,\dots,\bm{C}_m$. For each $j=1,\dots,m$, let $W_j$ be the subspace spanned by all columns except $C_j$. Then,
    \begin{equation*}
        \sigma_{\rm{min}}(\mathbf{M}) \geq \frac{1}{\sqrt{m}} \min_{j=1,\dots,m} {\rm{dist}}( \bm{C}_j, W_j),
    \end{equation*}
    where $\sigma_{\rm{min}}(\mathbf{M})$ denotes the smallest singular value of $\mathbf{M}$ and ${\rm{dist}}(\bm{C}_j, W_j)$ is the Euclidean distance from $\bm{C}_j$ to the subspace $W_j$.
\end{lemma}
\begin{proof}
    The distance between a column $\bm{C}_j$ and the subspace $W_j$ is defined as
    \begin{equation*}
        {\rm{dist}}(\bm{C}_j, W_j) = \min_{\yy \in W_j} \norm[2]{\bm{C}_j - \sum_{\substack{\alpha=1,\dots, m\\ \alpha \neq j}}y_{\alpha} \bm{C}_{\alpha}}.
    \end{equation*}
    For each $\xx \in \R^{m}$ such that $\norm[2]{\xx} = 1$, there exists $i = 1,\dots,m$ such that $\vert x_i \vert \geq \frac{1}{\sqrt{m}}$. Thus, we have
    \begin{equation*}
        \norm[2]{\mathbf{M}\xx} = \vert x_i \vert \norm[2]{ \bm{C}_i + \sum_{\substack{\alpha=1,\dots,m\\ \alpha \neq i}} \frac{x_{\alpha}}{x_i} \bm{C}_{\alpha}} \geq \frac{1}{\sqrt{m}}{\rm{dist}}(\bm{C}_i, W_i) \geq \frac{1}{\sqrt{m}} \min_{j=1\dots,m} {\rm{dist}}(\bm{C}_j, W_j). 
    \end{equation*}
    We conclude the proof just by noting that $\sigma_{\rm{min}}(\mathbf{M}) = \displaystyle{\min_{\xx \in \R^m:\ \norm[2]{\xx} = 1} \norm[2]{\mathbf{M}\xx}}$.
\end{proof}

Using a QR-factorization with column pivoting of $\mathbf{D}$, we are maximizing such a distance, thereby providing the best lower bound for $\sigma_{\rm{min}}(\mathbf{A})$.

Nonetheless, the stability bound \eqref{eq:function_regularity_post} depends on the condition number of the matrix $\mathbf{D}$. This matrix may become ill-conditioned for badly shaped elements and for higher orders of the method \cite{Mascotto2018, Teora2023}. For this reason, we always perform this post-processing phase on the inertia elements, whose corresponding matrix $\mathbf{D}$ is better conditioned with respect to the related matrix computed over the original elements, as proved in \cite{Teora2024}. The dependence of the stability bound \eqref{eq:function_regularity_post} with respect to the condition number of the matrix $\mathbf{D}$ will be further investigated in Section~\ref{sec:test2}.

Finally, we observe that the LU factorization can be replaced by the more stable
QR factorization with column pivoting, which preserves the same properties of the
underlying basis functions (see \ref{app:qr}). In this case, a single factorization
can be used to both determine the DOFs ordering and guarantee polynomial completeness,
whereas the LU-based approach requires two distinct factorizations. Nevertheless, we
focus on the LU factorization in this paper because it is, in our view, the most
natural way to introduce and present the algorithm.


\section{A priori error estimates}\label{sec:apriori}

For a generic function $w \in \con{0}{\overline{\Omega}} \cap \V$, let us introduce the following interpolator
\begin{equation*}
    I^{\NAPEM}w = \sum_{i=1}^{\Ndof} w_i  \,\varphi_{i}^{\NAPEM},
\end{equation*}
and let $\{w_i\}_{i=1}^{\Ndof}$ collect the values of the application of the operators $\dof_i$ to the function $w$. 

\begin{lemma}\label{lemma:pnavem_continuity_interp}
 There exists a constant $C$, independent of the mesh size $h$, such that, for any function $w\in \sob{k+1}{E}$,
\begin{equation}
    \norm[\leb{2}{E}]{\nabla I^{\NAPEM} w} \leq C ( h_E^k \seminorm[\sob{k+1}{E}]{w} + h^{-1}_E \norm[\leb{2}{E}]{w}), 
\end{equation}
holds true.
\end{lemma}
\begin{proof}
Given a generic function $w\in \sob{k+1}{E}$, by using Proposition \ref{prop:function_regularity_pnavem} and Remark~\ref{rem:inftynorm}, we have
\begin{align}
    \norm[\leb{2}{E}]{\nabla I^{\NAPEM} w}
    &= \norm[\leb{2}{E}]{\sum_{j=1}^{\Ndof[E]} w_j \, \nabla \varphi_{j,E}^{\NAPEM}}
    \le \left(\max_{j} \vert w_j \vert\right) \sum_{j=1}^{\Ndof[E]} \norm[\leb{2}{E}]{\nabla \varphi_{j,E}^{\NAPEM}} \le C \, \norm[\leb{\infty}{E}]{w}.
    \label{eq:continuita_interp}
\end{align}

Under Assumptions \ref{ass:meshassumption}, the element is star-shaped with respect to a ball $B$. Hereafter, we follow the same steps as in the proof of Lemma 4.3.4 in \cite{brennerscott}. In particular, we observe that
\begin{align*}
    \norm[\leb{\infty}{E}]{w} \leq \norm[\leb{\infty}{E}]{w - Q^{k + 1} w} + \norm[\leb{\infty}{E}]{Q^{k + 1} w}, 
\end{align*}
where $Q^{k + 1} w$ is the Taylor polynomial of degree $k + 1$ averaged over the ball $B$. 
Using Proposition 4.3.2 in \cite{brennerscott}, we have
\begin{align*}
\norm[\leb{\infty}{E}]{w - Q^{k+1} w} \leq C h_E^{k} \seminorm[\sob{k + 1}{E}]{w}.
\end{align*}
Moreover, following standard scaling arguments \cite{brennerscott} and the Remark 1.37 in \cite{ErnGuermond2022}, we obtain
\begin{align*}
    \norm[\leb{\infty}{E}]{Q^{k+1} w} \le C h_E^{-2} \norm[\leb{1}{E}]{w} \leq C h^{-1}_E \norm[\leb{2}{E}]{w}.
\end{align*}
Thus, we get
\begin{equation}
   \norm[\leb{\infty}{E}]{w} \leq C ( h_E^k \seminorm[\sob{k+1}{E}]{w} + h^{-1}_E \norm[\leb{2}{E}]{w}), 
\end{equation}
which concludes the proof.
\end{proof}

\begin{theorem}[$\sob{1}{}$-estimate]\label{theorem:h1_estimate}
Let $u \in \sob{k+1}{\Omega}$ be the exact solution of Problem \eqref{eq:continuous_poisson} and $u_h^{\NAPEM}$ be the solution of the discrete system \eqref{eq:discrete_var_pb_p} assembled with the new functions $\{\pvarphi_i\}_{i=1}^{\Ndof}$. The following a priori error estimate holds
\begin{equation*}
    \norm[\leb{2}{\Omega}]{\nabla u - \nabla u_h^{\NAPEM}} \le C h^k \seminorm[\sob{k+1}{\Omega}]{u},
\end{equation*}
for a constant $C$ independent of $h$.
\end{theorem}

\begin{proof}

By C{\'e}a's Lemma, we have
\begin{align}
     \norm[\leb{2}{\Omega}]{\nabla u - \nabla u_h^{\NAPEM}} &\leq C  \inf_{v_h^{\NAPEM} \in \pVh{k}} \norm[\leb{2}{\Omega}]{\nabla u - \nabla v_h^{\NAPEM}} 
    \le C\norm[\leb{2}{\Omega}]{\nabla u - \nabla I^{\NAPEM} u}.
        \label{eq:inizio_pnavem}
\end{align}

Let us introduce the $\leb{2}{}$-projection $\proj{E}{k}$ onto the polynomial space $\Poly{k}{E}$. By the triangle inequality, we have:
\begin{align}
    \norm[\leb{2}{E}]{\nabla u - \nabla I^{\NAPEM} u}
        &\leq 
        \norm[\leb{2}{E}]{\nabla u - \nabla \proj{E}{k} u} + \norm[\leb{2}{E}]{\nabla \proj{E}{k}u - \nabla I^{\NAPEM} \proj{E}{k} u}+ \norm[\leb{2}{E}]{\nabla I^{\NAPEM}( \proj{E}{k} u -  u)}.
        \label{eq:split_locale}
\end{align}
The first term can be bounded by using the Bramble-Hilbert Theorem \cite{brennerscott}:
\begin{equation}
    \norm[\leb{2}{E}]{\nabla u - \nabla \proj{E}{k} u} \le Ch^k\seminorm[\sob{k+1}{E}]{u}.
    \label{eq:bound_a1}
\end{equation}
Concerning the second term, it is zero since $I^{\NAPEM} q = q$ for each polynomial $q \in \Poly{k}{E}$, which follows from $\Poly{k}{E} \subseteq \pVh[E]{k}$. 

Moreover, by using Lemma \ref{lemma:pnavem_continuity_interp} applied to $w = \proj{E}{k} u -  u$ along with the Bramble-Hilbert Theorem, the third term in \eqref{eq:split_locale} can be bounded as
\begin{align}
    \norm[\leb{2}{E}]{\nabla I^{\NAPEM} (\proj{E}{k} u -  u)} &\leq C (h_E^k  \seminorm[\sob{k+1}{\Omega}]{u} + h_E^{-1} \norm[\leb{2}{E}]{\proj{E}{k} u -  u})\leq C h_E^{k} \seminorm[\sob{k+1}{E}]{  u}.
    \label{eq:bound_a3}
\end{align}

Combining the previous results, we get
\begin{align}
    \nonumber \norm[\leb{2}{E}]{\nabla u- \nabla I_h^{\NAPEM} u}
    &\leq
    Ch_E^k\seminorm[\sob{k+1}{E}]{u}.
\end{align}
Squaring and summing over the elements, and recalling \eqref{eq:inizio_pnavem}, yields the result.
\end{proof}

We observe that including polynomials inside the space helps us to remove the additional limiting factor that is present in the pure B-NAVEM method (see Theorem \ref{theorem:bnavem-estimate}).

Moreover, we observe that the Galerkin orthogonality condition holds true:
\begin{equation*}
    \dbilin{u - u_h^{\ast}}{v_h^{\ast}}  = 0 \quad \forall v_h^{\ast} \in V_{h,k}^{\ast}\quad \forall \ast \in \{\BNAVEM, \NAPEM\}.
\end{equation*}

\begin{theorem}[$\leb{2}{}$-estimate]
Let $u \in \sob{k+1}{\Omega} $ be the exact solution of Problem \eqref{eq:continuous_poisson} and $u_h^{\NAPEM}$ be the solution of the discrete system \eqref{eq:discrete_var_pb_p} assembled with the new functions $\{\pvarphi_i\}_{i=1}^{\Ndof}$. The following a priori error estimate holds true
\begin{equation*}
    \norm[\leb{2}{\Omega}]{ u -  u_h^{\NAPEM}} \le C h^{k+1} \seminorm[\sob{k+1}{\Omega}]{u},
\end{equation*}
for a constant $C$ independent of $h$.
\end{theorem}
\begin{proof}
   The result follows by using a standard duality argument and the Galerkin orthogonality property.
\end{proof}

\section{$\con{0}{}$-conforming basis functions elemental construction}\label{sec:boundary_funcs}

In this section, we show how to define the operator 
\begin{equation}
    {\mathcal B}_{j,E}:\con{0}{E}\rightarrow \con{0}{E},
    \label{eq:operator_B_def}
\end{equation}
that can be used to exactly enforce Property \eqref{eq:vem:linearity} when approximating $\varphi_{j,E}^V$ for $j=1,\dots,\Ndof[E]$ and $E \in \Th$. 

There are two different ways to build such an operator \cite{Sukumar2021, Pintore2023, Sukumar2026}. The former, introduced in \cite{Sukumar2026} and valid for convex elements, defines the operator $\mathcal{B}_{j,E}$ as 
\begin{equation}
    \left(\mathcal{B}_{j,E}({v})\right)({\xx}) := \left(v({\xx}) - \mathcal{L}(v(\xx))\right) + (\mathcal{L}(g_{j,E}^{\partial}))({\xx}),\qquad \forall {\xx} \in {E},\quad \forall 
    {v} \in \con{0}{E},
    \label{eq:bndry_operator_2}
\end{equation} 
where the lifting operator $\mathcal{L}$ is defined in \eqref{eq:ti_wachspress} by using Wachspress functions.
This construction ensures that the Laplacian of $\bvarphi_{j,E} = \mathcal{B}_{j,E}(\mathcal{N}_{j,E})$ is bounded and $\Delta \bvarphi_{j,E} \in \leb{2}{E}$. Thus, by inspecting the structure of \eqref{eq:ti_wachspress} and looking at $\mathcal{L}(\mathcal{N}_{j,E})$ and $\mathcal{L}(g_{j,E}^{\partial})$ as sum of products of functions belonging to $\sob{2}{E}$ with norms uniformly bounded, the Assumption~\ref{ass:navem-regularity} is satisfied.

The latter, introduced in \cite{Sukumar2021} and adopted in \cite{PintoreTeora2026} for defining the lowest-order version of B-NAVEM, defines the operator as 
\begin{equation}
    \left(\mathcal{B}_{j,E}({v})\right)({\xx}) := {\psi}^0_E({\xx})\, {v}({\xx}) + {\psi}_{j,E}({\xx}),\qquad \forall {\xx} \in {E},\quad \forall 
    {v} \in \con{0}{E}.
    \label{eq:bndry_operator}
\end{equation} 
Here, ${\psi}^0_E\in \con{0}{E}$ is a bubble function vanishing on $\partial {E}$, whereas ${\psi}_{j,E}\in \con{0}{E}$ is the transfinite interporlant of the boundary data. In this way, $\mathcal{B}_{j,E}({v})$ satisfies \eqref{eq:boundary_prop_boperator}, i.e. it coincides with ${\varphi}_{j,E}^{\VEM}$ on $\partial{E}$ for all $v \in \con{0}{E}$. In \cite{PintoreTeora2026}, both the definitions of the bubble function and the lifting operator were based on the concept of the \textit{Approximate Distance Function} (ADF) \cite{BiswasShapiro2004}, leading to a definition of the operator $\mathcal{B}_{j,E}$ such that the Laplacian of $\mathcal{B}_{j,E}({v})$ is generally unbounded at the vertices of a polygonal domain (it blows up as the reciprocal of the distance from vertices. See, e.g., the Appendix A in \cite{Pintore2023}). Hence, to ensure accurate model training, the collocation points needed to compute the discrete counterpart of the loss function \eqref{eq:bnavem:loss-jE} could not be chosen very close to polygon vertices. Therefore, this construction does not ensure that B-NAVEM functions satisfy Assumption~\ref{ass:navem-regularity}, but it can also be applied to generic concave elements.


In this Section, we thus introduce some ingredients that help us to build different families of basis functions, described in Section \ref{sec:numeri_experiements}, that do not contain polynomials, and that allow us to show the ability of the {\postprocess} algorithm in restoring polynomial completeness for the related spaces.

Let us denote by $\{e_1, \dots, e_{\Nv[E]}\}$ the set of edges of $E$, where $e_n$ connects the vertices $\vertex_n$ and $\vertex_{n+1}$ of $E$, with $\vertex_{\Nv[E] + 1} = \vertex_1$. We further define $\nn_n$ as the unit outward normal vector to the edge $e_n$ with respect to the polygon $E$.

\begin{figure}[!ht]
    \centering
    \begin{subfigure}{0.32\textwidth}
        \includegraphics[width=\linewidth]{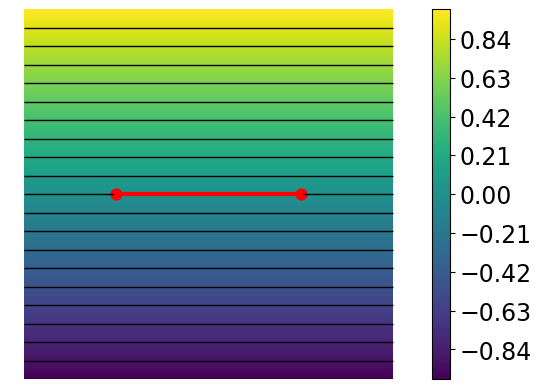}
        \caption{}
        \label{fig:bubble:distance}
    \end{subfigure}
    \begin{subfigure}{0.32\textwidth}
        \includegraphics[width=\linewidth]{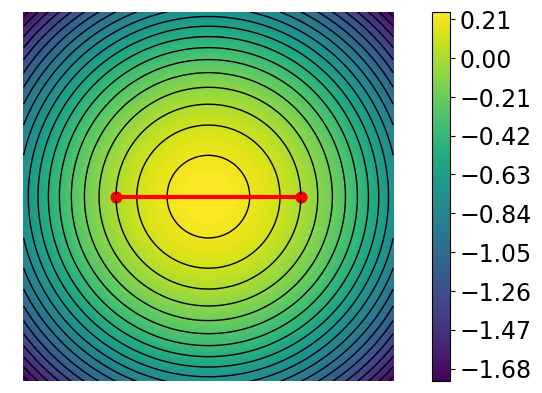}
        \caption{}
        \label{fig:bubble:trim}
    \end{subfigure}
    \begin{subfigure}{0.31\textwidth}
        \includegraphics[width=\linewidth]{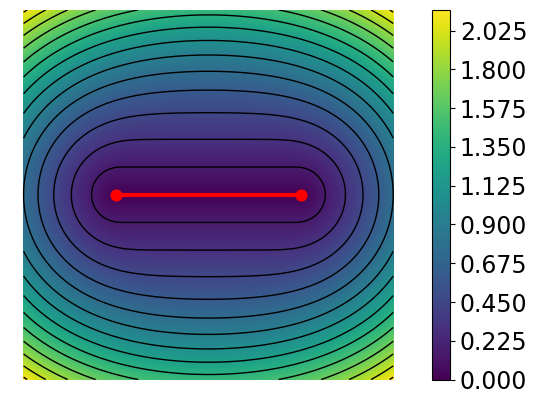}
        \caption{}
        \label{fig:bubble:adf_segment}
    \end{subfigure}
    \caption{Left: Distance function \eqref{eq:def_phi_line} to the line lying on the segment of extrema $(-0.5, 0.0)$ and $(0.5, 0.0)$. Center: trimming function defined in \eqref{eq:trim_circle} related the edge of extrema $(-0.5, 0.0)$ and $(0.5, 0.0)$. Right: the ADF defined in \eqref{eq:def_phi_segment} associated with the edge of extrema $(-0.5, 0.0)$ and $(0.5, 0.0)$.}
    \label{fig:bubble:edge}
\end{figure}

\subsection{A standard bubble function on generic polygons}\label{sec:bubble}

Our goal is to build a bubble function $\bubble$ that vanishes on $\partial E$, is strictly positive inside $E$, and has non-vanishing inward normal derivative on $\partial E$ \cite{RvachevSheiko1995RFunctionsIB}.

For this purpose, we first aim to build a function $\bubble[n]$, for each $n =1,\dots, \Nv[E]$, that is computable and strictly positive on $E$ and vanishes only on the edge $e_n$.

\begin{figure}[!ht]
    \centering
    \begin{subfigure}{0.32\textwidth}
        \includegraphics[width=\linewidth]{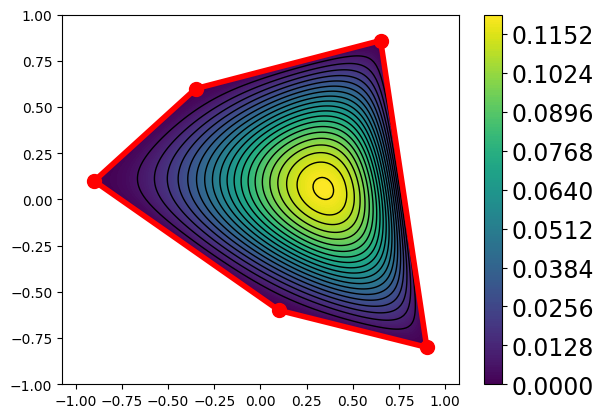}
        \caption{}
        \label{fig:bubble:convex}
    \end{subfigure}
    \begin{subfigure}{0.32\textwidth}
        \includegraphics[width=\linewidth]{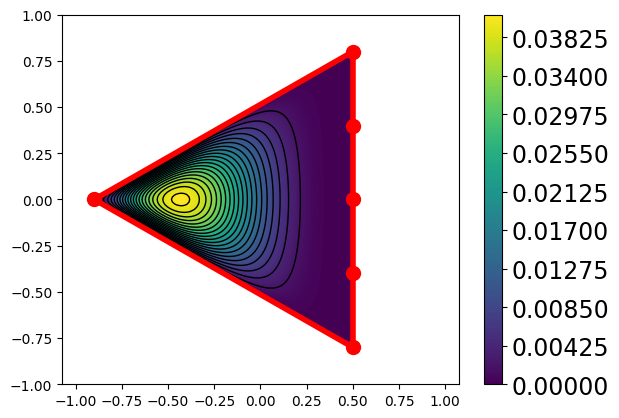}
        \caption{}
        \label{fig:bubble:hanging}
    \end{subfigure}
    \begin{subfigure}{0.32\textwidth}
        \includegraphics[width=\linewidth]{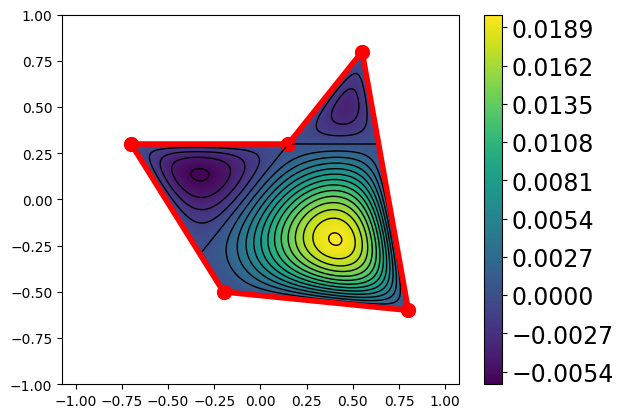}
        \caption{}
        \label{fig:bubble:concave}
    \end{subfigure}
    \begin{subfigure}{0.32\textwidth}
        \includegraphics[width=\linewidth]{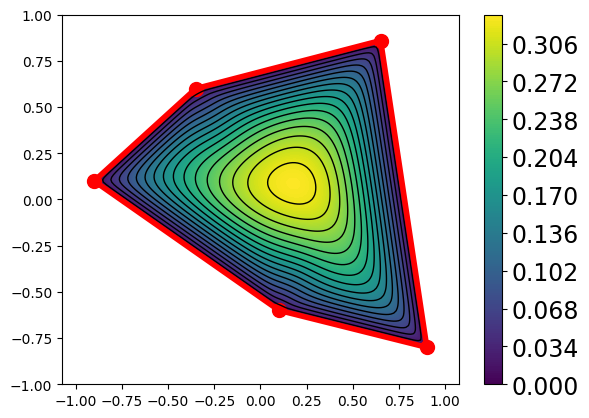}
        \caption{}
        \label{fig:adf:convex}
    \end{subfigure}
    \begin{subfigure}{0.32\textwidth}
        \includegraphics[width=\linewidth]{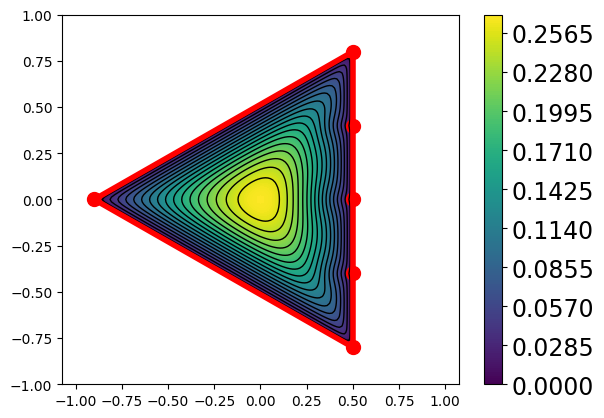}
        \caption{}
        \label{fig:adf:hanging}
    \end{subfigure}
    \begin{subfigure}{0.32\textwidth}
        \includegraphics[width=\linewidth]{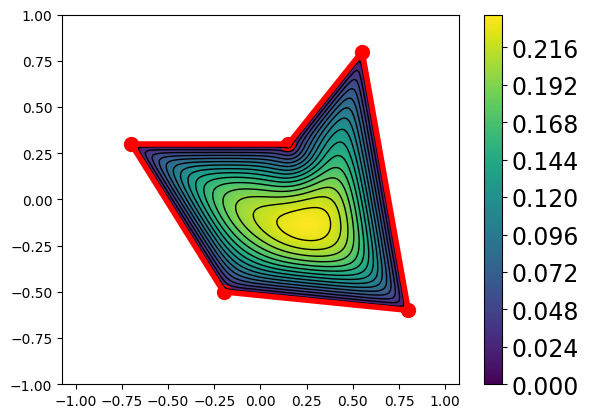}
        \caption{}
        \label{fig:adf:concave}
    \end{subfigure}
    \caption{Different alternatives to define the bubble function $\psi_{E}^0$ for a convex polygon (left), a triangle with hanging nodes (center), and a concave polygon (right). The top row refers to definition \eqref{eq:bubble_convex}, whereas the bottom row refers to definition \eqref{eq:def_phi_segment}.}
    \label{fig:bubble}
\end{figure}

Different strategies, that vary according to the shape of the polygon $E$, can be pursued to build these functions. If the element $E$ is convex, we can define the function $\bubble[n]$ as the \textit{signed distance function $d_n$} from the line $l_n$ where the edge $e_n$ lies, i.e. 
\begin{equation}
    d_n(\xx) \coloneq \left(\xx - \vertex_n\right)\cdot \nn_n \quad \forall \xx \in E \quad \forall n =1,\dots,\Nv[E].
    \label{eq:def_phi_line}
\end{equation}
which is strictly positive in $E$, and is zero on the whole $l_i$ and, in particular, on $e_i$. For an illustration, see Figure \ref{fig:bubble:distance}. Given these distance functions, we can finally compute
\begin{equation}
    \bubble(\xx) = \prod_{n=1}^{N_E^v} \bubble[n](\xx), \quad \text{with}\quad \bubble[n] \equiv d_n,
    \label{eq:bubble_convex}
\end{equation}
which vanishes on the boundary of $\partial E$ and is strictly positive inside $E$ (see Figure \ref{fig:bubble:convex}). 

This construction ensures that all the derivatives of $\bubble$ exist and are bounded in $E$, although $\bubble$ may be very small in regions close to many edges of $E$ and has zero normal derivative. See, for instance, Figure \ref{fig:bubble:hanging}. Moreover, if the element $E$ is concave, the line $l_n$ where the signed distance function $d_n$ is zero could cross the polygon, making $\bubble$ no longer strictly positive inside $E$ as requested (see Figure \ref{fig:bubble:concave}). 

\subsection{The approximate distance function}

If the element $E$ is concave, an example of a function that vanishes only on $e_i$ is represented by the \textit{Approximate Distance Function} (ADF) \cite{Pintore2023}. More specifically, let us define the trimming function $t_n^{\rm{c}}$ as
\begin{equation}
    t_n^{\rm{c}}(\xx) = \frac {1}{\vert e_n\vert} \left[\left(\frac {\vert e_n \vert}{2}\right)^2 - \norm[2]{\xx - \frac{\vertex_n+\vertex_{n+1}}{2}}^2\right], \quad \forall \xx \in E.
    \label{eq:trim_circle}
\end{equation}
Note that $t_n^{\rm{c}}\ge 0$ defines a circle of center $\dfrac{\vertex_n+\vertex_{n+1}}{2}$ as shown in Figure \ref{fig:bubble:trim}. Now, we can define $\bubble[n]$ as the ADF related to $e_n$, i.e.
\begin{equation}
    \bubble[n] = \sqrt{d_n^2(\xx) + \left(\frac{\sqrt{(t_n^{\rm{c}}(\xx))^2+d_n^4(\xx)}-t_n^{\rm{c}}(\xx)}{2}\right)^2} \quad \forall \xx \in  E\quad \forall n =1,\dots,\Nv[E].
    \label{eq:def_phi_segment}
\end{equation}
In this way, we obtain a function that is zero only on the edge $e_n$ and is strictly positive elsewhere. For an illustration of such a function, see Figure \ref{fig:bubble:adf_segment}.

Finally, the bubble function $\bubble$ can be chosen as the ADF to $\partial E$, normalized up to order $2$, i.e. as
\begin{equation}
    \bubble(\xx) = \dfrac{1}{\sqrt{{\mathlarger{\mathlarger{\sum}}}_{n=1}^{\Nv[E]} \dfrac{1}{\left(\psi_{n,E}^0(\xx)\right)^{2}}} }\quad \text{with $\psi_{n,E}$ defined in \eqref{eq:def_phi_segment}.}
    \label{eq:adf_boundary}
\end{equation}
With this definition, $\bubble$ vanishes at the boundary, is strictly positive inside $E$, and has a unitary inward normal derivative at the boundary, independently of the shape of the polygon $E$. See Figures \ref{fig:adf:convex}, \ref{fig:adf:hanging}, and \ref{fig:adf:concave} for an illustration of the ADF for different kinds of polygons.
Moreover, it is $\con{2}{}$ in all the points away from the boundary $\partial E$, while its Laplacian blows up at the vertices $\vertex_n$ (i.e. at the points belonging to $l_n$ such that $t_n^{\rm{c}}(\xx)=0$). 

Since the combination \eqref{eq:def_phi_segment}-\eqref{eq:adf_boundary} is applicable to all classes of polygons and has been numerically observed to provide better accuracy, it is therefore employed regardless of the specific polygon under consideration, i.e. for all the convex and concave polygons, as well as for polygons characterized by hanging nodes.

\begin{remark}
    In the case of convex polygonal elements, one could employ the standard polynomial bubble function defined in \eqref{eq:bubble_convex}. However, since this bubble may attain very small values in large portions of the element, particularly at points $z$ that are close to many edges, the contribution of the neural network output $\mathcal{N}(\xx_0)$ is strongly damped. As a consequence, learning an effective interior correction becomes difficult, because the product $\bubble(\zz) \mathcal{N}(\xx_0)$ has only a marginal influence on the final neural approximation $\left({\mathcal B}_{j,E}(\mathcal{N}(\xx_0))\right)(\zz)$.

    The bubble function \eqref{eq:bubble_convex} cannot be used to learn internal basis functions because both its Laplacian and its derivatives vanish at the vertices of the polygon. Consequently, it cannot approximate functions whose Laplacian is nonzero at the vertices.
\end{remark}

\subsection{The transfinite interpolation}\label{sec:boundary_funz}

Let us consider a point $\xx \in E$ and let $\Pi_{e_n} \xx$ be the orthogonal projection of $\xx$ onto the line $l_n$ where the edge $e_n$ lies. Let us introduce $s_n : e_n \to \R$, where $s_n(\Pi_{e_n} \xx)$ represents the curvilinear coordinate of $\Pi_{e_n} \xx$ with respect to the origin $\vertex_n$ and the direction $(\vertex_{n+1}-\vertex_n)$, whose expression is
\begin{equation*}
    s_n(\Pi_{e_n} \xx) = \frac{(\Pi_{e_n} \xx - \vertex_n) \cdot (\vertex_{n+1} - \vertex_n)}{\norm[2]{\vertex_{n+1} - \vertex_n}^2} = \frac{(\xx - \vertex_n) \cdot (\vertex_{n+1} - \vertex_n)}{\norm[2]{\vertex_{n+1} - \vertex_n}^2},
\end{equation*}
where the last equality has been proved in \cite{PintoreTeora2026}.

For any $j=1,\dots,\Ndof[E]$ and for any $\xx\in  E$, we introduce the functions
\begin{equation}
    \ti[n](\xx) = g^{\partial}_{j,E}(s_n(\Pi_{e_n} \xx)) \quad \forall n =1,\dots,\Nv[E].
    \label{eq:g_transfinite_interpolant}
\end{equation}
The transfinite interpolation $\ti$ of the function $\vvarphi_{j,E}$ is defined as \cite{Sukumar2021}
\begin{equation}
    \ti(\xx) = \sum_{n=1}^{\Nv[E]} w_{n,E}(\xx) \ti[n](\xx) \quad \forall j =1,\dots,\Ndof[E],
    \label{eq:transfite_interpolant}
\end{equation}
where
\begin{equation*}
    w_{n,E}(\xx) = \dfrac{\prod_{\gamma=1; \gamma \ne n}^{\Nv[E]} \bubble[n](\xx)}{\sum_{r=1}^{\Nv[E]} \prod_{\gamma=1;\gamma\ne r}^{\Nv[E]} \bubble[\gamma](\xx)} \quad \forall n =1,\dots,\Nv[E],
\end{equation*}
and the functions $\bubble[n]$ are defined in \eqref{eq:def_phi_segment}.

\subsection{Wachspress coordinates on convex polygons}\label{eq:wachspress}

Given a convex element $E$, the Wachspress coordinates are defined as \cite{Wachspress2016}
\begin{equation}
    \lambda_n^E(\xx) = \frac{w^{\lambda,E}_n(\xx)}{\sum_{\gamma=1}^{\Nv[E]} w^{\lambda,E}_\gamma(\xx)} \quad \forall n =1,\dots,\Nv[E]\quad \forall \xx \in E,
    \label{eq:wachspress_func}
\end{equation}
where
\begin{equation}
    w^{\lambda,E}_n(\xx) = \det(\nn_{n-1}, \nn_{n}) \prod_{\gamma=1, \gamma \neq n}^{\Nv[E]} d_\gamma(\xx),
\end{equation}
$d_{\gamma}$ are defined in \eqref{eq:def_phi_line}, and $\det(\nn_{n-1}, \nn_{n}) \coloneq \norm[2]{\nn_{n-1} \times \nn_{n}}$ is the Euclidean norm of the cross product of $\nn_{n-1}$ and $\nn_{n}$.

Wachspress coordinates are rational functions that, for $k=1$, satisfy the property \eqref{eq:vem:linearity} and are able to reproduce polynomials of degree $k \leq 1$ \cite{LBe13}. Therefore, in the lowest-order case, these functions can be used directly as basis functions and do not need {\postprocess}. For higher orders, Wachspress coordinates can be used to define the lifting operator $\mathcal{L}$ used in \eqref{eq:bndry_operator_2} to enforce the boundary conditions exactly along the element boundary. Since they are smoother and possess a bounded Laplacian, the resulting neural network approximation is also smoother and exhibits a bounded Laplacian \cite{Sukumar2026}. Moreover, thanks to {\postprocess}, the resulting family of basis functions recovers polynomial completeness for all orders $k \geq 1$.

Finally, the lifting operator $\mathcal{L}$, when applied to a generic function $v \in \con{0}{E}$, can be defined as \cite{Sukumar2026}:
\begin{equation}
    (\mathcal{L}(v))(\xx) = \sum_{n=1}^{\Nv[E]} \lambda_n^E(\xx)\Big[\alpha_{n,j,E}(\lambda_{n+1}^E(\xx)) + \alpha_{n-1,j,E}(1-\lambda_{n-1}^E(\xx)) - \alpha_{n,j,E}(0)\Big] \quad \forall j=1,\dots,\Ndof[E],
    \label{eq:ti_wachspress}
\end{equation}
where $\alpha_{n,j,E}$ is a parametrization of $v$ over the edge $e_n$.

\section{Numerical results}\label{sec:numeri_experiements}

In this section, we provide two numerical experiments aimed at numerically confirming the theoretical results proven in Section \ref{sec:apriori}. Moreover, we investigate the influence of the numerical quadrature in the error estimates (see Remark \ref{rem:quadrature}).

In both experiments, we test our algorithm for $k=1,2,3$, by considering different families of basis functions:
\begin{enumerate}[label=\textbf{FT}.\arabic*]
    \item\label{function:C} A B-NAVEM family of functions built by the operator \eqref{eq:bndry_operator_2} that exploits Wachspress coordinates. The neural networks are trained on 100 elements using an architecture with $4$ inner layers, each with $40$ neurons, and with regularization with parameter $1e\text{-}3$.
    \item\label{function:B} A B-NAVEM family of functions built with the operator \eqref{eq:bndry_operator}, whose corresponding $\boperator$ operator is defined by using
    \begin{enumerate}
        \item For $k=1$, the transfinite interpolant \eqref{eq:transfite_interpolant} and the ADF function \eqref{eq:adf_boundary}. We predict these functions using the same neural networks employed in \cite{PintoreTeora2026}.
        \item For $k = 2$ and $k=3$, the Wachspress lifting \eqref{eq:ti_wachspress} and the ADF function \eqref{eq:adf_boundary}. In these cases, we train new neural networks for $k=2$ and $k=3$ by using the same characteristics of the neural networks described in \cite{PintoreTeora2026}. Nonetheless, by preserving the same architecture to predict more variable functions, the corresponding neural networks will be, of course, less accurate. 
    \end{enumerate}
    \item\label{function:A} A non-neural family of basis. For this purpose, we adopt the lifting operator \eqref{eq:ti_wachspress} as vertex and edge basis functions, since the edge parametrizations of these functions are well-known. Then, the set is completed by adding
    \begin{enumerate}
        \item for $k=2$, the standard bubble function \eqref{eq:bubble_convex} as the only internal basis function.
        \item For $k=3$, the standard bubble function multiplied by $1$, $x_1$, and $x_1^2$ as the three internal basis functions.
    \end{enumerate}
    For $k=1$, we remark that this set of basis functions generates a polynomial-complete discrete space also without applying the {\postprocess}.
\end{enumerate}

We highlight that the set \ref{function:A} does not require the introduction of any neural networks and does not belong to the B-NAVEM family of methods. Nonetheless, we add this set of basis functions to show the broader applicability of {\postprocess}.

In the following, each set of functions will be denoted by $\{\bvarphi_{j,E}\}_{j=1}^{\Ndof[E]}$. To assess the effectiveness of the post-processing algorithm {\postprocess}, we compare the results obtained using the original basis functions with those obtained using the post-processed basis functions $\{\pvarphi_{j,E}\}_{j=1}^{\Ndof[E]}$, resulting from the application of the procedure described in Section \ref{sec:pol_reproducibility} to $\{\bvarphi_{j,E}\}_{j=1}^{\Ndof[E]}$.

We observe that all these functions satisfy, for each $E \in \Th$,
\begin{equation*}
    \dof_i(\varphi_{j,E}^{\ast}) = \delta_{ij}\quad \forall j =1,\dots,\Ndof[E]\quad \forall i =1,\dots,k\Nv[E],
\end{equation*}
but, generally,
\begin{equation*}
    \dof_i(\varphi_{j,E}^{\ast}) \neq \delta_{ij}\quad \forall j =1,\dots,\Ndof[E]\quad \forall i =k\Nv[E] + 1,\dots,\Ndof[E],
\end{equation*}
for $\ast \in \{\BNAVEM, \NAPEM\}$, where we are assuming a natural order for the degrees of freedom.

\begin{figure}[!ht]
    \centering
    \begin{subfigure}{0.26\textwidth}
        \includegraphics[width=1\linewidth]{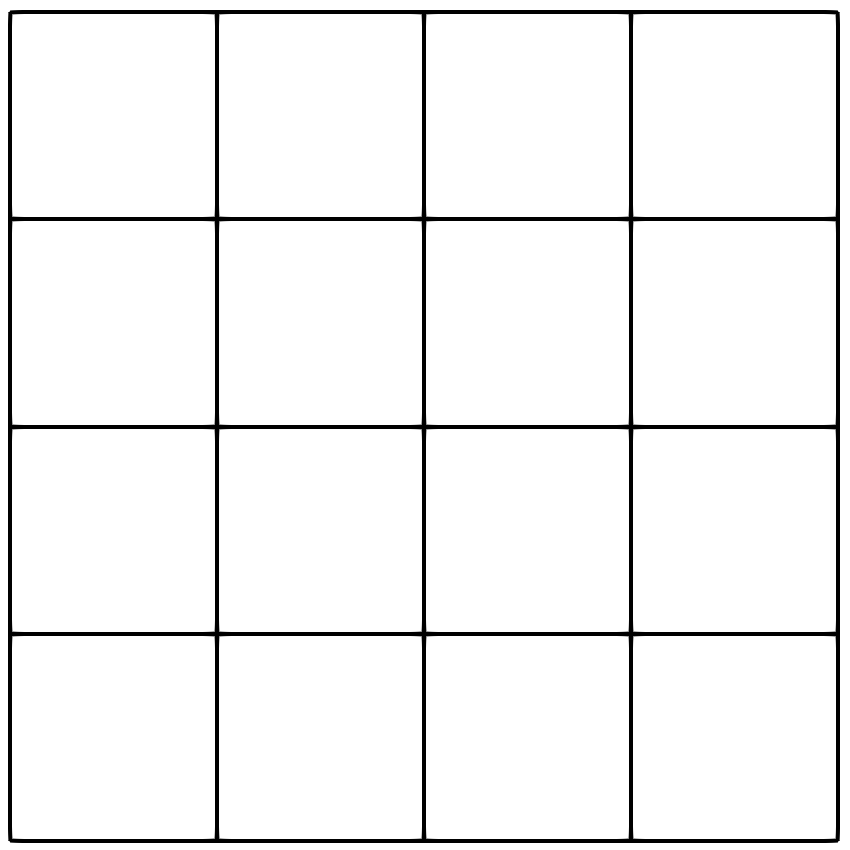}
        \caption{}
    \end{subfigure}
    \begin{subfigure}{0.26\textwidth}
        \includegraphics[width=1\linewidth]{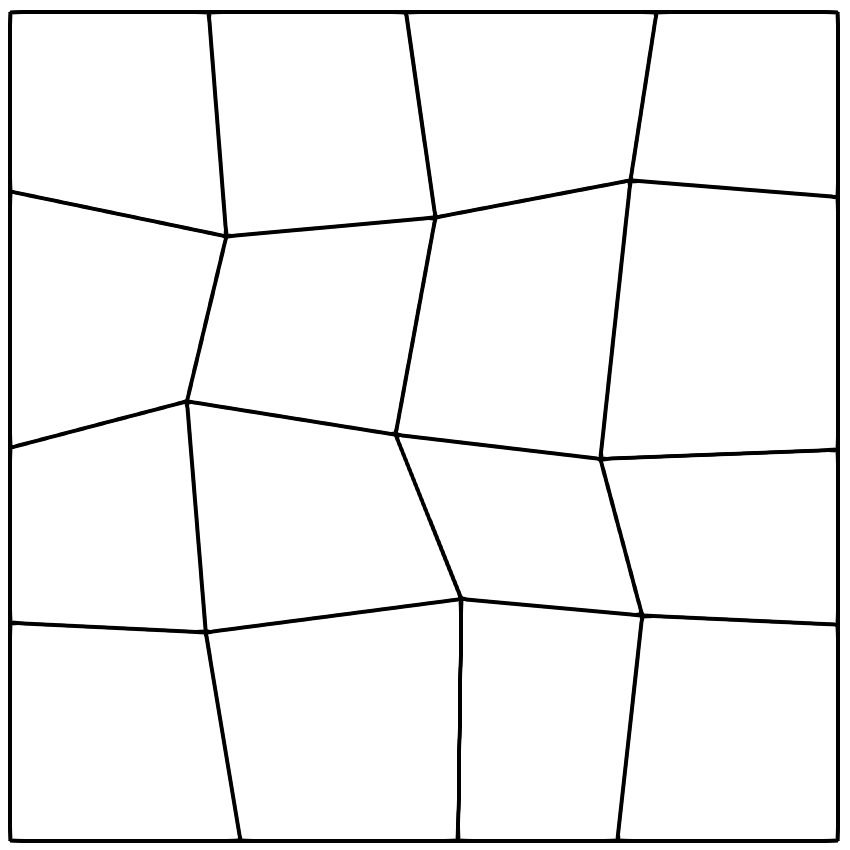}
        \caption{}
    \end{subfigure}
    \begin{subfigure}{0.26\textwidth}
        \includegraphics[width=1\linewidth]{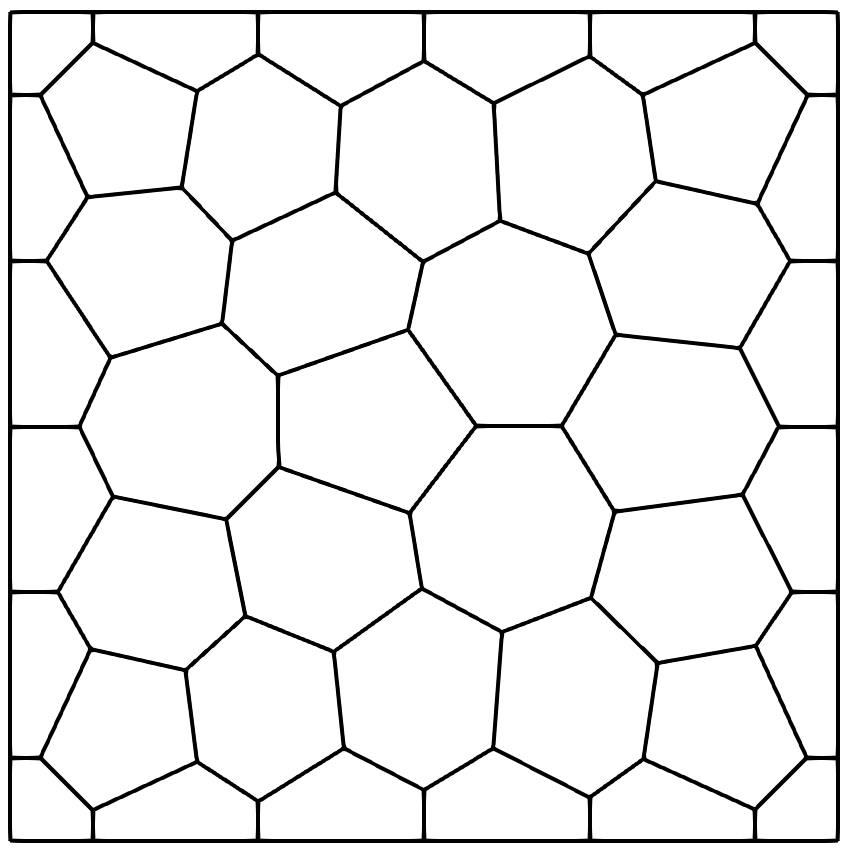}
        \caption{}
    \end{subfigure}
    \caption{Test 1, 2: Left: Mesh \ref{mesh:A}. Center: Mesh \ref{mesh:B}. Right: Mesh \ref{mesh:C}.}
    \label{fig:test1:mesh}
\end{figure}

To assess performance, we consider three different mesh types:
\begin{enumerate}[label=\textbf{MT}.\arabic*]
    \item\label{mesh:A} A Cartesian mesh.
    \item\label{mesh:B} A randomly distorted quadrilateral mesh.
    \item\label{mesh:C} A Voronoi mesh.
\end{enumerate}
Examples of these meshes are shown in Figure \ref{fig:test1:mesh} and are generated through the PolyDiM library \cite{Polydim}.

\begin{remark}
    We observe that all the tested cases concern uniquely convex elements, since the basis functions of types \ref{function:C}, \ref{function:A}, and \ref{function:B} for $k\geq 2$ we consider in this section are defined only over convex elements. Nonetheless, the algorithm we propose is general and can be applied independently of the shape of the underlying element. We further observe that, for functions of type \ref{function:B} and for $k\geq 2$, the transfinite interpolant \eqref{eq:transfite_interpolant} can also be employed to extend B-NAVEM to non-convex elements. Since our experiments are however restricted to convex elements, we adopt Wachspress lifting to build B-NAVEM functions of type \ref{function:B} and for $k\geq 2$. In this way, we exploit the smoothness of lifting to achieve higher accuracy at a lower computational effort.
\end{remark}

\subsection{Test 1: Polynomial patch test}

Let us consider the Poisson problem \eqref{eq:continuous_poisson}, modified by allowing non-zero Dirichlet boundary conditions. Let the source term $f$ and Dirichlet boundary condition be chosen in accordance with the following exact solution
\begin{equation*}
    u(x_1, x_2) = (x_1 + x_2 + 1)^k,
\end{equation*}
where $k$ represents the order of the method.

\begin{table}[!ht]
\centering
\caption{Test 1: Errors \eqref{eq:quadrature_errors} and \eqref{eq:errors} computed, for different values of the quadrature order $k_q$ and for $k=1,2,3$, both by using the original functions ($\BNAVEM$) and the post-processed functions ($\NAPEM$). Basis function type \ref{function:C}.}
\label{tab:test1:bndry_operator_2}
\resizebox{\textwidth}{!}{%
\begin{tabular}{lcccc|ccc|ccc|ccc|}
$k$ & Mesh Type    & \multicolumn{3}{c|}{$\ast = \BNAVEM$ ($k_q= 2k$)}           & \multicolumn{3}{c|}{$\ast = \BNAVEM$ ($k_q = 30$)}          & \multicolumn{3}{c|}{$\ast = \NAPEM$ ($k_q= 2k$)}            & \multicolumn{3}{c|}{$\ast = \NAPEM$ ($k_q = 30$)}           \\ \hline
    &              & ${\rm{err}}_{0}$ & ${\rm{err}}_{\nabla}$ & ${\rm{err}}_{q}$ & ${\rm{err}}_{0}$ & ${\rm{err}}_{\nabla}$ & ${\rm{err}}_{q}$ & ${\rm{err}}_{0}$ & ${\rm{err}}_{\nabla}$ & ${\rm{err}}_{q}$ & ${\rm{err}}_{0}$ & ${\rm{err}}_{\nabla}$ & ${\rm{err}}_{q}$ \\
1   & \ref{mesh:A} & 2.45e-06         & 7.37e-05              & 3.23e-06         & 1.82e-06         & 7.54e-05              & 1.61e-15         & 8.04e-08         & 1.04e-06              & 3.23e-06         & 3.00e-16         & 2.11e-15              & 1.55e-15         \\
2   & \ref{mesh:A} & 7.30e-03         & 3.07e-01              & 1.38e-02         & 5.73e-03         & 2.93e-01              & 1.70e-09         & 3.02e-05         & 1.06e-03              & 7.71e-03         & 3.69e-12         & 8.46e-11              & 2.40e-09         \\
3   & \ref{mesh:A} & 2.26e-02         & 7.65e-01              & 1.36e-02         & 2.21e-02         & 8.27e-01              & 6.29e-09         & 1.16e-04         & 4.68e-03              & 2.92e-03         & 9.46e-10         & 3.44e-08              & 1.33e-08         \\ \hline
1   & \ref{mesh:B} & 1.05e-03         & 2.36e-02              & 8.53e-04         & 1.06e-03         & 2.37e-02              & 1.62e-15         & 1.05e-05         & 1.44e-04              & 1.01e-03         & 4.88e-16         & 3.93e-15              & 1.88e-15         \\
2   & \ref{mesh:B} & 8.98e-03         & 3.39e-01              & 1.44e-02         & 7.70e-03         & 3.35e-01              & 1.93e-07         & 1.12e-04         & 2.65e-03              & 8.96e-03         & 6.92e-10         & 2.01e-08              & 1.02e-07         \\
3   & \ref{mesh:B} & 2.24e-02         & 8.11e-01              & 4.90e-02         & 2.20e-02         & 8.47e-01              & 1.40e-07         & 1.31e-04         & 5.56e-03              & 4.13e-03         & 2.43e-09         & 8.50e-08              & 2.58e-07         \\ \hline
1   & \ref{mesh:C} & 1.34e+00         & 3.46e+01              & 1.50e+00         & 1.41e+00         & 3.09e+01              & 1.08e-03         & 2.24e-03         & 7.79e-02              & 1.96e+00         & 1.86e-06         & 3.99e-05              & 1.02e-03         \\
2   & \ref{mesh:C} & 6.72e-02         & 2.42e+00              & 2.75e+00         & 6.94e-02         & 2.60e+00              & 2.91e-02         & 3.00e-04         & 1.51e-02              & 3.60e-02         & 1.59e-05         & 7.56e-04              & 1.62e-03         \\
3   & \ref{mesh:C} & 1.04e-01         & 6.48e+00              & 3.35e+00         & 1.79e-01         & 7.57e+00              & 6.66e-02         & 8.48e-04         & 8.24e-02              & 7.28e-02         & 5.01e-06         & 4.70e-04              & 4.48e-04         \\ \hline
\end{tabular}%
}
\end{table}

\begin{table}[!ht]
\centering
\caption{Test 1: Errors \eqref{eq:quadrature_errors} and \eqref{eq:errors} computed, for different values of the quadrature order $k_q$ and for $k=1,2,3$, both by using the original functions ($\BNAVEM$) and the post-processed functions ($\NAPEM$). Basis function type  \ref{function:B}.}
\label{tab:test1:bnavem}
\resizebox{\textwidth}{!}{%
\begin{tabular}{lcccc|ccc|ccc|ccc|}
$k$ & Mesh Type    & \multicolumn{3}{c|}{$\ast = \BNAVEM$ ($k_q= 2k$)}           & \multicolumn{3}{c|}{$\ast = \BNAVEM$ ($k_q = 30$)}          & \multicolumn{3}{c|}{$\ast = \NAPEM$ ($k_q= 2k$)}            & \multicolumn{3}{c|}{$\ast = \NAPEM$ ($k_q = 30$)}           \\ \hline
    &              & ${\rm{err}}_{0}$ & ${\rm{err}}_{\nabla}$ & ${\rm{err}}_{q}$ & ${\rm{err}}_{0}$ & ${\rm{err}}_{\nabla}$ & ${\rm{err}}_{q}$ & ${\rm{err}}_{0}$ & ${\rm{err}}_{\nabla}$ & ${\rm{err}}_{q}$ & ${\rm{err}}_{0}$ & ${\rm{err}}_{\nabla}$ & ${\rm{err}}_{q}$ \\
1   & \ref{mesh:A} & 6.95e-03         & 1.26e-01              & 2.29e-03         & 6.30e-03         & 1.37e-01              & 5.34e-06         & 6.91e-05         & 8.97e-04              & 2.29e-03         & 3.74e-08         & 4.88e-07              & 5.34e-06         \\
2   & \ref{mesh:A} & 3.17e-01         & 6.99e+00              & 1.82e-02         & 2.66e-01         & 6.33e+00              & 1.17e-05         & 4.76e-05         & 1.10e-03              & 1.65e-02         & 1.66e-08         & 4.91e-07              & 5.29e-06         \\
3   & \ref{mesh:A} & 5.25e-01         & 1.27e+01              & 1.76e-02         & 5.46e-01         & 1.32e+01              & 5.92e-05         & 2.85e-03         & 1.14e-01              & 3.16e-02         & 2.10e-06         & 9.37e-05              & 3.69e-05         \\ \hline
1   & \ref{mesh:B} & 4.63e-03         & 1.05e-01              & 4.14e-03         & 4.33e-03         & 1.15e-01              & 1.66e-05         & 7.18e-05         & 8.67e-04              & 4.74e-03         & 3.01e-07         & 3.10e-06              & 1.94e-05         \\
2   & \ref{mesh:B} & 4.07e-01         & 7.76e+00              & 2.01e-02         & 3.59e-01         & 7.15e+00              & 3.29e-05         & 2.95e-04         & 7.11e-03              & 1.73e-02         & 9.94e-08         & 2.64e-06              & 7.30e-06         \\
3   & \ref{mesh:B} & 6.39e-01         & 1.43e+01              & 3.71e-02         & 6.63e-01         & 1.49e+01              & 1.43e-04         & 2.50e-03         & 8.45e-02              & 2.59e-02         & 2.07e-06         & 9.19e-05              & 3.82e-05         \\ \hline
1   & \ref{mesh:C} & 2.86e-02         & 6.74e-01              & 7.78e-03         & 3.29e-02         & 6.59e-01              & 1.30e-04         & 7.95e-05         & 2.39e-03              & 9.45e-03         & 9.14e-07         & 2.59e-05              & 1.72e-04         \\
2   & \ref{mesh:C} & 1.06e-02         & 9.22e-01              & 3.24e-02         & 1.21e-02         & 1.08e+00              & 6.88e-04         & 1.52e-04         & 9.34e-03              & 3.13e-02         & 5.28e-06         & 2.46e-04              & 6.23e-04         \\
3   & \ref{mesh:C} & 9.30e-02         & 6.03e+00              & 9.26e-02         & 8.74e-02         & 6.11e+00              & 3.42e-03         & 7.05e-04         & 6.17e-02              & 7.08e-02         & 4.96e-05         & 4.31e-03              & 4.15e-03         \\ \hline
\end{tabular}%
}
\end{table}

\begin{table}[!ht]
\centering
\caption{Test 1: Errors \eqref{eq:quadrature_errors} and \eqref{eq:errors} computed, for different values of the quadrature order $k_q$ and for $k=1,2,3$, both by using the original functions ($\BNAVEM$) and the post-processed functions ($\NAPEM$). Basis function type \ref{function:A}.}
\label{tab:test1:wachspress}
\resizebox{\textwidth}{!}{%
\begin{tabular}{lcccc|ccc|ccc|ccc|}
$k$ & Mesh Type    & \multicolumn{3}{c|}{$\ast = \BNAVEM$ ($k_q= 2k$)}           & \multicolumn{3}{c|}{$\ast = \BNAVEM$ ($k_q = 30$)}          & \multicolumn{3}{c|}{$\ast = \NAPEM$ ($k_q= 2k$)}            & \multicolumn{3}{c|}{$\ast = \NAPEM$ ($k_q = 30$)}           \\ \hline
    &              & ${\rm{err}}_{0}$ & ${\rm{err}}_{\nabla}$ & ${\rm{err}}_{q}$ & ${\rm{err}}_{0}$ & ${\rm{err}}_{\nabla}$ & ${\rm{err}}_{q}$ & ${\rm{err}}_{0}$ & ${\rm{err}}_{\nabla}$ & ${\rm{err}}_{q}$ & ${\rm{err}}_{0}$ & ${\rm{err}}_{\nabla}$ & ${\rm{err}}_{q}$ \\
1   & \ref{mesh:A} & 5.28e-16         & 3.29e-15              & 4.84e-16         & 5.01e-16         & 3.07e-15              & 7.28e-16         & 2.57e-16         & 1.64e-15              & 8.81e-16         & 2.73e-16         & 1.71e-15              & 7.93e-16         \\
2   & \ref{mesh:A} & 5.94e-15         & 1.24e-13              & 1.57e-15         & 1.56e-14         & 2.71e-13              & 1.50e-14         & 1.14e-14         & 2.70e-13              & 1.75e-15         & 1.21e-14         & 7.98e-14              & 1.52e-14         \\
3   & \ref{mesh:A} & 1.02e-04         & 3.87e-03              & 2.47e-04         & 1.02e-04         & 3.88e-03              & 2.62e-14         & 6.99e-07         & 2.43e-05              & 2.72e-04         & 3.98e-14         & 1.06e-12              & 2.87e-14         \\ \hline
1   & \ref{mesh:B} & 6.18e-06         & 8.38e-05              & 6.98e-04         & 6.15e-16         & 6.54e-15              & 9.49e-16         & 6.18e-06         & 8.38e-05              & 6.98e-04         & 6.48e-16         & 5.86e-15              & 9.99e-16         \\
2   & \ref{mesh:B} & 5.37e-04         & 1.69e-02              & 1.55e-03         & 5.28e-04         & 1.72e-02              & 1.50e-14         & 8.98e-06         & 2.13e-04              & 1.47e-03         & 1.05e-14         & 2.42e-13              & 1.48e-14         \\
3   & \ref{mesh:B} & 2.08e-03         & 7.44e-02              & 3.00e-04         & 2.08e-03         & 7.44e-02              & 2.59e-14         & 3.44e-06         & 1.27e-04              & 3.11e-04         & 4.24e-14         & 1.43e-12              & 2.68e-14         \\ \hline
1   & \ref{mesh:C} & 1.18e-04         & 4.16e-03              & 1.41e-02         & 1.39e-15         & 1.31e-14              & 2.04e-15         & 1.18e-04         & 4.16e-03              & 1.41e-02         & 1.35e-15         & 1.19e-14              & 2.09e-15         \\
2   & \ref{mesh:C} & 5.48e-04         & 3.00e-02              & 7.54e-02         & 5.55e-04         & 2.95e-02              & 7.58e-14         & 2.18e-04         & 1.09e-02              & 7.50e-02         & 8.79e-15         & 1.84e-13              & 7.77e-14         \\
3   & \ref{mesh:C} & 2.47e-03         & 1.38e-01              & 7.92e-02         & 2.55e-03         & 1.41e-01              & 4.12e-12         & 2.49e-04         & 1.77e-02              & 7.94e-02         & 6.79e-14         & 3.63e-12              & 4.15e-12         \\ \hline
\end{tabular}%
}
\end{table}

To compute integrals over an element $E \in \Th$, we employ a Gauss quadrature rule of order $20$ along each edge of $E$. In the element interior, we consider two different quadrature strategies: a Gauss quadrature rule for triangles, applied over a sub-triangulation of $E$, with order $2k$, and another with order $30$. This comparison allows us to assess the influence of the internal quadrature rule on the convergence errors. Specifically, the value $30$ corresponds to the maximum degree of available quadrature order for the reduced Gaussian rule in the used library PolyDiM \cite{Polydim} and it is used here to reproduce ``exact'' quadrature.

To assess the error introduced by such internal quadrature formulas, we measure the following error
\begin{equation}
    {\rm{err}}_{q}^2 = \max_{E \in \Th}   \sum_{\alpha=1}^{n_k}\sum_{i=1}^{\Ndof[E]} \Big[\int_E \Delta m_{\alpha} \varphi^{\ast}_{i,E} + \int_E \nabla m_{\alpha} \cdot \nabla \varphi^{\ast}_{i,E} - \int_{\partial E} \nabla m_{\alpha} \cdot \nn\  \varphi^{\ast}_{i,E} \Big]^2 \quad \forall \ast \in \{\BNAVEM, \NAPEM\}.
    \label{eq:quadrature_errors}
\end{equation}

Moreover, given $u_h^{\ast}$ as the solution of the discrete system, we compute the following discrete errors:
\begin{equation}
    {\rm{err}}_{0}^2 = \sum_{E \in \Th} \norm[\leb{2}{E}]{ u_h^{\ast} - u }^2\quad {\rm{err}}_{\nabla}^2 = \sum_{E \in \Th} \norm[\leb{2}{E}]{ \nabla u_h^{\ast} - \nabla u }^2 \quad \forall \ast \in \{\BNAVEM, \NAPEM\}.
    \label{eq:errors}
\end{equation}

Tables \ref{tab:test1:bndry_operator_2}, \ref{tab:test1:bnavem}, and \ref{tab:test1:wachspress} report the method errors \eqref{eq:errors} and the quadrature error \eqref{eq:quadrature_errors} computed by employing basis functions of type \ref{function:C}, \ref{function:B}, and \ref{function:A}, resp., over meshes shown in Figure \ref{fig:test1:mesh}. From these tables, we can observe that the basis functions of type \ref{function:A} exactly reproduce polynomials for $k=1$, both with and without the application of the algorithm proposed in Section \ref{sec:pol_reproducibility}, since the corresponding approximation space contains all linear polynomials in both cases. Indeed, the errors \eqref{eq:quadrature_errors} are zero up to the precision imposed by the quadrature error \eqref{eq:quadrature_errors}. For $k = 2, 3$, functions of type \ref{function:C} reproduce only linear polynomials: RePAlg is therefore essential to recover all polynomials up to degree k, and thus the exact solution, on every mesh element. Concerning the basis functions of types \ref{function:C} and \ref{function:B}, the proposed algorithm consistently improves the accuracy of the B-NAVEM method until the error reaches a lower bound determined by the quadrature accuracy. We further remark that basis functions of type \ref{function:A} are consistently more accurate and, since they require no neural networks, also cheaper to compute. However, constructing closed-form functions that satisfy further interior properties, such as a prescribed Laplacian, in addition to the boundary conditions, is generally difficult without an optimization procedure that enforces these properties. This is precisely the role played by the neural-network training in B-NAVEM. Moreover, it is the proximity of the B-NAVEM functions to the VEM basis functions (see Theorem \ref{theor:bnavem:local_basis}) that makes the theoretical analysis of Sections~\ref{sec:discrete_space} and~\ref{sec:apriori} possible. Functions of type \ref{function:A} have no relation to the VEM space, so this analysis does not apply to them, and their stability and convergence properties are assessed only numerically.

\subsection{Test 2: Convergence test}\label{sec:test2}

Let us consider Problem \eqref{eq:continuous_poisson} on the domain $\Omega=(0,1)^2$, where the Dirichlet boundary condition and the source term are chosen in accordance with the following exact solution: 
\begin{equation*}
    u(\xx) = \sin(\pi x_1) \exp(x_2).
\end{equation*}

In this test, we evaluate the performance of the method by computing the errors \eqref{eq:errors} as the mesh size parameter $h$ decreases, for $k=1,2,3$ and for each function type \ref{function:C}, \ref{function:B}, and \ref{function:A}. For this experiment, we consider three families of five meshes each, obtained by iteratively refining the meshes shown in Figure \ref{fig:test1:mesh}.

\begin{table}[!ht]
\centering
\caption{Test 2: Convergence rates for $k =1, 2, 3$ and for the function type \ref{function:C}, computed on the data in Figure \ref{fig:test2:fc_k1}.}
\label{tab:test2:fc}
\begin{tabular}{lccc|cc|cc|cc|}
$k$ & Mesh Type    & \multicolumn{2}{c|}{$\ast = \BNAVEM$ ($k_q= 2k$)} & \multicolumn{2}{c|}{$\ast = \BNAVEM$ ($k_q= 30$)} & \multicolumn{2}{c|}{$\ast = \NAPEM$ ($k_q= 2k$)} & \multicolumn{2}{c|}{$\ast = \NAPEM$ ($k_q= 30$)} \\ \hline
    &              & ${\rm{err}}_{0}$      & ${\rm{err}}_{\nabla}$     & ${\rm{err}}_{0}$      & ${\rm{err}}_{\nabla}$     & ${\rm{err}}_{0}$     & ${\rm{err}}_{\nabla}$     & ${\rm{err}}_{0}$     & ${\rm{err}}_{\nabla}$     \\
1   & \ref{mesh:A} & 1.99                  & 0.99                      & 1.99                  & 0.99                      & 1.99                 & 0.99                      & 1.99                 & 0.99                      \\
2   & \ref{mesh:A} & 0.49                  & -0.63                     & 0.36                  & -0.62                     & 3.00                 & 1.99                      & 3.00                 & 1.99                      \\
3   & \ref{mesh:A} & 0.06                  & -0.75                     & -0.02                 & -0.79                     & 2.62                 & 1.51                      & 4.00                 & 3.00                      \\ \hline
1   & \ref{mesh:B} & 1.92                  & 0.93                      & 1.92                  & 0.93                      & 2.02                 & 1.01                      & 2.02                 & 1.01                      \\
2   & \ref{mesh:B} & 0.05                  & -0.78                     & -0.08                 & -0.85                     & 2.95                 & 1.93                      & 2.97                 & 1.96                      \\
3   & \ref{mesh:B} & 0.08                  & -0.67                     & 0.02                  & -0.69                     & 2.55                 & 1.29                      & 4.12                 & 3.08                      \\ \hline
1   & \ref{mesh:C} & -0.21                 & -0.43                     & -0.15                 & -0.45                     & 1.75                 & 1.11                      & 1.84                 & 1.19                      \\
2   & \ref{mesh:C} & -1.23                 & -0.73                     & -1.44                 & -0.73                     & 2.61                 & 1.43                      & 2.91                 & 1.89                      \\
3   & \ref{mesh:C} & -1.76                 & -0.67                     & -1.78                 & -0.79                     & 0.87                 & 0.01                      & 3.97                 & 2.87                      \\ \hline
\end{tabular}
\end{table}

\begin{table}[!ht]
\centering
\caption{Test 2: Convergence rates for $k =1, 2, 3$ and for the function type \ref{function:B}, computed on the data in Figure \ref{fig:test2:fb_k1}.}
\label{tab:test2:fb}
\begin{tabular}{lccc|cc|cc|cc|}
$k$ & Mesh Type    & \multicolumn{2}{c|}{$\ast = \BNAVEM$ ($k_q= 2k$)} & \multicolumn{2}{c|}{$\ast = \BNAVEM$ ($k_q= 30$)} & \multicolumn{2}{c|}{$\ast = \NAPEM$ ($k_q= 2k$)} & \multicolumn{2}{c|}{$\ast = \NAPEM$ ($k_q= 30$)} \\ \hline
    &              & ${\rm{err}}_{0}$      & ${\rm{err}}_{\nabla}$     & ${\rm{err}}_{0}$      & ${\rm{err}}_{\nabla}$     & ${\rm{err}}_{0}$     & ${\rm{err}}_{\nabla}$     & ${\rm{err}}_{0}$     & ${\rm{err}}_{\nabla}$     \\
1   & \ref{mesh:A} & 2.24                  & 0.98                      & 2.14                  & 0.92                      & 1.99                 & 0.99                      & 1.99                 & 0.99                      \\
2   & \ref{mesh:A} & -1.31                 & -0.71                     & -1.33                 & -0.72                     & 2.99                 & 1.98                      & 3.00                 & 1.99                      \\
3   & \ref{mesh:A} & -1.58                 & -0.84                     & -1.54                 & -0.83                     & 0.99                 & -0.11                     & 3.99                 & 2.99                      \\ \hline
1   & \ref{mesh:B} & 2.03                  & 0.93                      & 1.95                  & 0.87                      & 2.02                 & 1.01                      & 2.02                 & 1.01                      \\
2   & \ref{mesh:B} & -1.28                 & -0.78                     & -1.32                 & -0.79                     & 2.81                 & 1.79                      & 2.97                 & 1.96                      \\
3   & \ref{mesh:B} & -1.49                 & -0.85                     & -1.47                 & -0.83                     & 1.01                 & -0.09                     & 4.13                 & 3.08                      \\ \hline
1   & \ref{mesh:C} & 0.26                  & -0.23                     & 0.27                  & -0.23                     & 1.94                 & 1.07                      & 1.94                 & 1.07                      \\
2   & \ref{mesh:C} & -0.67                 & -0.53                     & -1.13                 & -0.63                     & 2.71                 & 1.43                      & 2.91                 & 1.89                      \\
3   & \ref{mesh:C} & -1.76                 & -0.74                     & -1.79                 & -0.71                     & 1.00                 & 0.09                      & 2.77                 & 1.32                      \\ \hline
\end{tabular}%
\end{table}

\begin{table}[!ht]
\centering
\caption{Test 2: Convergence rates for $k =1, 2, 3$ and for the function type \ref{function:A}, computed on the data in Figure \ref{fig:test2:fa_k1}.}
\label{tab:test2:fa}
\begin{tabular}{lccc|cc|cc|cc|}
$k$ & Mesh Type    & \multicolumn{2}{c|}{$\ast = \BNAVEM$ ($k_q= 2k$)} & \multicolumn{2}{c|}{$\ast = \BNAVEM$ ($k_q= 30$)} & \multicolumn{2}{c|}{$\ast = \NAPEM$ ($k_q= 2k$)} & \multicolumn{2}{c|}{$\ast = \NAPEM$ ($k_q= 30$)} \\ \hline
    &              & ${\rm{err}}_{0}$      & ${\rm{err}}_{\nabla}$     & ${\rm{err}}_{0}$      & ${\rm{err}}_{\nabla}$     & ${\rm{err}}_{0}$     & ${\rm{err}}_{\nabla}$     & ${\rm{err}}_{0}$     & ${\rm{err}}_{\nabla}$     \\
1   & \ref{mesh:A} & 1.99                  & 0.99                      & 1.99                  & 0.99                      & 1.99                 & 0.99                      & 1.99                 & 0.99                      \\
2   & \ref{mesh:A} & 2.98                  & 1.99                      & 2.98                  & 1.99                      & 2.98                 & 1.99                      & 2.98                 & 1.99                      \\
3   & \ref{mesh:A} & 3.05                  & 2.03                      & 3.05                  & 2.03                      & 3.99                 & 2.99                      & 3.99                 & 2.99                      \\ \hline
1   & \ref{mesh:B} & 2.02                  & 1.01                      & 2.02                  & 1.01                      & 2.02                 & 1.01                      & 2.02                 & 1.01                      \\
2   & \ref{mesh:B} & 2.53                  & 1.40                      & 2.55                  & 1.39                      & 2.98                 & 1.99                      & 2.98                 & 1.99                      \\
3   & \ref{mesh:B} & 2.00                  & 0.96                      & 2.00                  & 0.96                      & 4.07                 & 3.02                      & 4.10                 & 3.07                      \\ \hline
1   & \ref{mesh:C} & 1.95                  & 1.03                      & 1.95                  & 1.03                      & 1.95                 & 1.03                      & 1.95                 & 1.03                      \\
2   & \ref{mesh:C} & 2.30                  & 1.29                      & 2.32                  & 1.29                      & 2.45                 & 1.35                      & 2.92                 & 1.94                      \\
3   & \ref{mesh:C} & 2.28                  & 1.28                      & 2.29                  & 1.29                      & 1.46                 & 0.45                      & 4.03                 & 2.98                      \\ \hline
\end{tabular}%
\end{table}

\begin{figure}[!ht]
    \centering
    \begin{subfigure}{0.3\textwidth}
        \includegraphics[width=1\textwidth]{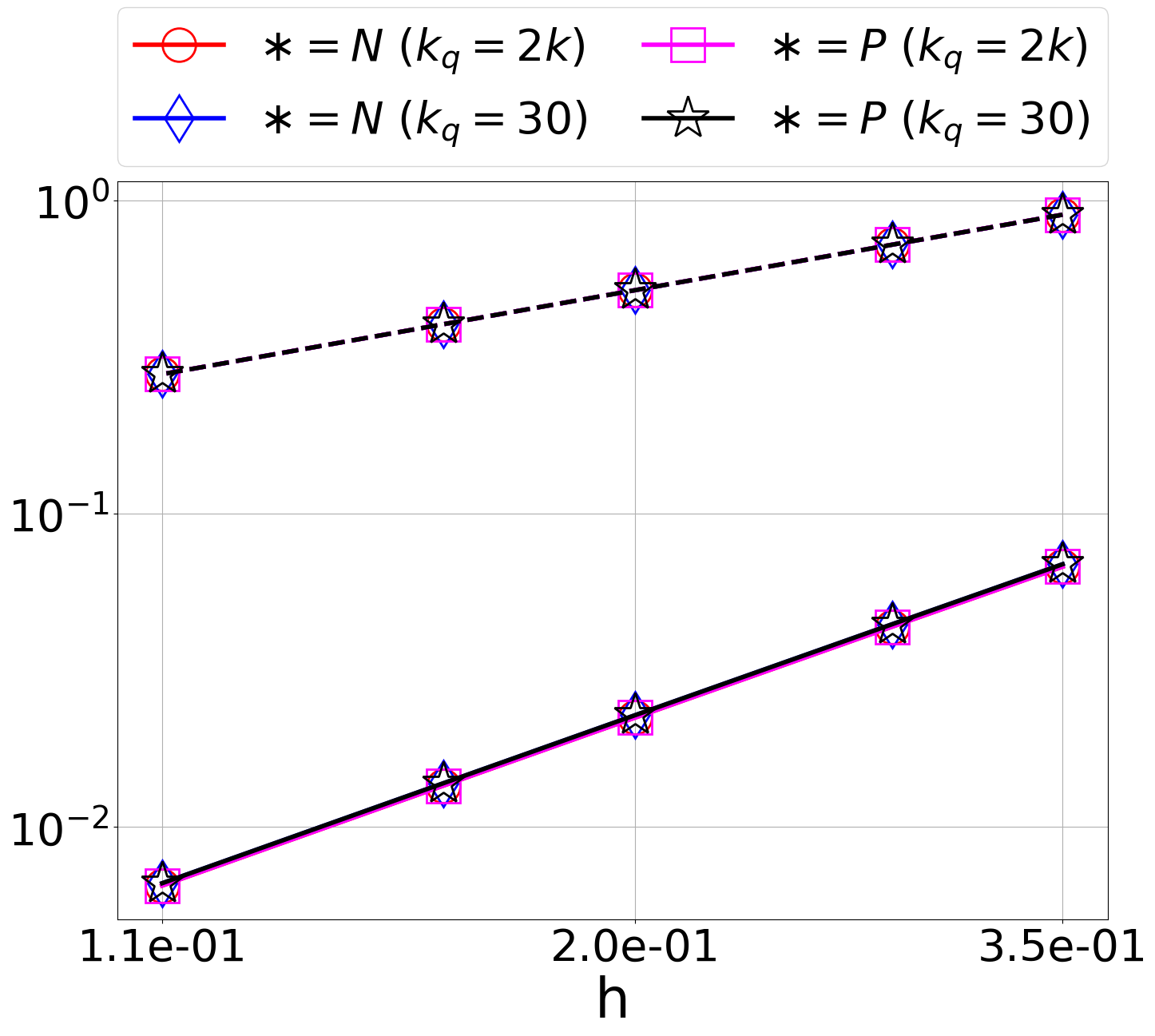}
        \caption{}
    \end{subfigure}
    \begin{subfigure}{0.3\textwidth}
        \includegraphics[width=1\textwidth]{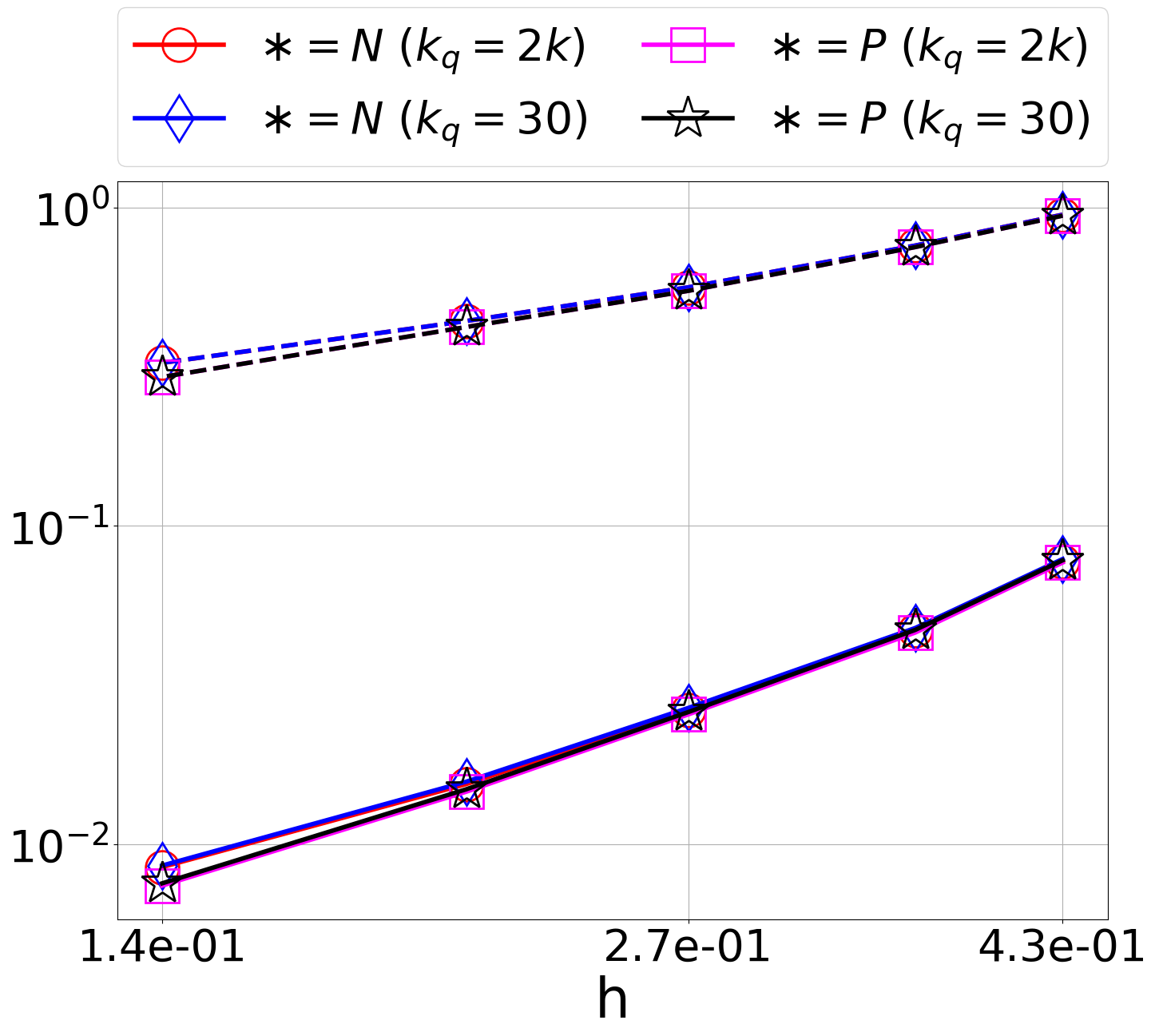}
        \caption{}
    \end{subfigure}
    \begin{subfigure}{0.3\textwidth}
        \includegraphics[width=1\textwidth]{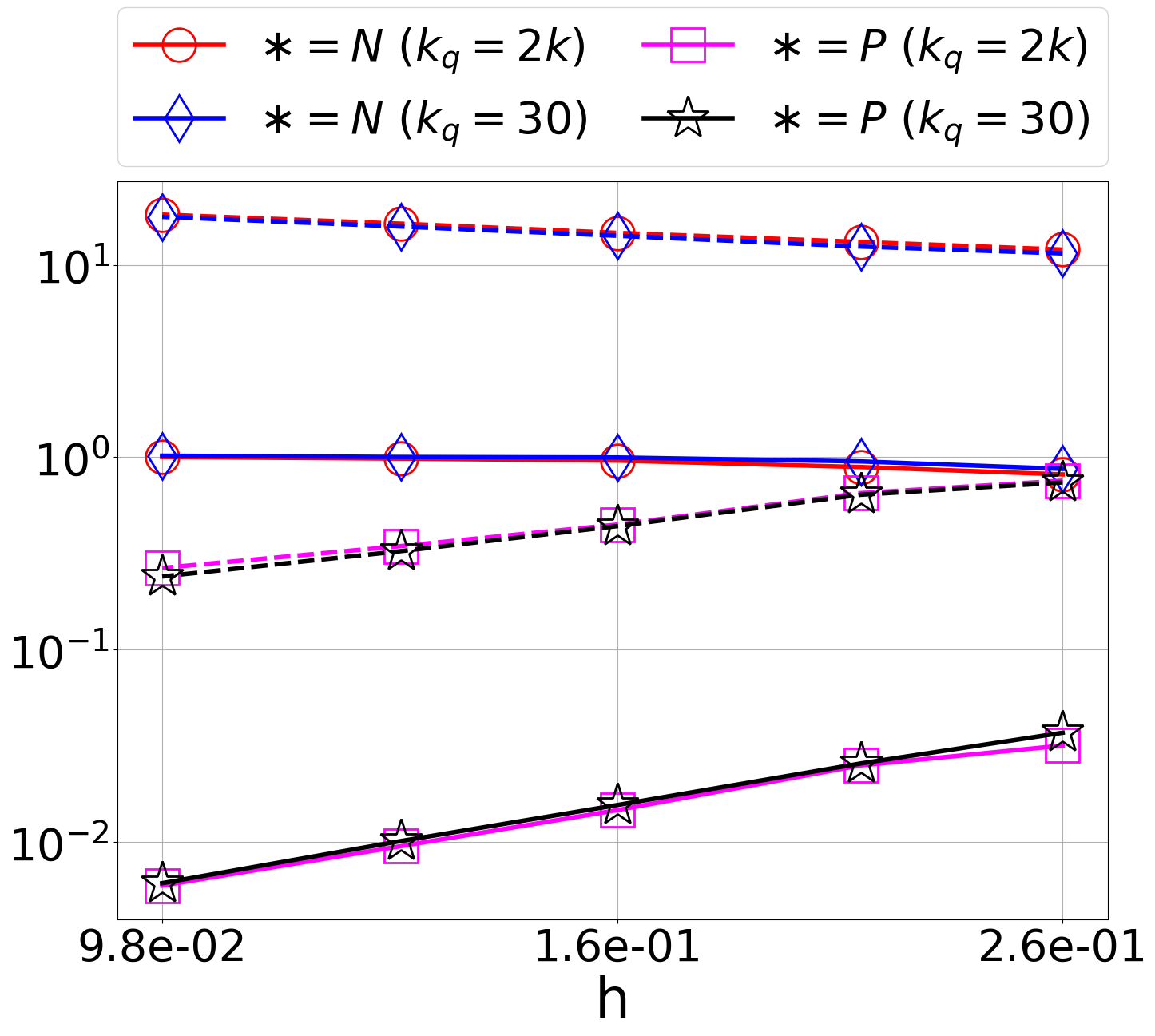}
        \caption{}
    \end{subfigure}
    \begin{subfigure}{0.3\textwidth}
        \includegraphics[width=1\textwidth]{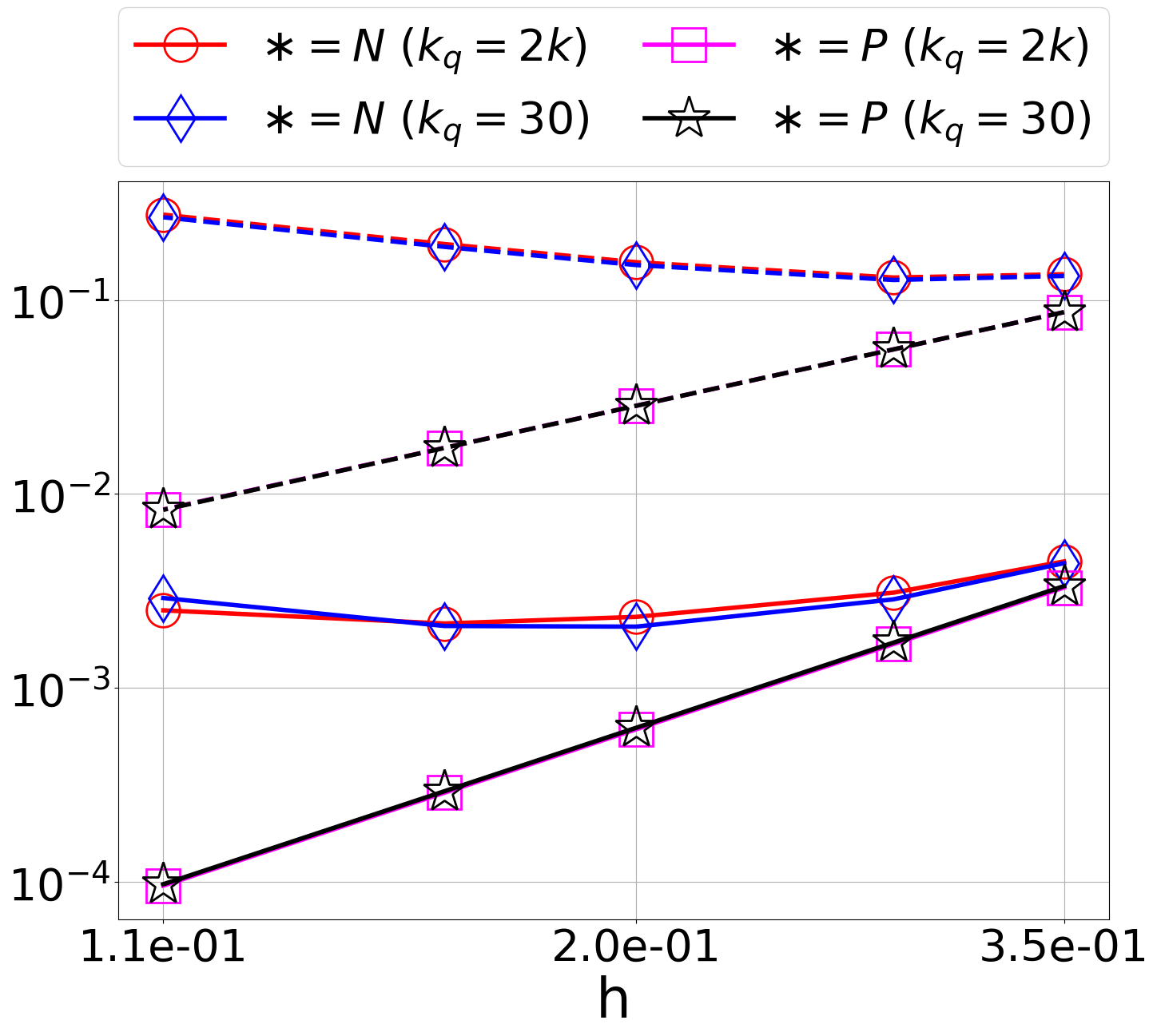}
        \caption{}
    \end{subfigure}
    \begin{subfigure}{0.3\textwidth}
        \includegraphics[width=1\textwidth]{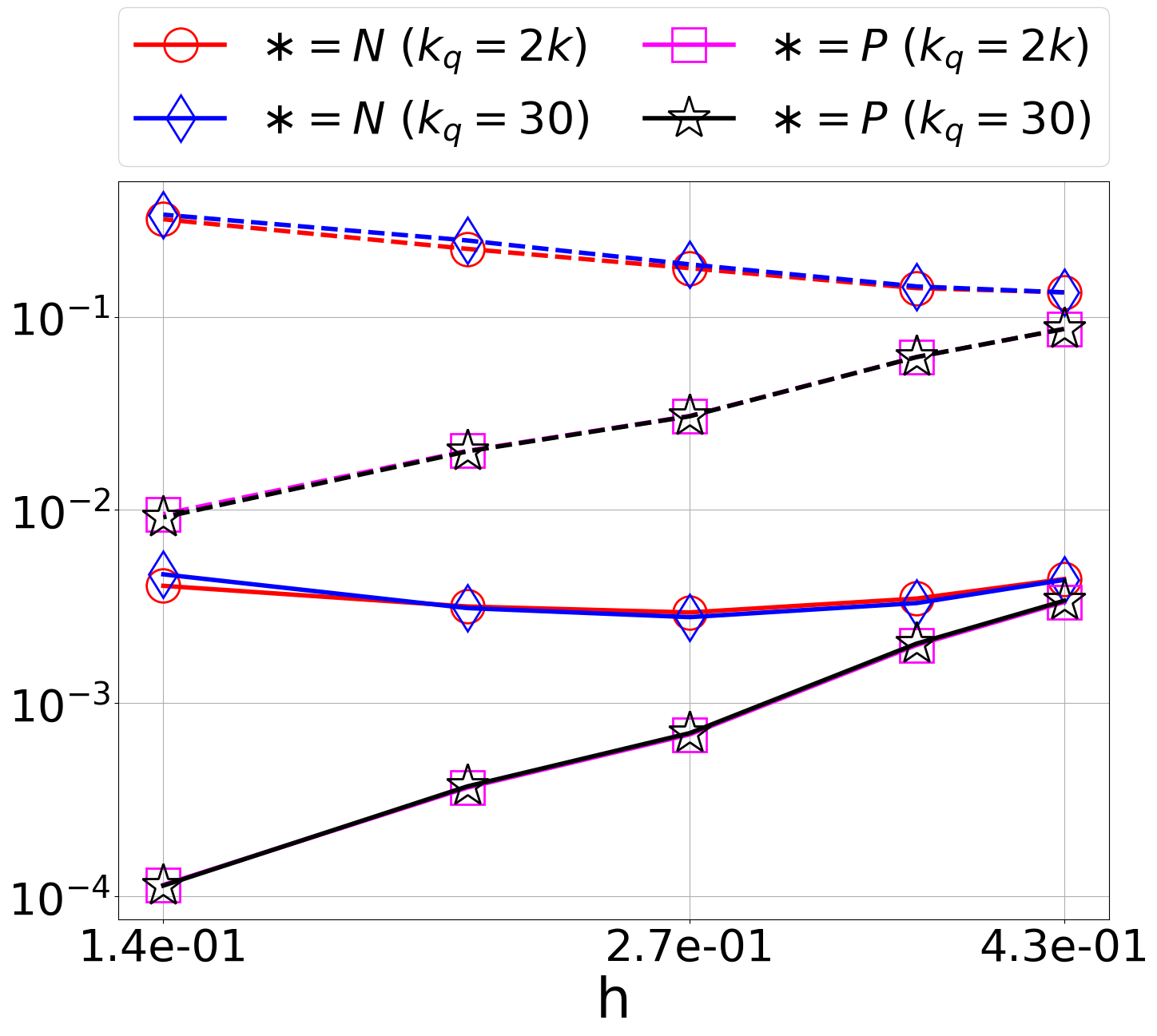}
        \caption{}
    \end{subfigure}
    \begin{subfigure}{0.3\textwidth}
        \includegraphics[width=1\textwidth]{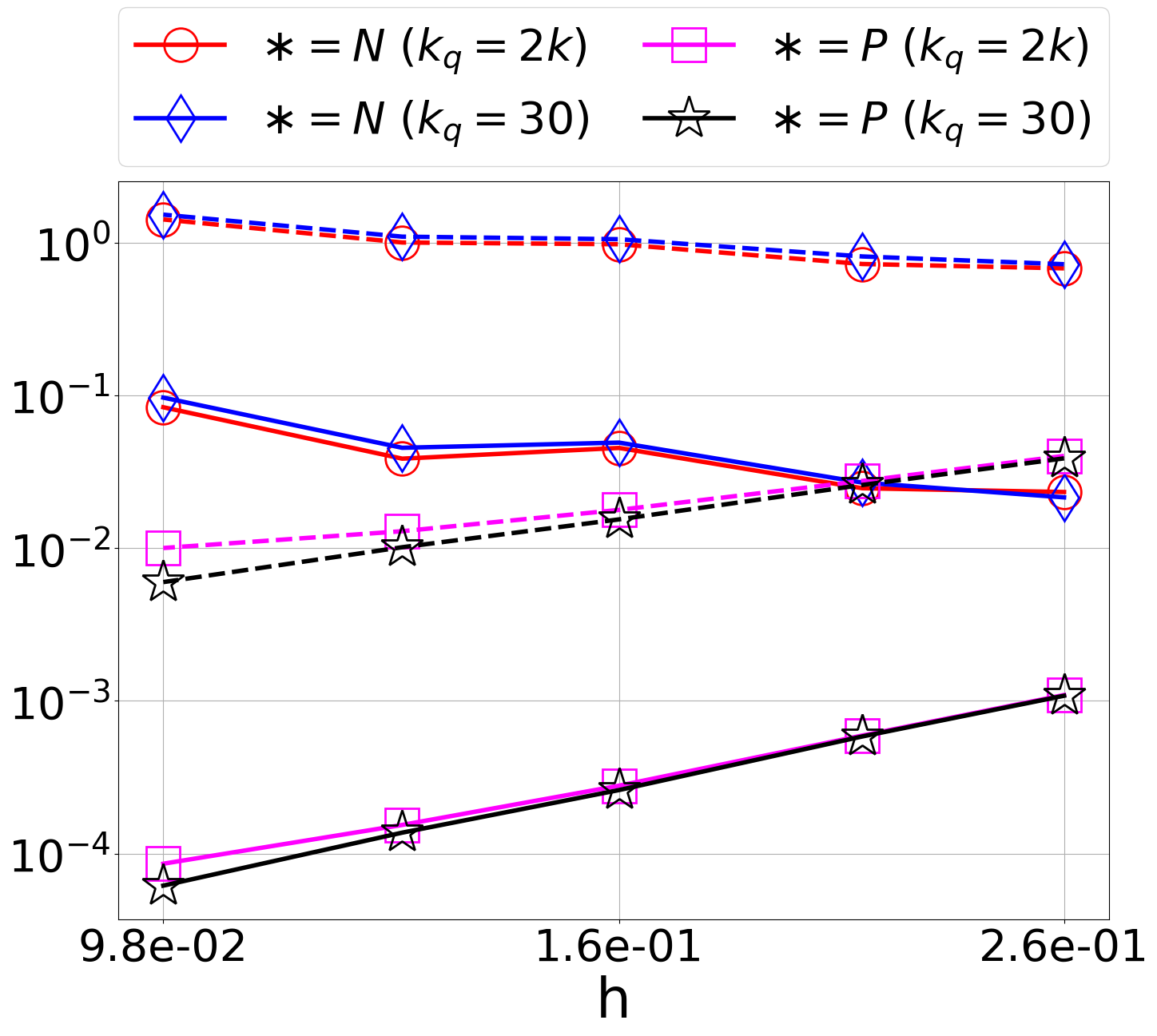}
        \caption{}
    \end{subfigure}
    \begin{subfigure}{0.3\textwidth}
        \includegraphics[width=1\textwidth]{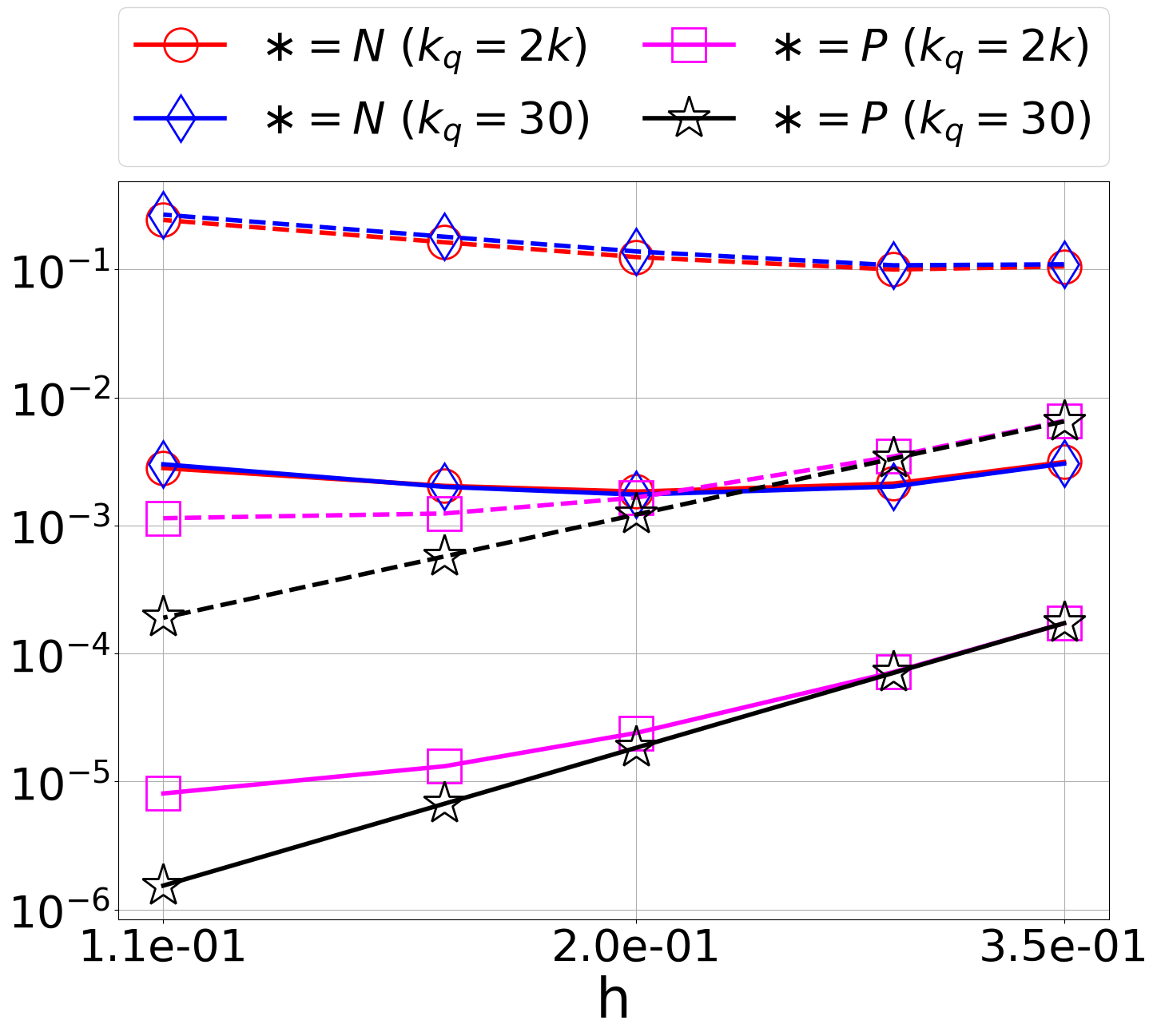}
        \caption{}
    \end{subfigure}
    \begin{subfigure}{0.3\textwidth}
        \includegraphics[width=1\textwidth]{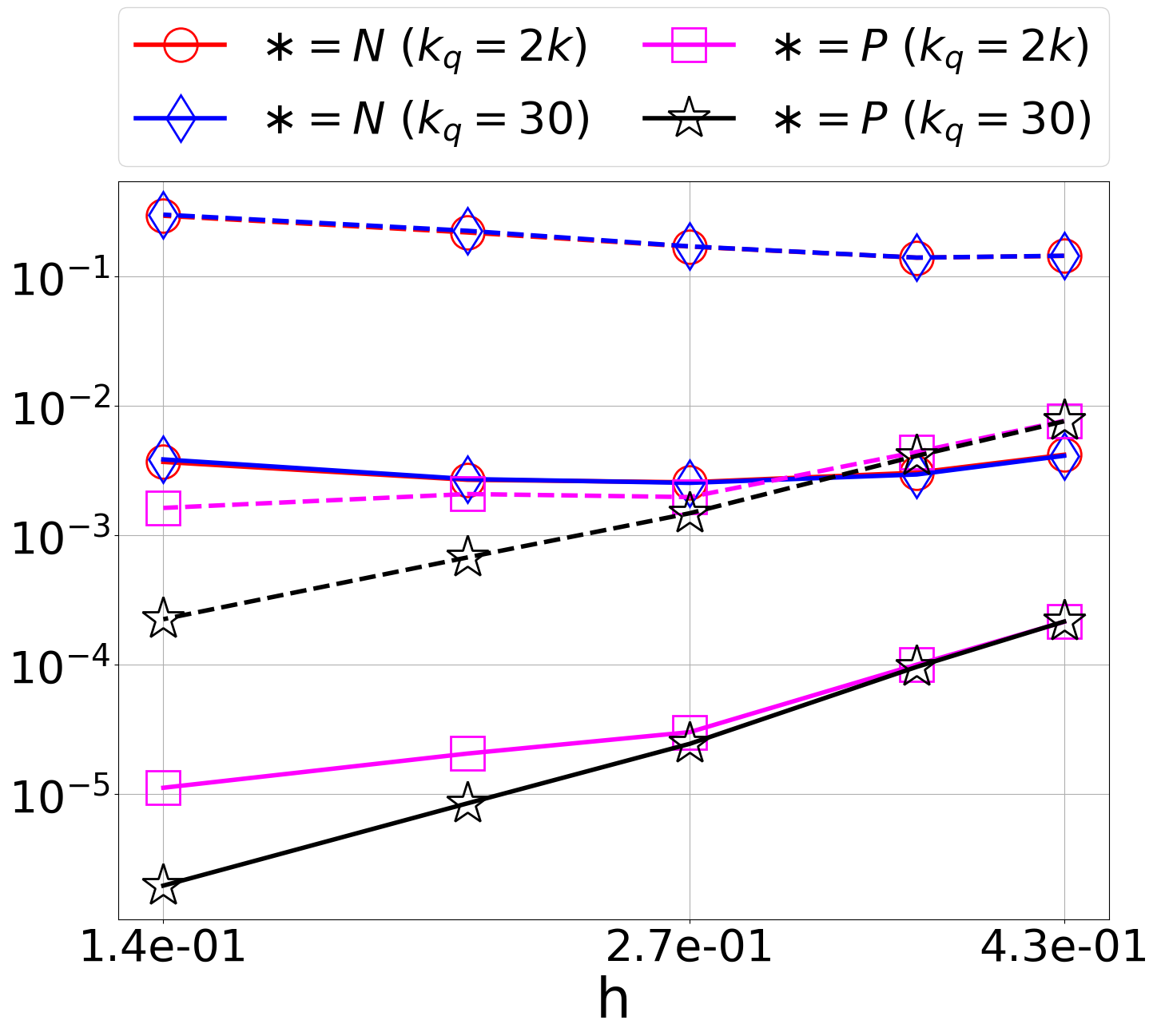}
        \caption{}
    \end{subfigure}
    \begin{subfigure}{0.3\textwidth}
        \includegraphics[width=1\textwidth]{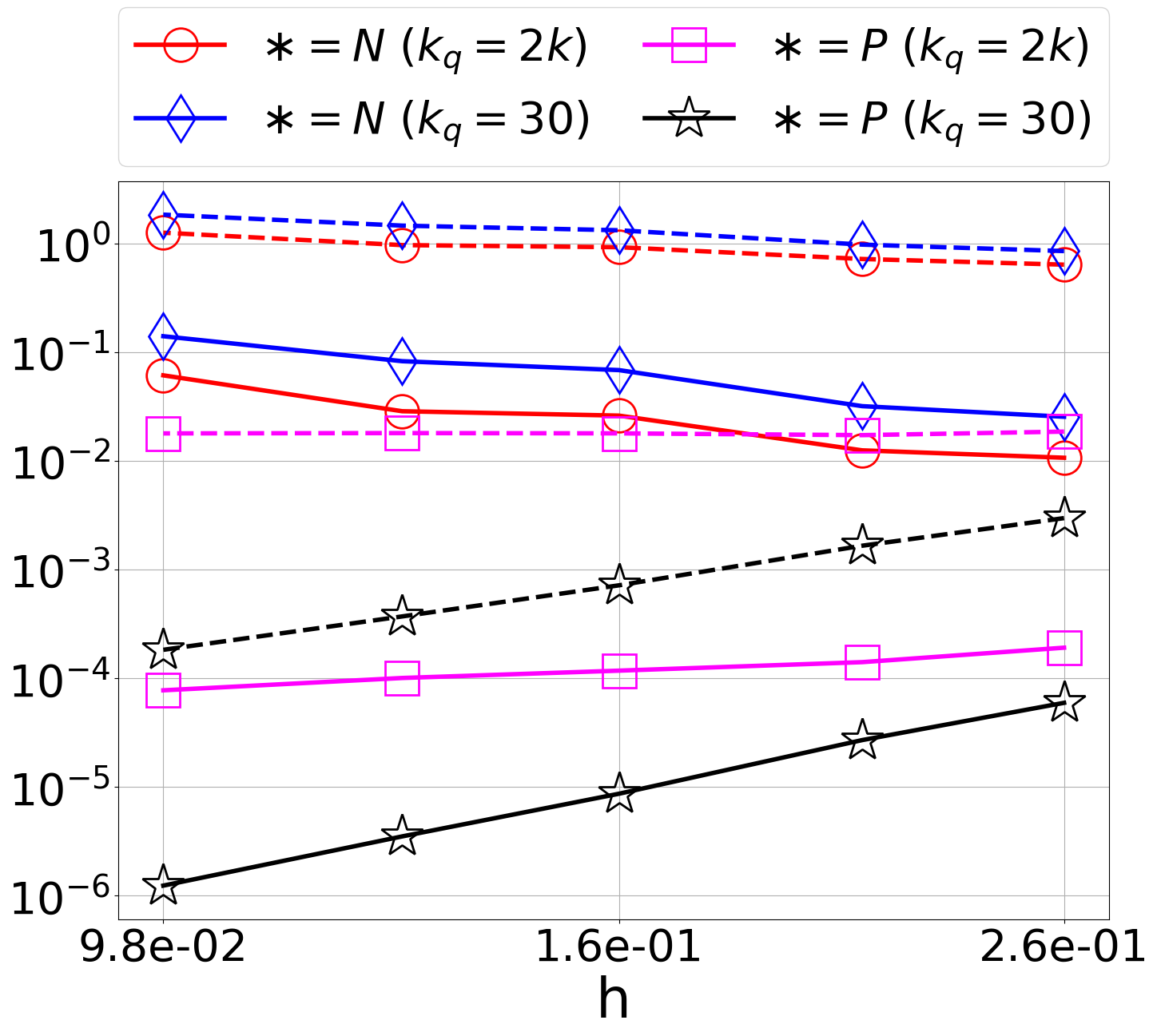}
        \caption{}
    \end{subfigure}
    \caption{Test 2: Behavior of the errors \eqref{eq:errors} as $h$ decreases for the function type \ref{function:C}. Solid lines represent error in $\leb{2}{}$-norm, while dashed lines denote the error in $\sob{1}{}$-seminorm. Left: \ref{mesh:A}. Center: \ref{mesh:B}. Right: \ref{mesh:C}. First row: $k=1$. Second row: $k=2$. Third row: $k=3$.}
    \label{fig:test2:fc_k1}
\end{figure}

\begin{figure}[!ht]
    \centering
    \begin{subfigure}{0.3\textwidth}
        \includegraphics[width=1\textwidth]{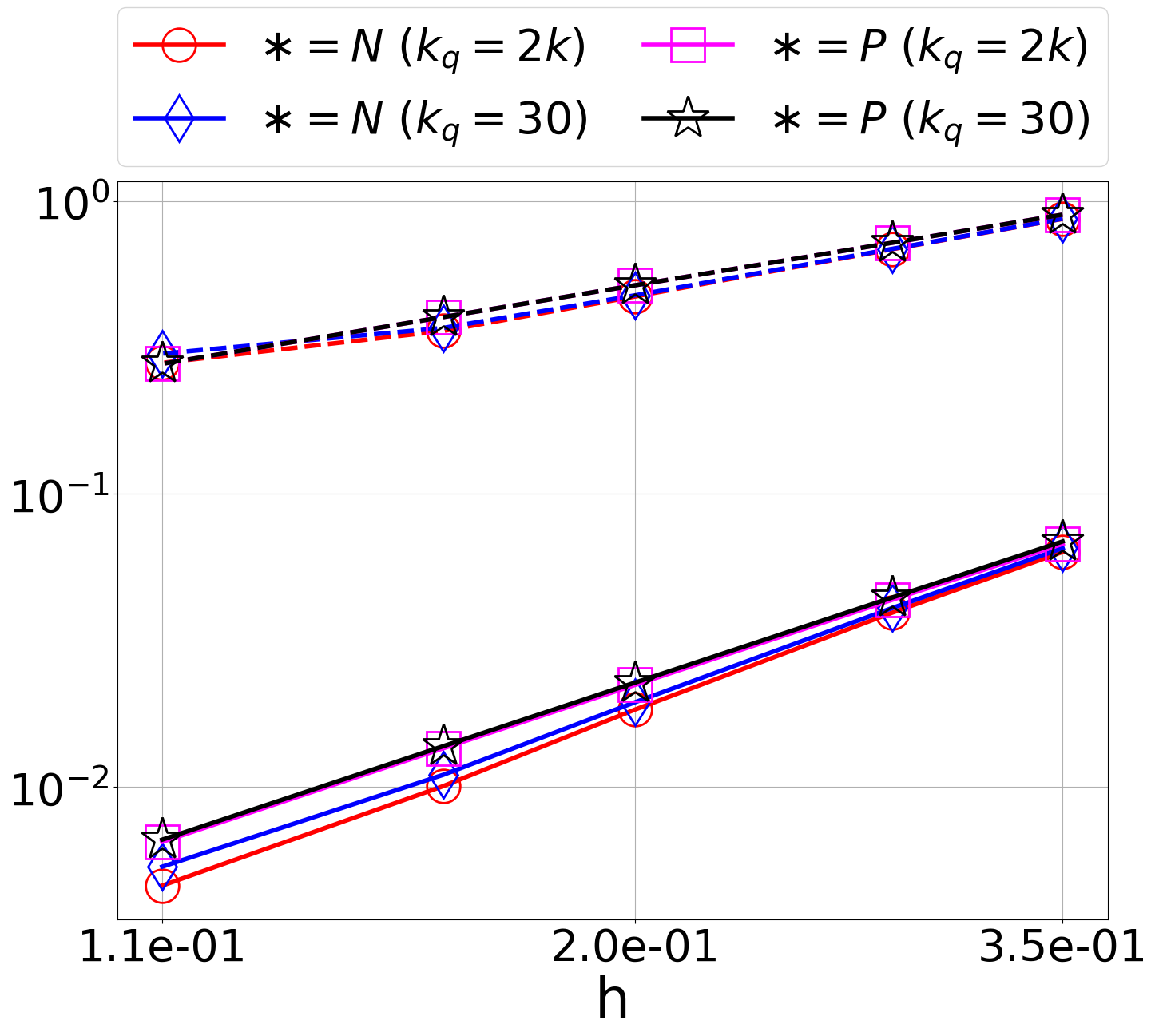}
        \caption{}
    \end{subfigure}
    \begin{subfigure}{0.3\textwidth}
        \includegraphics[width=1\textwidth]{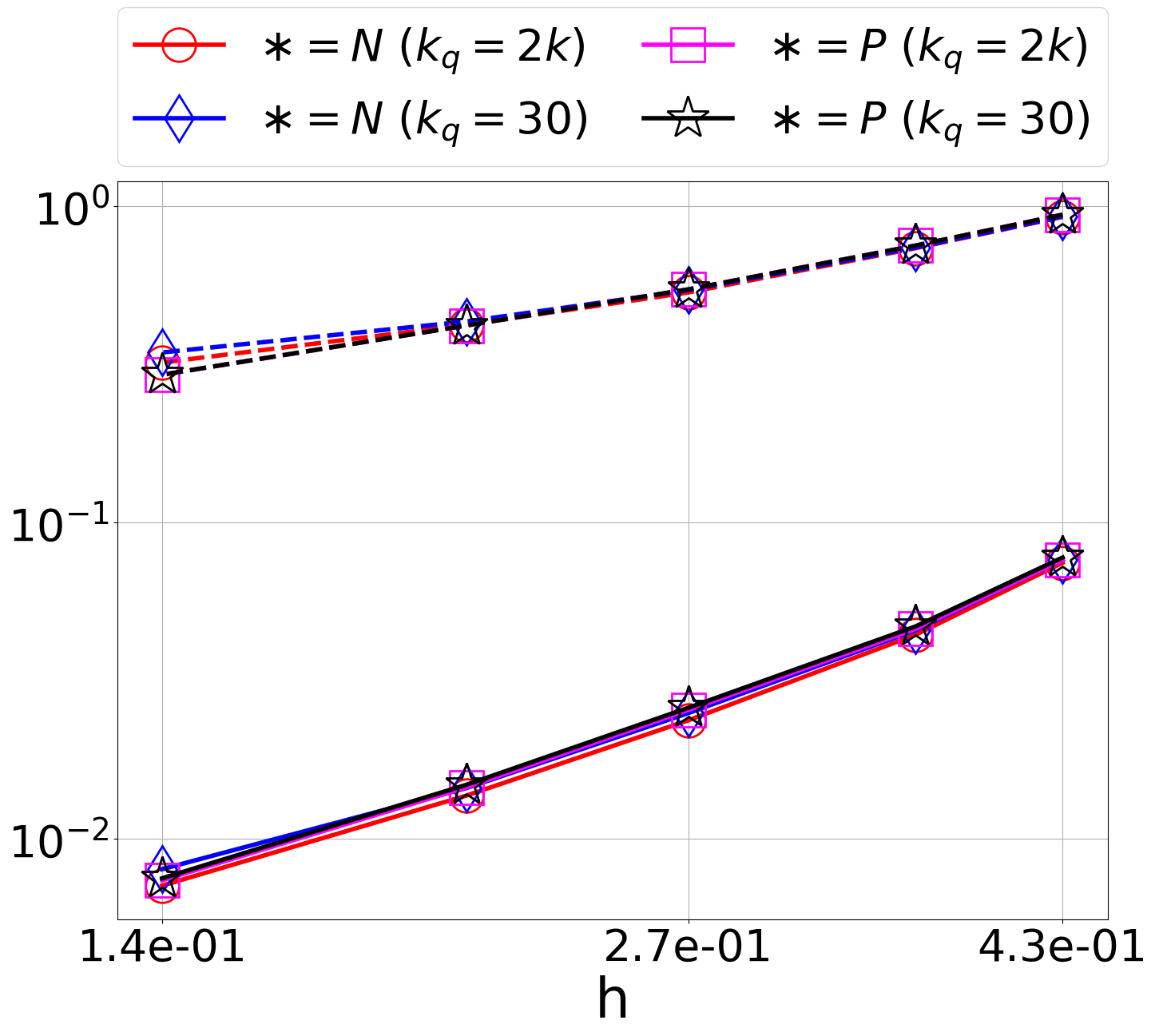}
        \caption{}
    \end{subfigure}
    \begin{subfigure}{0.3\textwidth}
        \includegraphics[width=1\textwidth]{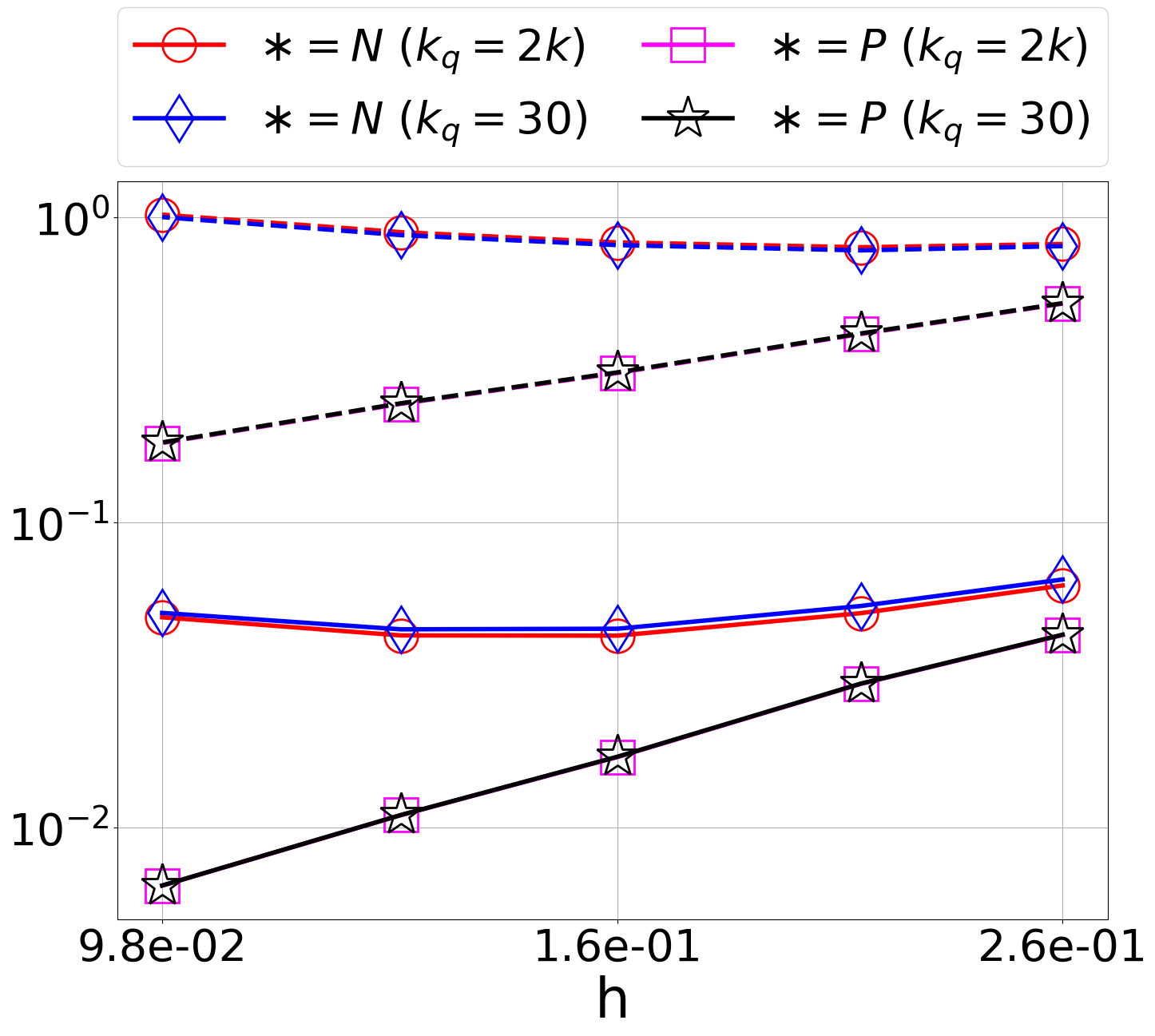}
        \caption{}
    \end{subfigure}
    \begin{subfigure}{0.3\textwidth}
        \includegraphics[width=1\textwidth]{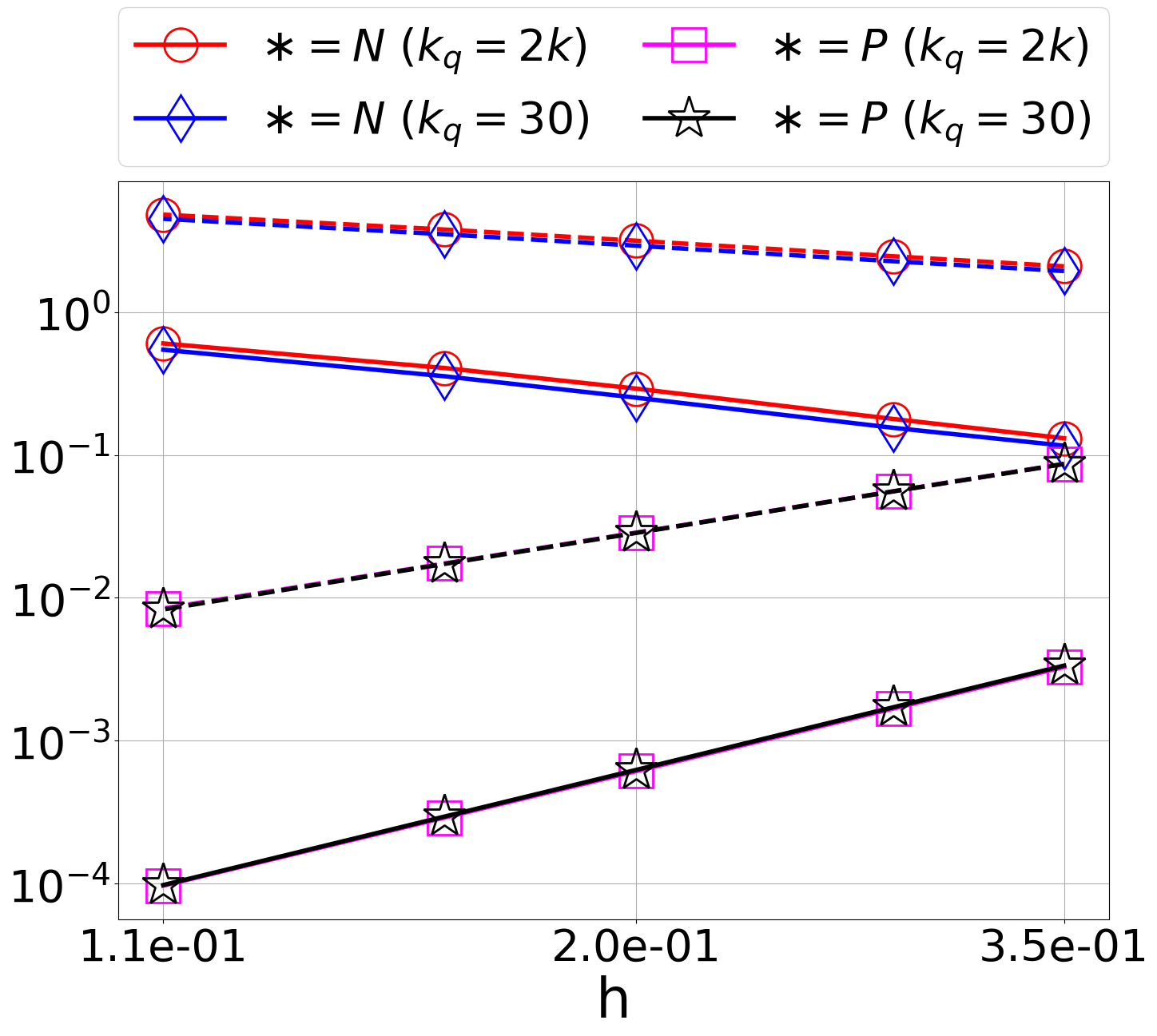}
        \caption{}
    \end{subfigure}
    \begin{subfigure}{0.3\textwidth}
        \includegraphics[width=1\textwidth]{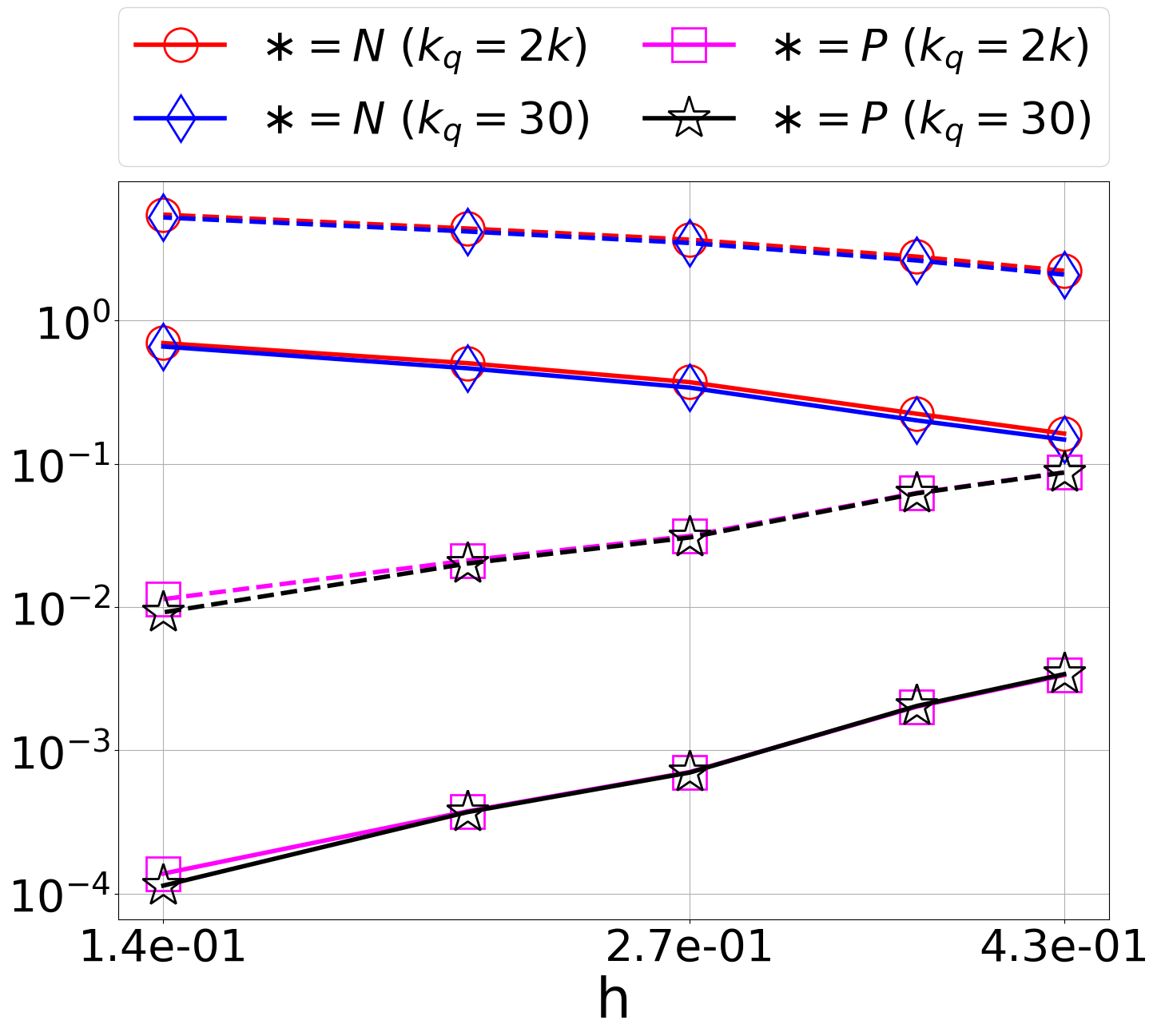}
        \caption{}
    \end{subfigure}
    \begin{subfigure}{0.3\textwidth}
        \includegraphics[width=1\textwidth]{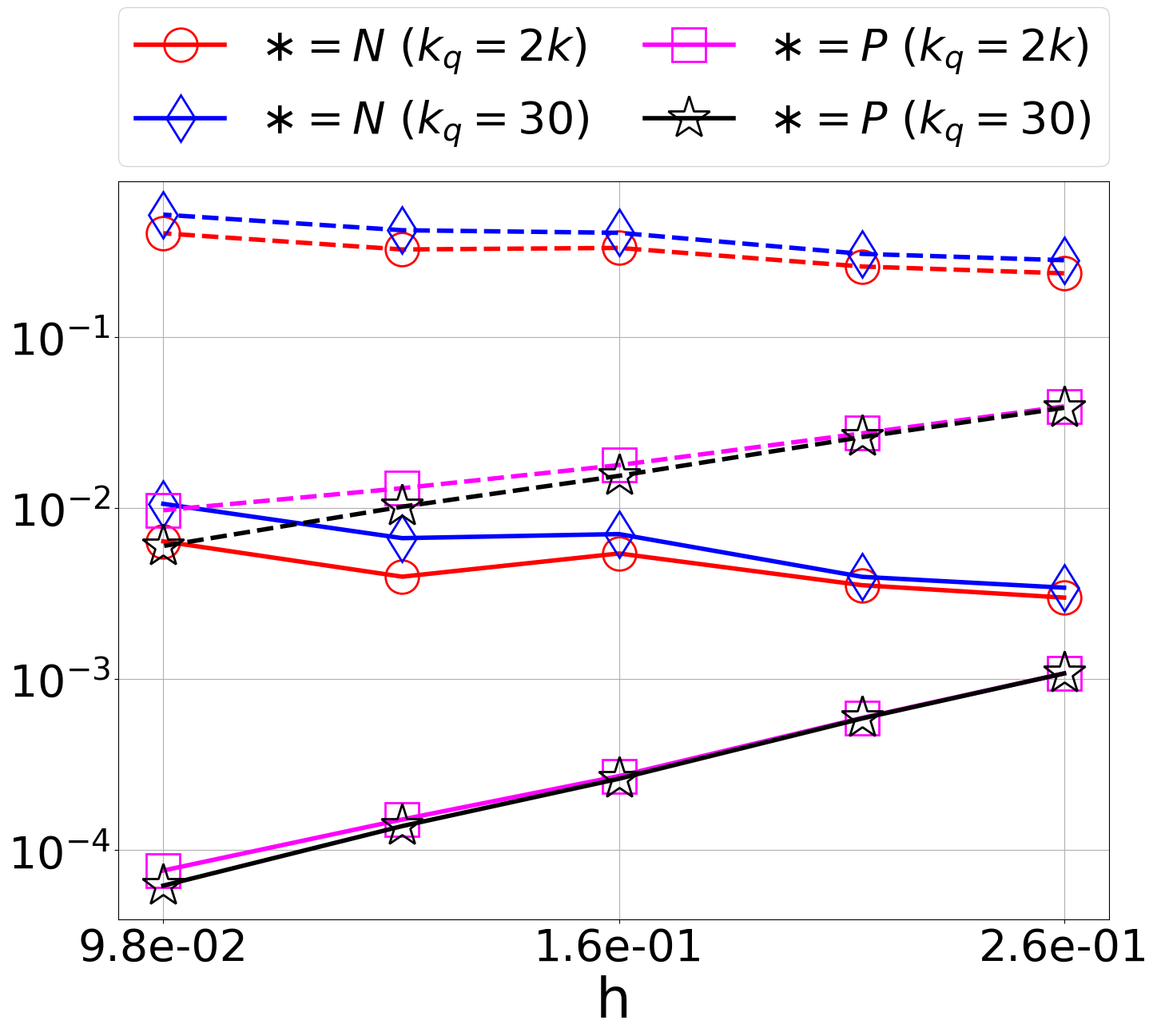}
        \caption{}
    \end{subfigure}
    \begin{subfigure}{0.3\textwidth}
        \includegraphics[width=1\textwidth]{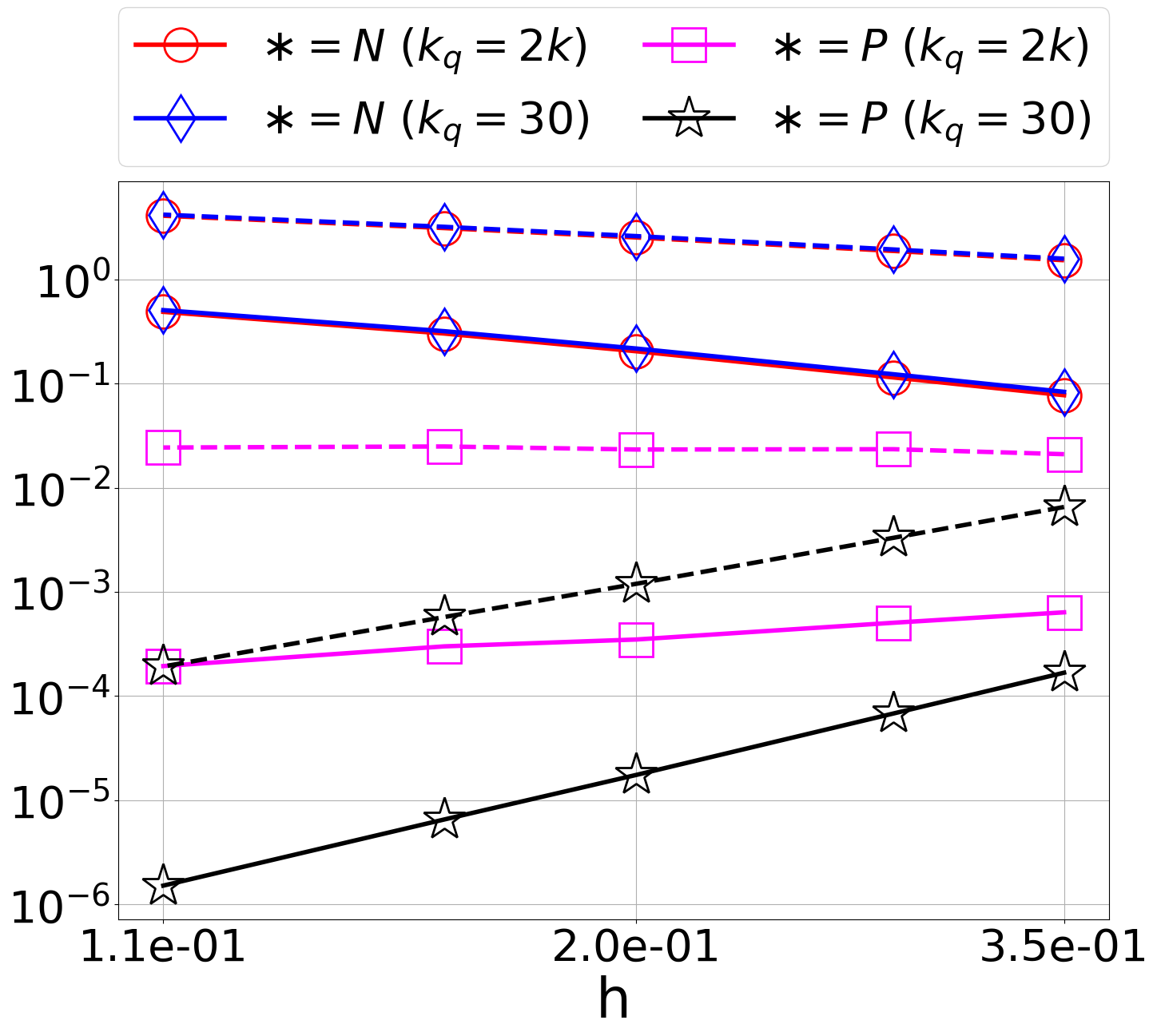}
        \caption{}
    \end{subfigure}
    \begin{subfigure}{0.3\textwidth}
        \includegraphics[width=1\textwidth]{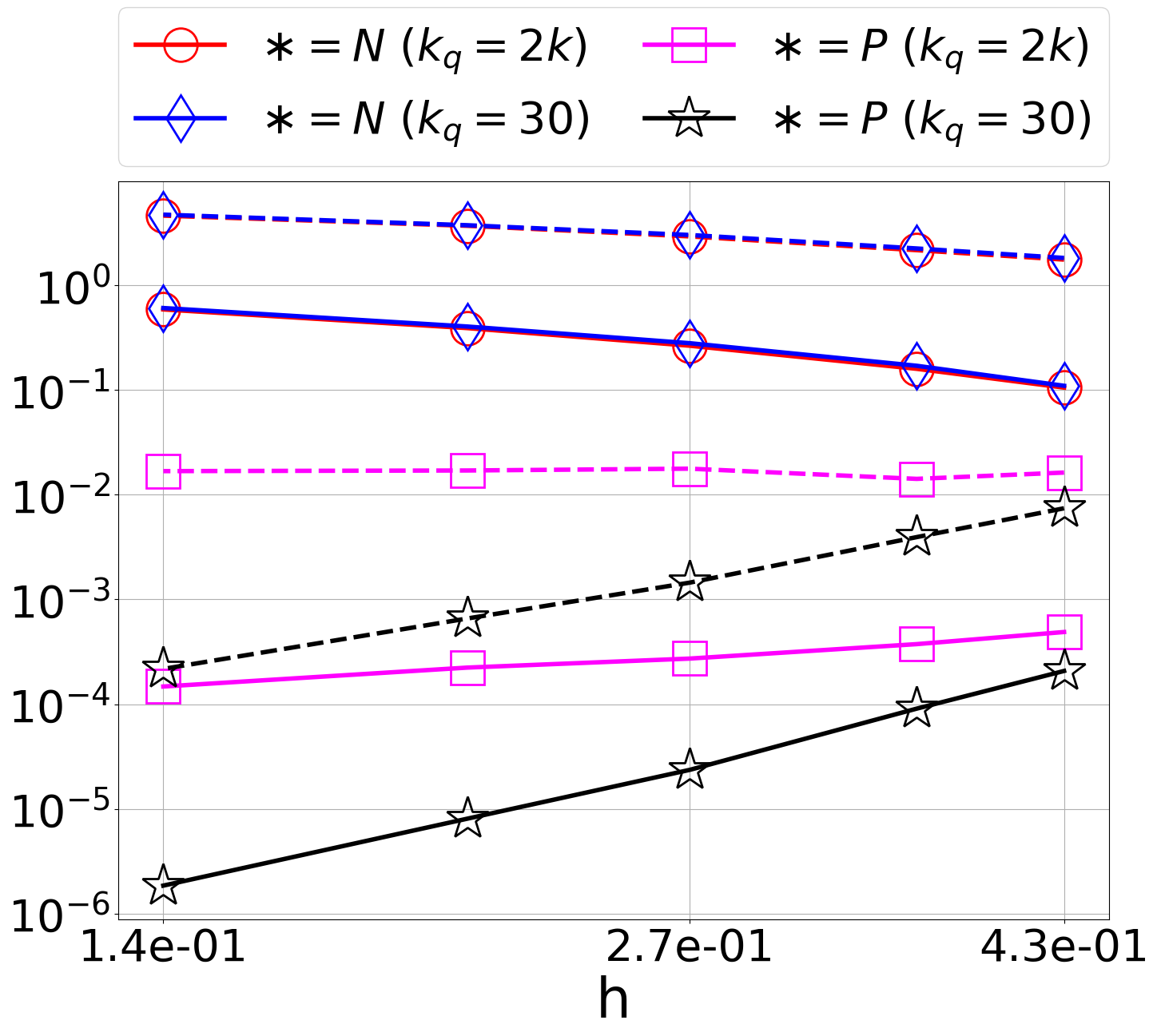}
        \caption{}
    \end{subfigure}
    \begin{subfigure}{0.3\textwidth}
        \includegraphics[width=1\textwidth]{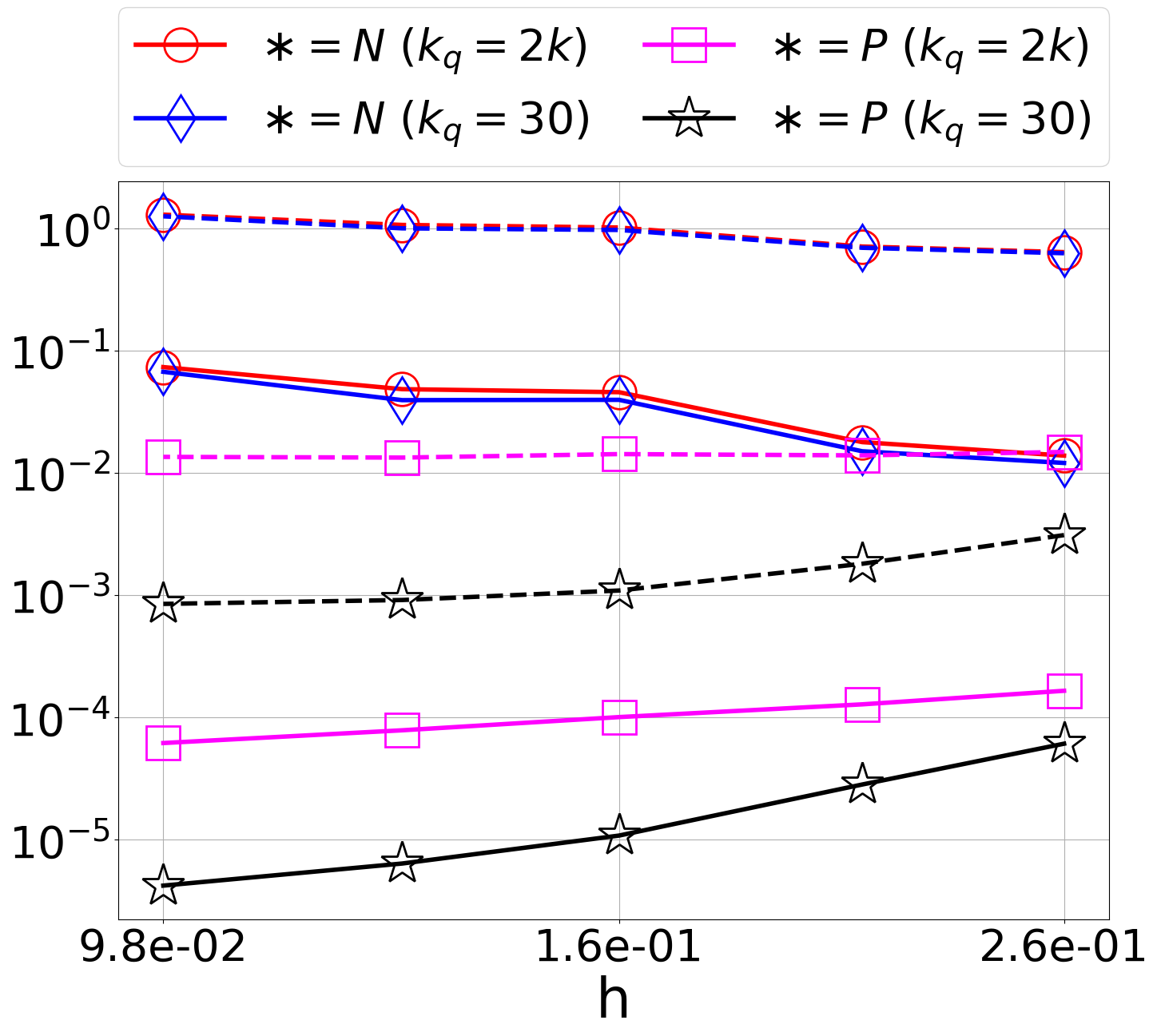}
        \caption{}
    \end{subfigure}
    \caption{Test 2: Behavior of the errors \eqref{eq:errors} as $h$ decreases for the function type \ref{function:B}. Solid lines represent error in $\leb{2}{}$-norm, while dashed lines denote the error in $\sob{1}{}$-seminorm. Left: \ref{mesh:A}. Center: \ref{mesh:B}. Right: \ref{mesh:C}. First row: $k=1$. Second row: $k=2$. Third row: $k=3$.}
    \label{fig:test2:fb_k1}
\end{figure}

\begin{figure}[!ht]
    \centering
    \begin{subfigure}{0.3\textwidth}
        \includegraphics[width=1\textwidth]{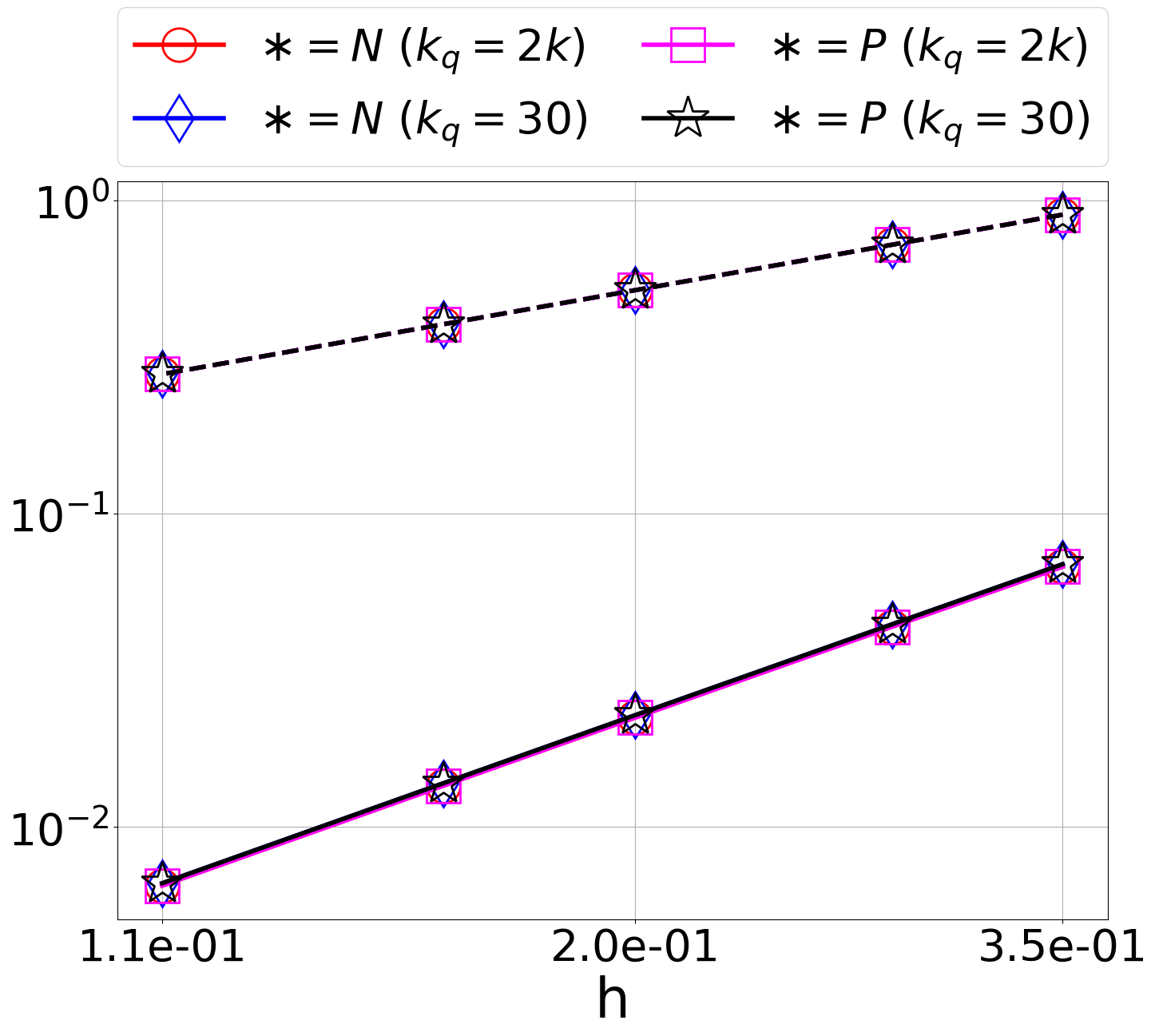}
        \caption{}
    \end{subfigure}
    \begin{subfigure}{0.3\textwidth}
        \includegraphics[width=1\textwidth]{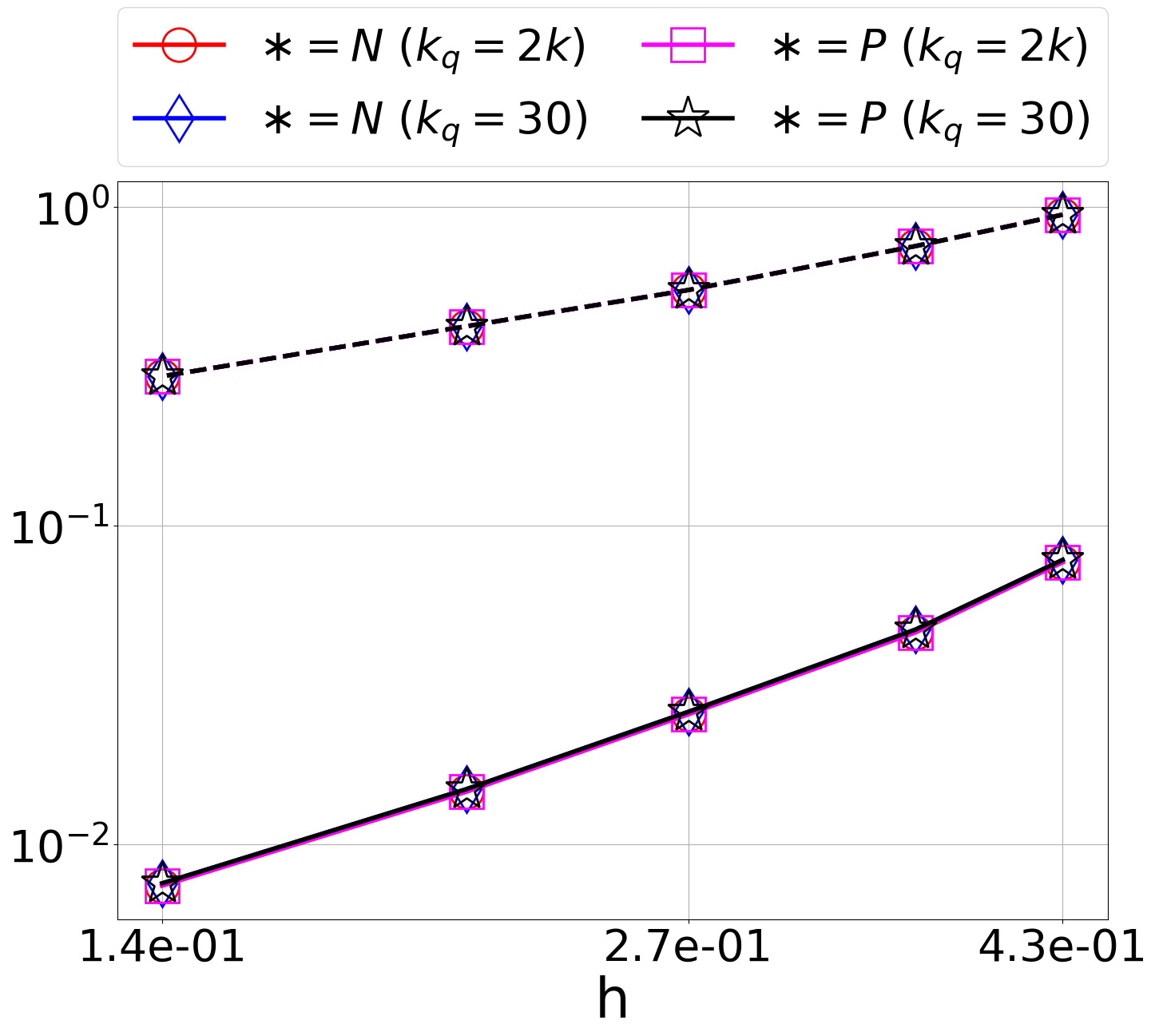}
        \caption{}
    \end{subfigure}
    \begin{subfigure}{0.3\textwidth}
        \includegraphics[width=1\textwidth]{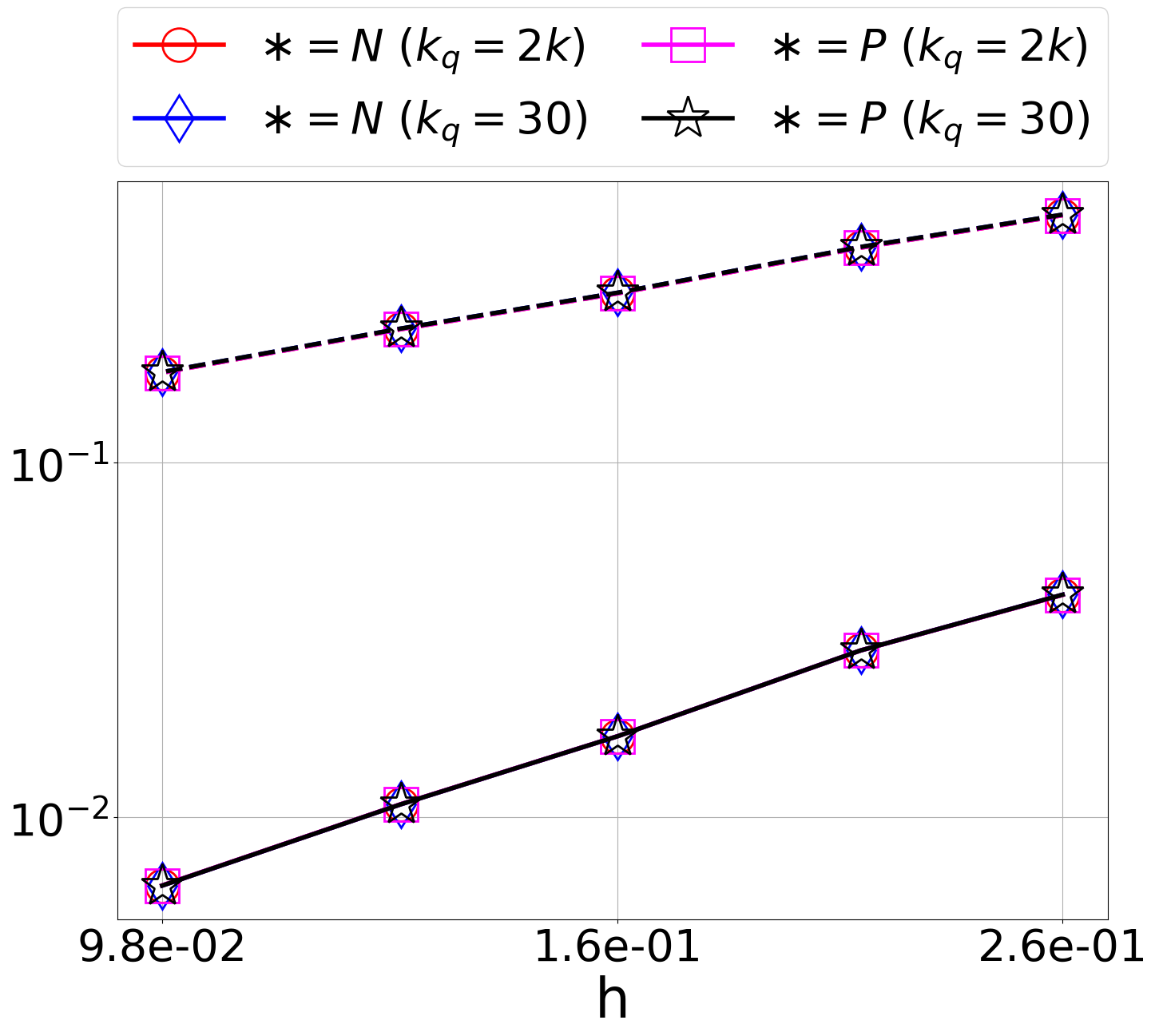}
        \caption{}
    \end{subfigure}
    \begin{subfigure}{0.3\textwidth}
        \includegraphics[width=1\textwidth]{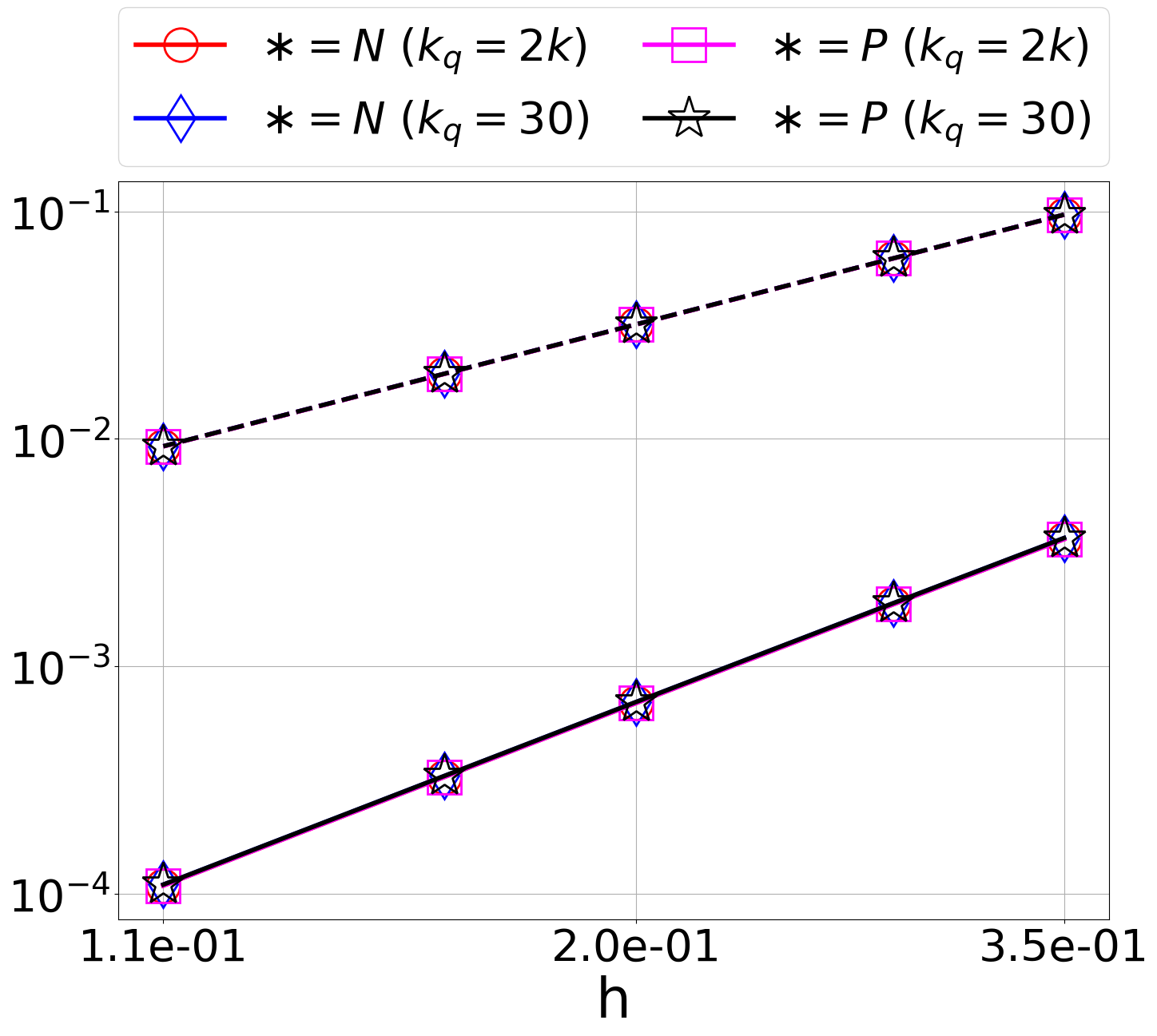}
        \caption{}
    \end{subfigure}
    \begin{subfigure}{0.3\textwidth}
        \includegraphics[width=1\textwidth]{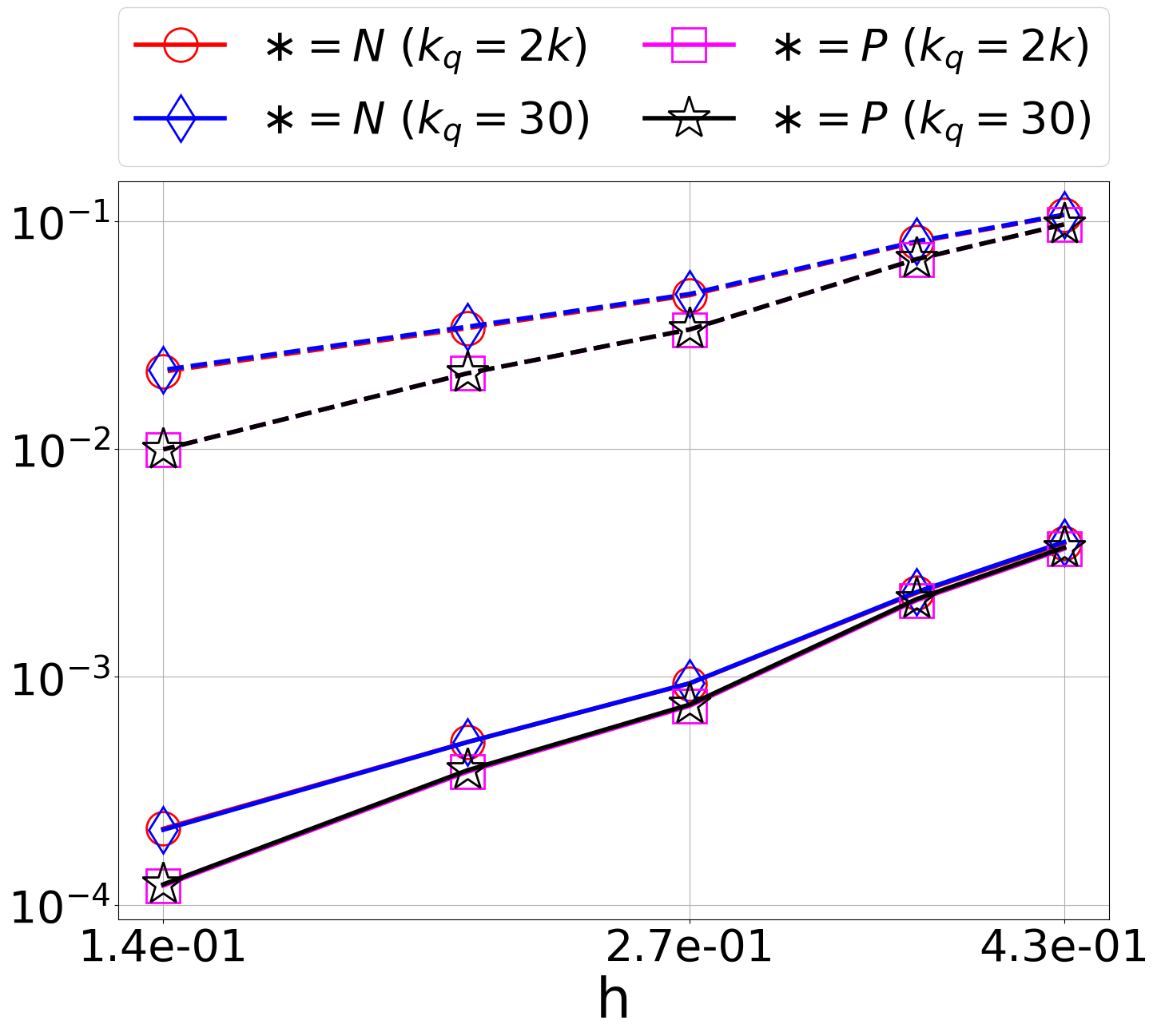}
        \caption{}
    \end{subfigure}
    \begin{subfigure}{0.3\textwidth}
        \includegraphics[width=1\textwidth]{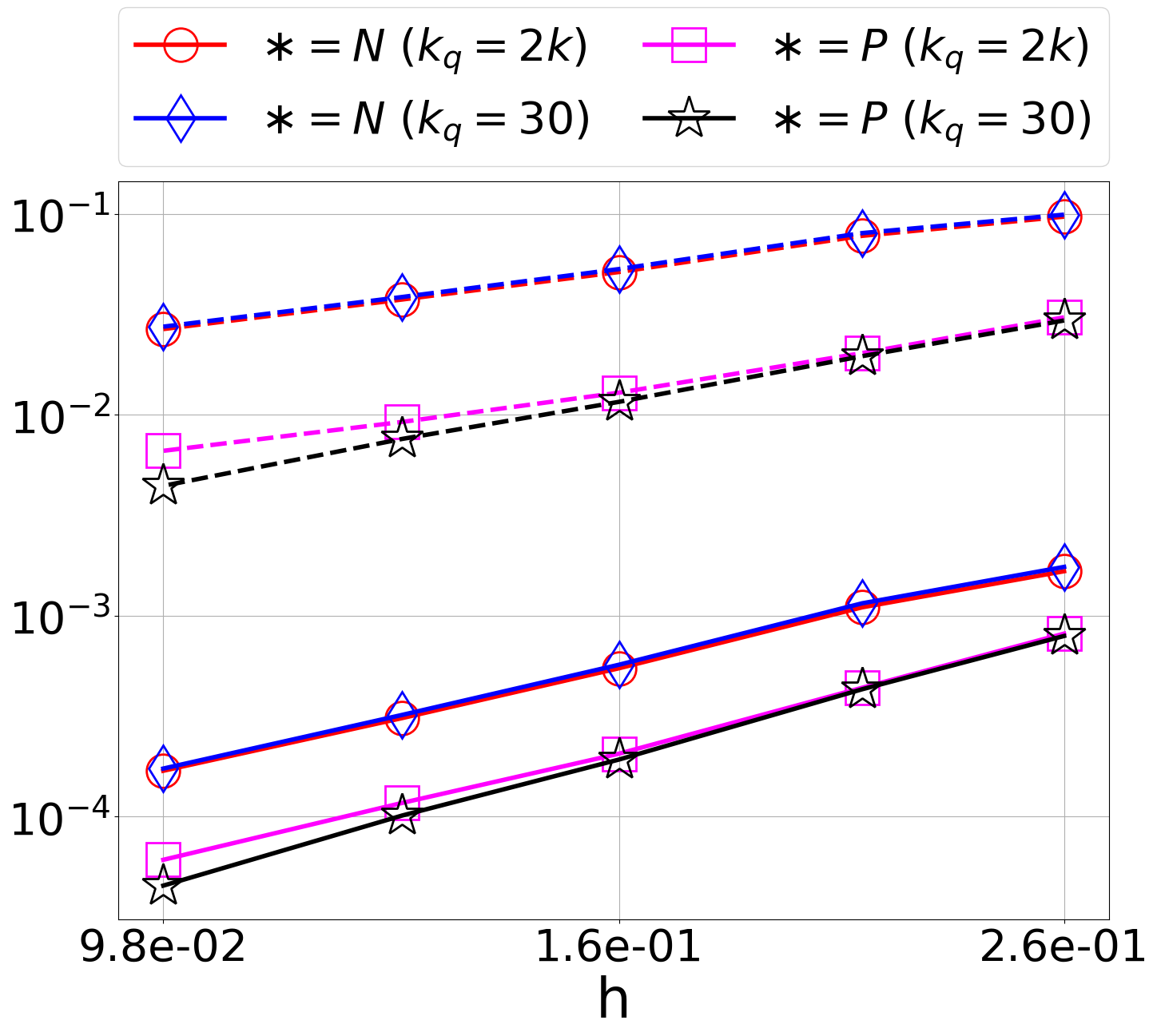}
        \caption{}
    \end{subfigure}
    \begin{subfigure}{0.3\textwidth}
        \includegraphics[width=1\textwidth]{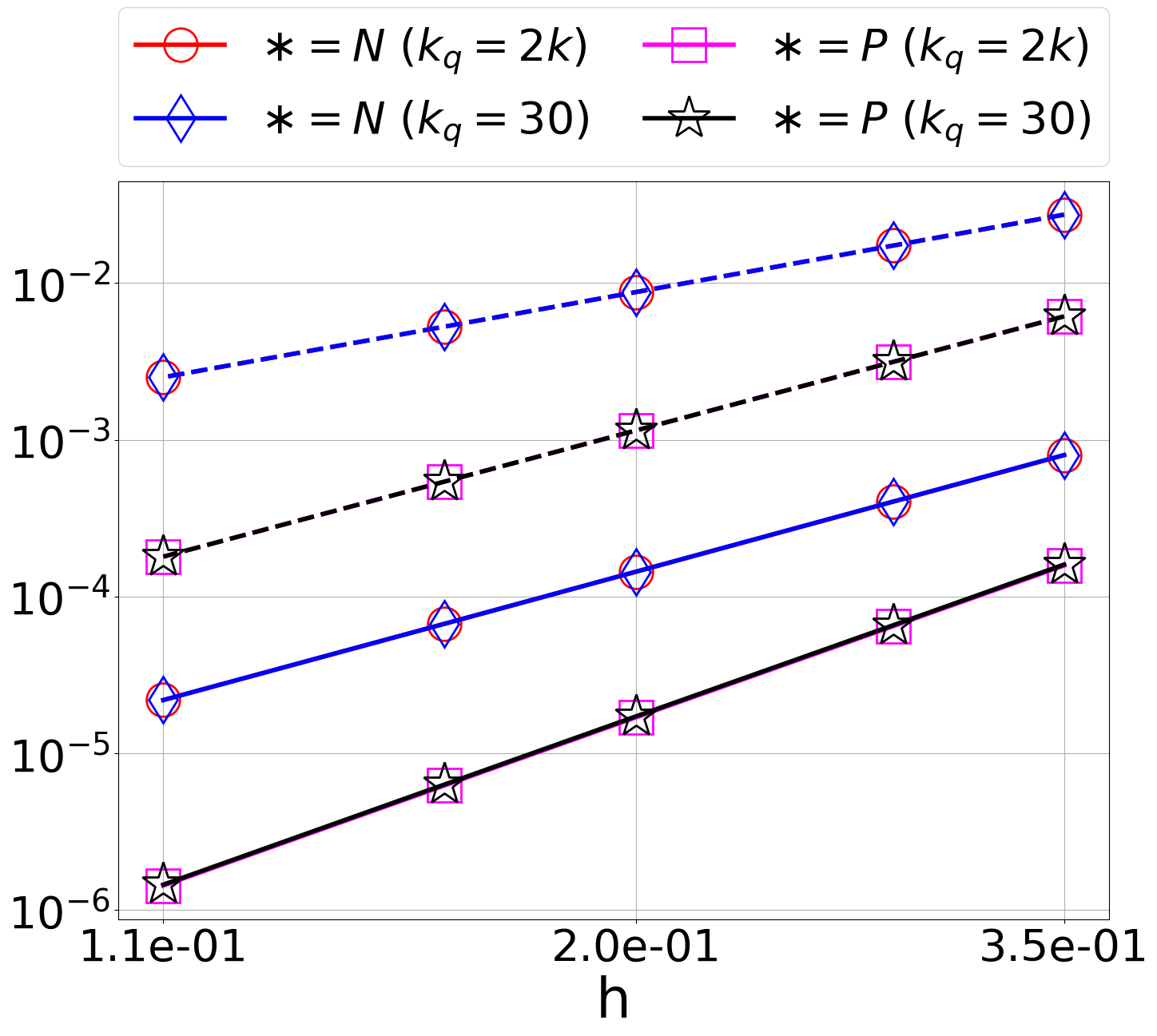}
        \caption{}
    \end{subfigure}
    \begin{subfigure}{0.3\textwidth}
        \includegraphics[width=1\textwidth]{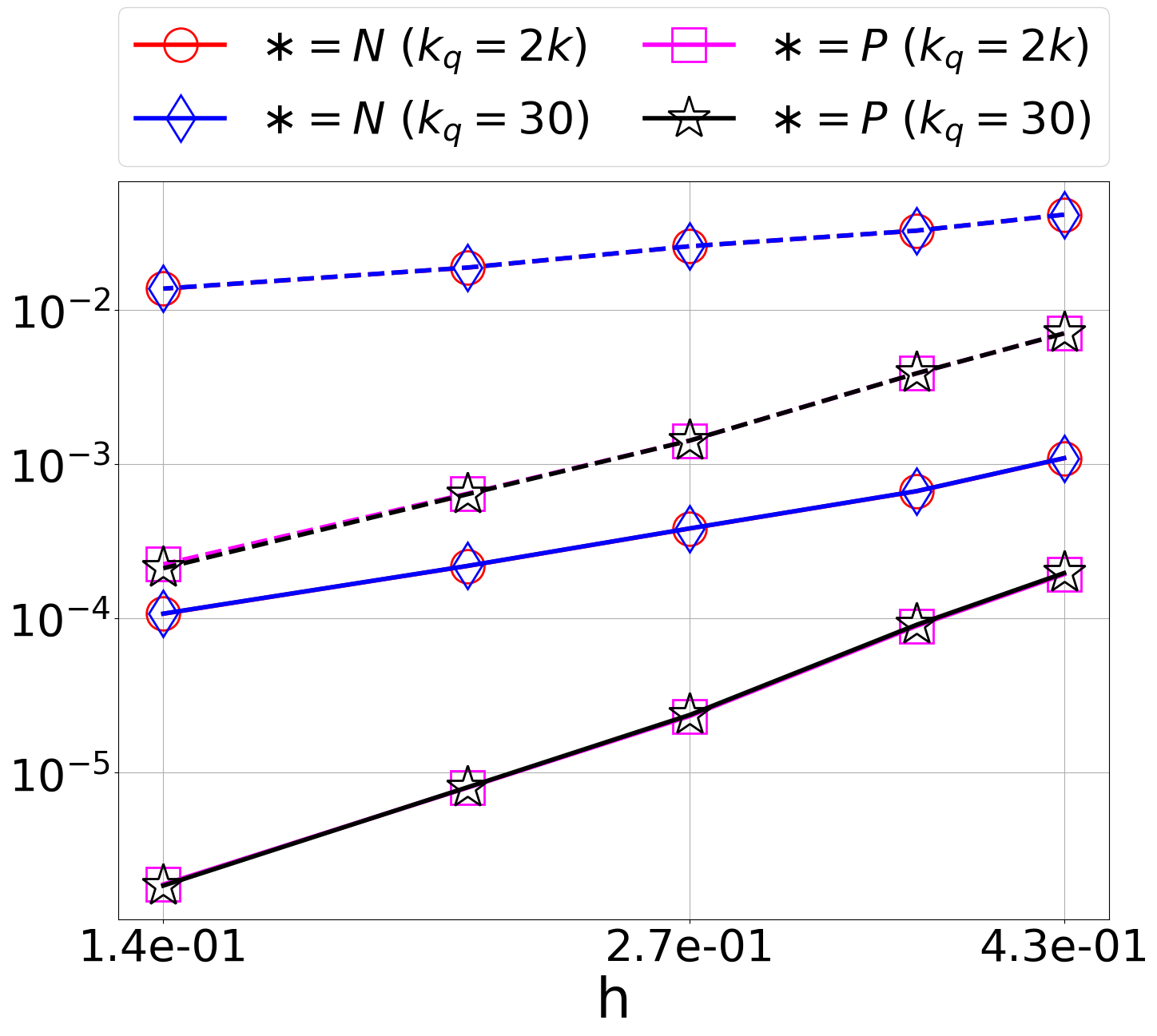}
        \caption{}
    \end{subfigure}
    \begin{subfigure}{0.3\textwidth}
        \includegraphics[width=1\textwidth]{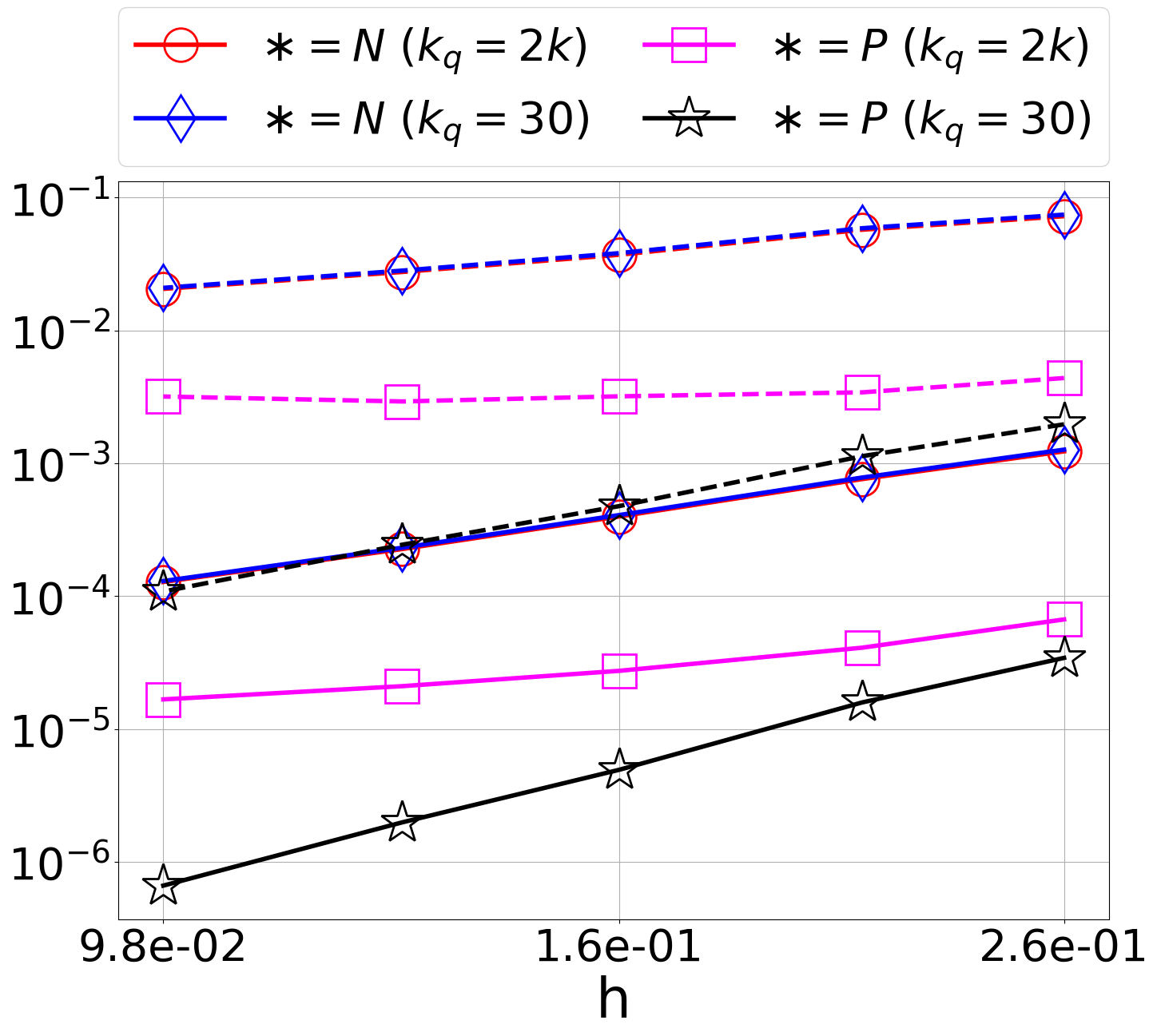}
        \caption{}
    \end{subfigure}
    \caption{Test 2: Behavior of the errors \eqref{eq:errors} as $h$ decreases for the function type \ref{function:A}. Solid lines represent error in $\leb{2}{}$-norm, while dashed lines denote the error in $\sob{1}{}$-seminorm. Left: \ref{mesh:A}. Center: \ref{mesh:B}. Right: \ref{mesh:C}. First row: $k=1$. Second row: $k=2$. Third row: $k=3$.}
    \label{fig:test2:fa_k1}
\end{figure}

In Figures \ref{fig:test2:fc_k1}, \ref{fig:test2:fb_k1}, and \ref{fig:test2:fa_k1}, we show the behavior of errors \eqref{eq:errors} as $h$ decreases for \ref{function:C}, \ref{function:B}, and \ref{function:A}, resp., whereas in Tables \ref{tab:test2:fc}, \ref{tab:test2:fb}, and \ref{tab:test2:fa}, we report the corresponding empirical convergence rates.

Let us first focus on the B-NAVEM performance, i.e. the basis functions of types \ref{function:C} and \ref{function:B}, whose related results are shown in Figures \ref{fig:test2:fc_k1} and \ref{fig:test2:fb_k1} and in Tables \ref{tab:test2:fc} and \ref{tab:test2:fb}. 
We observe that the two B-NAVEM variants exploiting different definitions of the $\boperator$ operators exhibit very similar behaviors. For $k=1$, the plain B-NAVEM method is able to guarantee convergence up to the neural network precision, by confirming numerical results obtained in \cite{PintoreTeora2026} and the theoretical estimates. This is particularly evident in the case of the Voronoi meshes \ref{mesh:C}, where the proposed algorithm is necessary to recover the optimal convergence rates and overcome neural network limitations. Moreover, we observe that, since the B-NAVEM basis functions do not oscillate excessively, the usage of the low-order quadrature rule is enough to preserve the accuracy of the method. 

For the higher-order cases, namely $k=2$ and $k=3$, the proposed algorithm becomes essential to recover the expected approximation properties. In our experiments, the neural networks employed to predict the higher-order basis functions share the same architecture as those used for the lowest-order case. Since higher-order basis functions exhibit significantly more complex features, this choice naturally leads to a lower approximation accuracy. We deliberately avoid increasing the network architecture or extending the training process, as our objective is not to maximize the prediction accuracy of the neural networks, but rather to assess the robustness of the proposed algorithm under challenging conditions. Consequently, the learned basis functions are far from spanning the desired VEM space, providing a stringent test for the proposed correction. Despite this limitation, the algorithm successfully restores polynomial completeness and recovers the optimal convergence rates, demonstrating that its effectiveness is largely independent of the quality of the neural-network approximation. Finally, we observe that, for high orders, in addition to having polynomials, a sufficiently high quadrature rule is needed to retrieve the right convergence rates. In particular, for meshes of type \ref{mesh:C} and for $k=3$, the convergence rates decay more slowly when the $\sob{1}{}$- error is close to $10^{-3}$, in accordance with the fact that the quadrature error on these meshes has the same order of magnitude (see Table \ref{tab:test1:bnavem}), thereby posing a limitation on the convergence of the method.

Let us now consider the solutions obtained with the basis functions of type \ref{function:A}. Results are provided in Figure~\ref{fig:test2:fa_k1} and in Table~\ref{tab:test2:fa}. As already discussed, for $k=1$ the local space already contains $\Poly{k}{E}$ and the basis functions have small gradients everywhere in $E$. This implies that {\postprocess} does not modify the local space (even if the basis functions may vary). As a consequence, all the curves in the first row of Figure \ref{fig:test2:fa_k1} are overlapped and aligned with that of the correct slope, as confirmed by the convergence rates reported in Table \ref{tab:test2:fa}. Moreover, we observe that the quadrature order $2k$ is already sufficient to properly integrate the functions without posing limitation to the convergence of the methods.  When $k=2$, instead, the post-processing procedure is superfluous only on the Cartesian meshes \ref{mesh:A}, whereas it is needed to retrieve the right order of convergence on \ref{mesh:B} and \ref{mesh:C}. The first effects of the loss of convergence rate due to the numerical quadrature begin to appear for the Voronoi meshes \ref{mesh:C}. Similar phenomena are visible and more accentuated in the last row, for $k=3$, where the post-processing is also required for the case of the Cartesian meshes.

\begin{figure}[!ht]
    \centering
    \begin{subfigure}{0.3\textwidth}
        \includegraphics[width=1\textwidth]{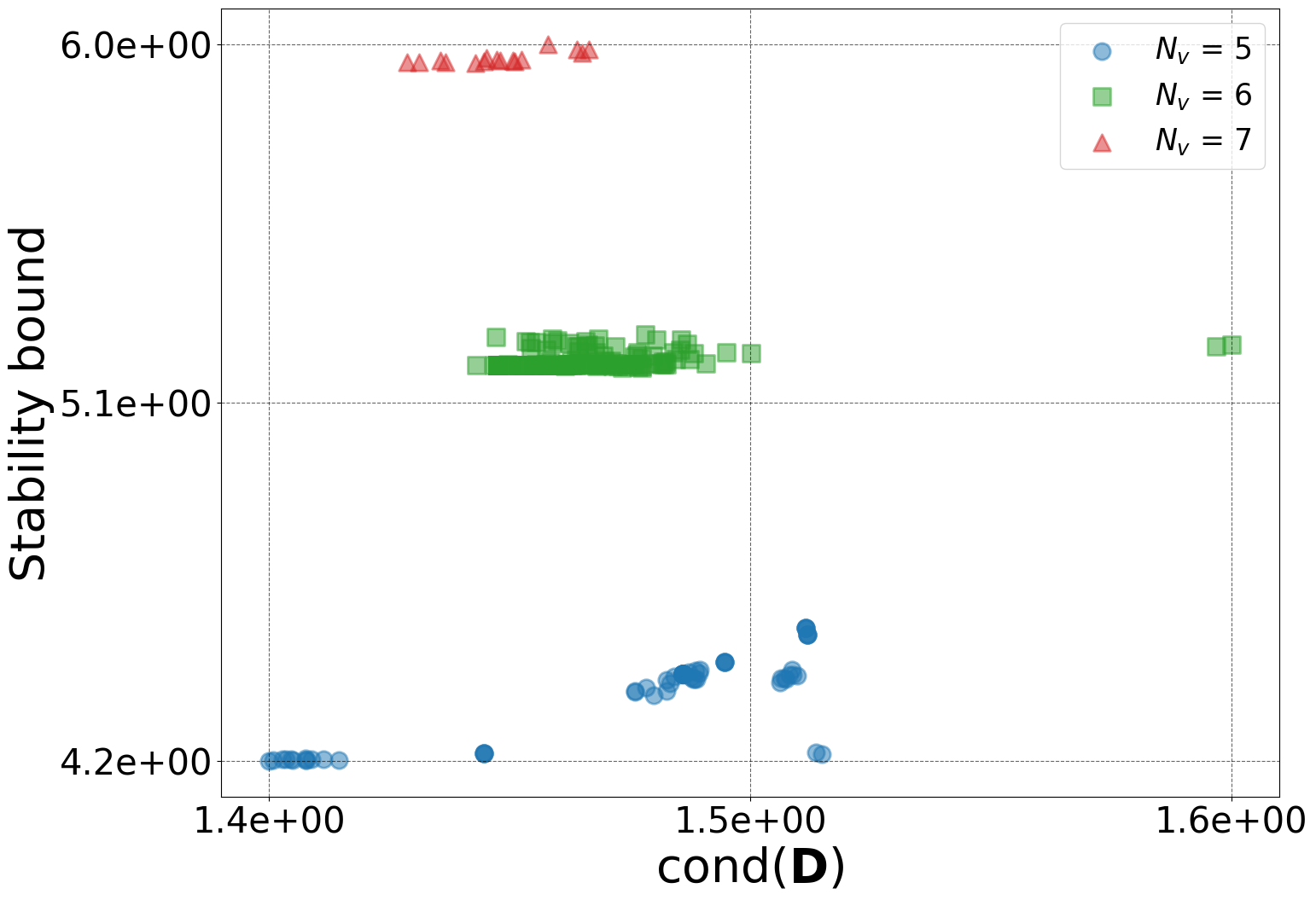}
        \caption{}
    \end{subfigure}
    \begin{subfigure}{0.3\textwidth}
        \includegraphics[width=1\textwidth]{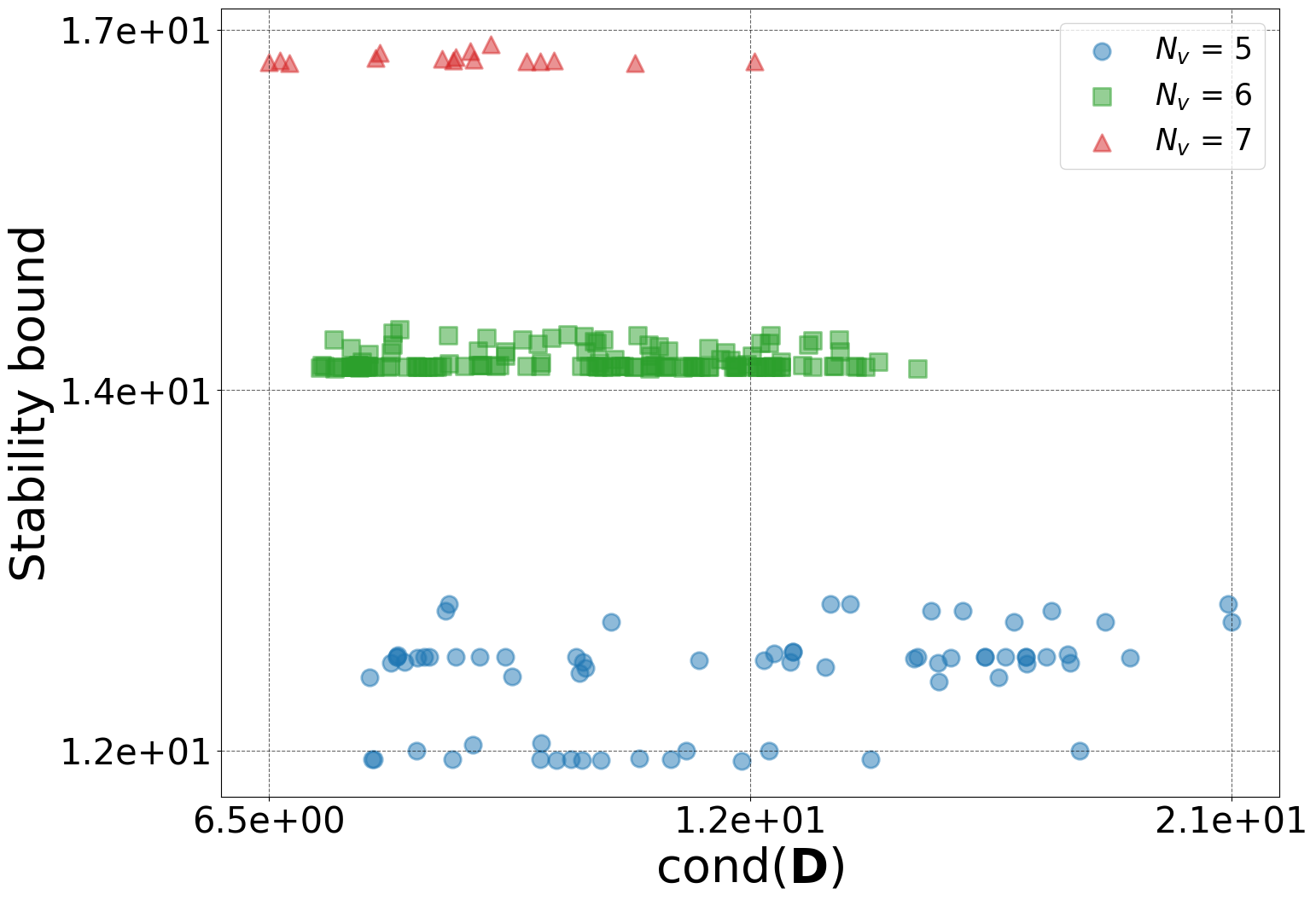}
        \caption{}
    \end{subfigure}
    \begin{subfigure}{0.3\textwidth}
        \includegraphics[width=1\textwidth]{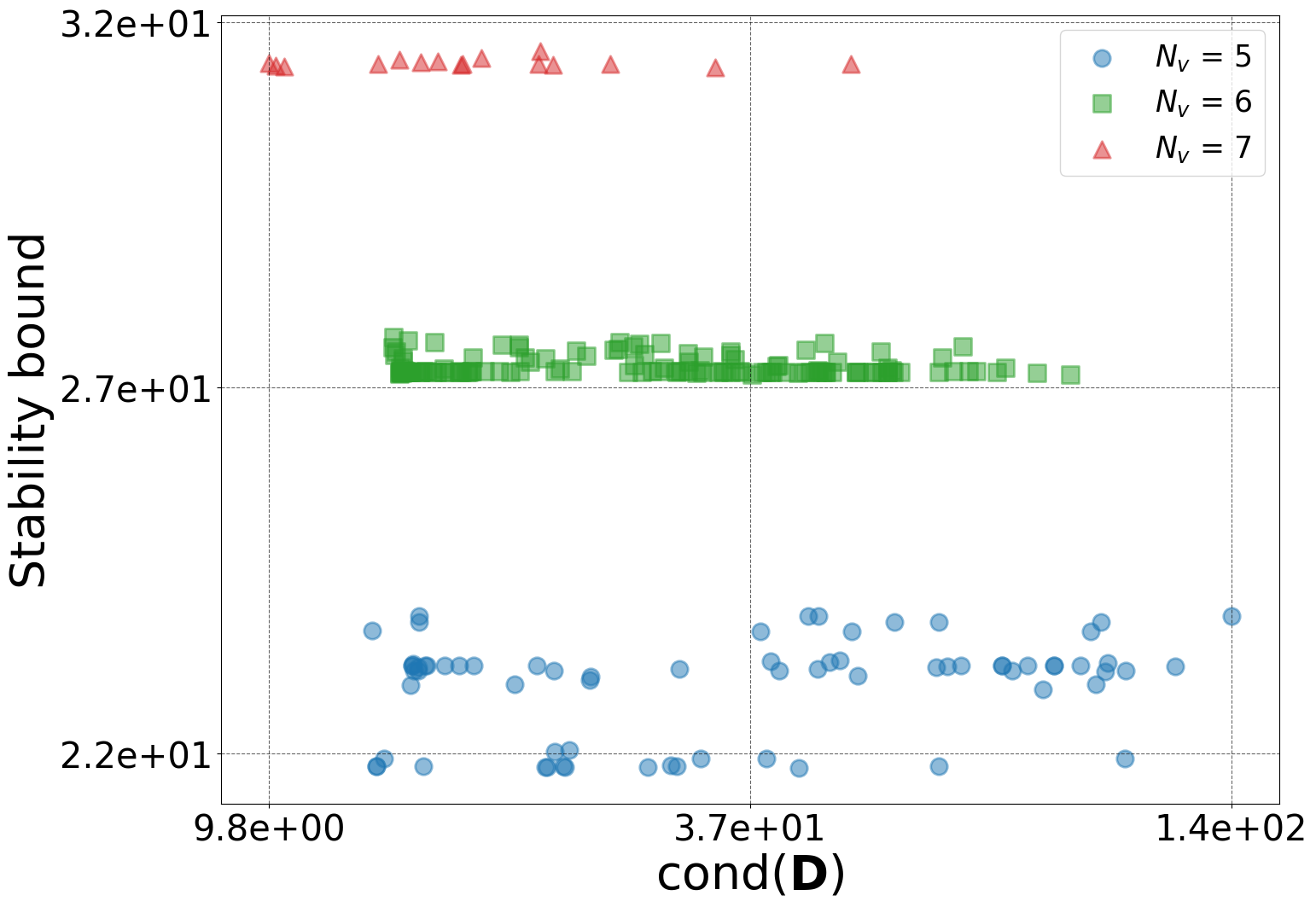}
        \caption{}
    \end{subfigure}
    \begin{subfigure}{0.3\textwidth}
        \includegraphics[width=1\textwidth]{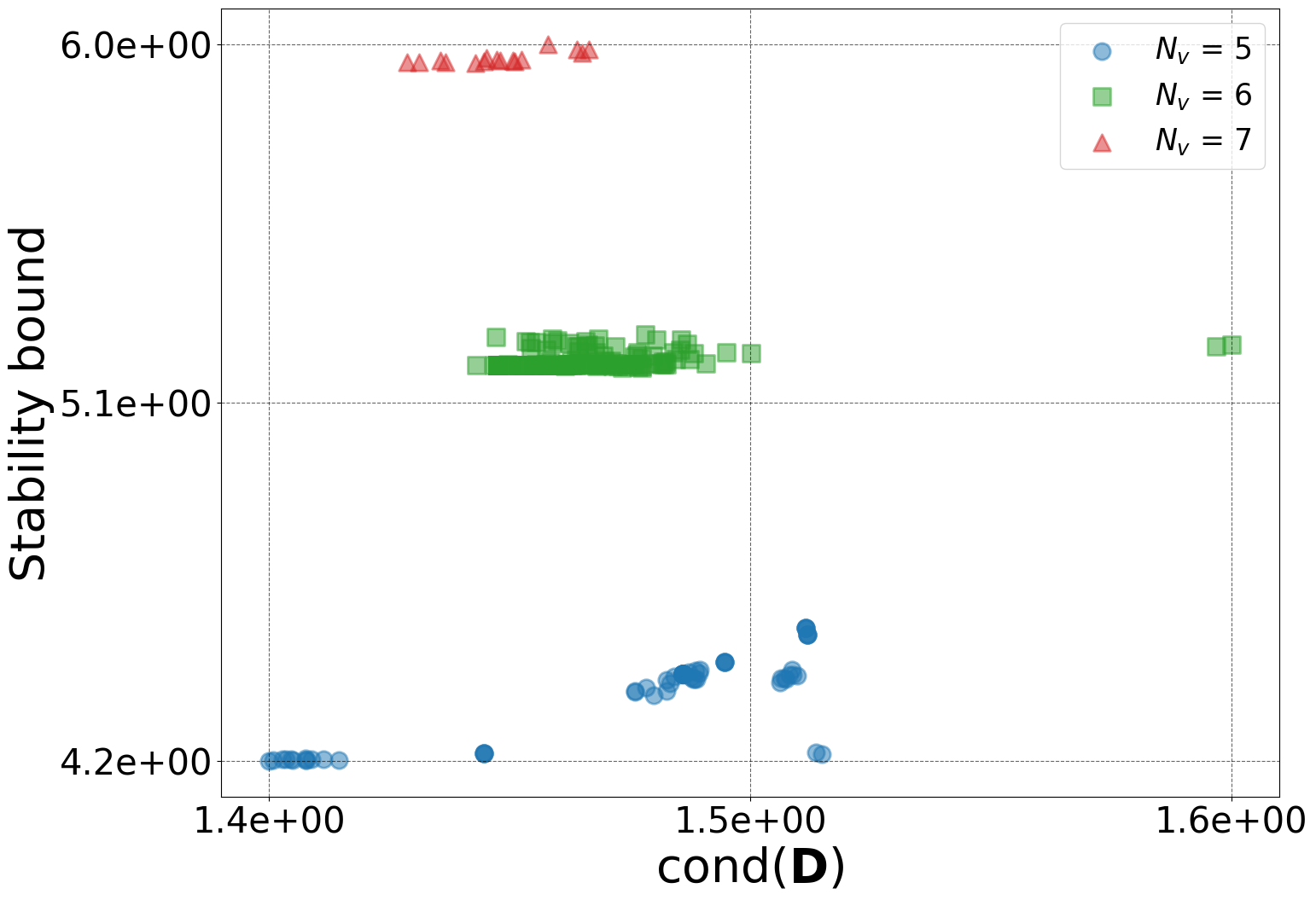}
        \caption{}
    \end{subfigure}
    \begin{subfigure}{0.3\textwidth}
        \includegraphics[width=1\textwidth]{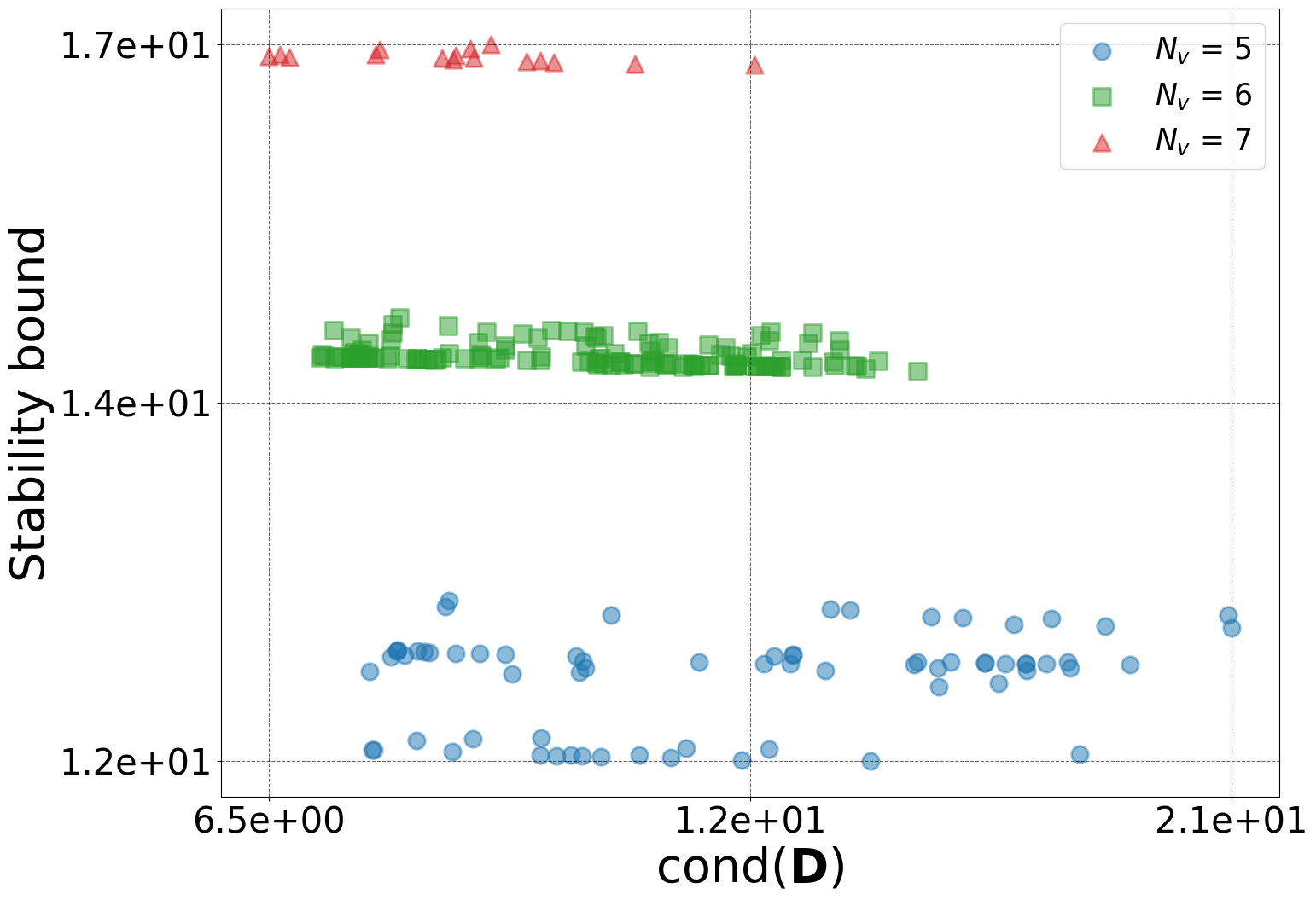}
        \caption{}
    \end{subfigure}
    \begin{subfigure}{0.3\textwidth}
        \includegraphics[width=1\textwidth]{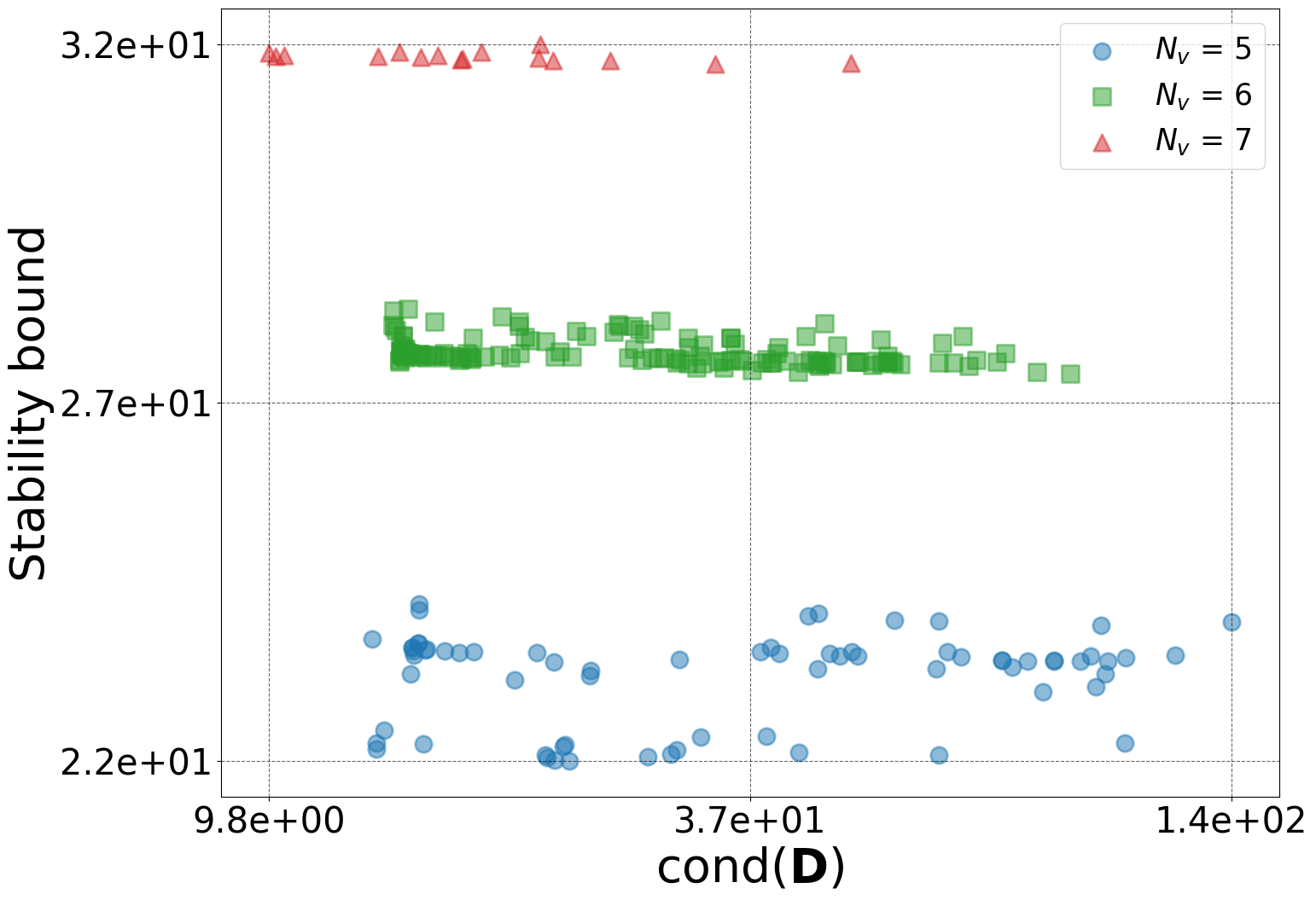}
        \caption{}
    \end{subfigure}
    \caption{Test 2: Behavior of the stability bounds \eqref{eq:function_regularity} (first row) and \eqref{eq:function_regularity_post} (second row) with respect to the condition number of the matrix $\mathbf{D}$ defined in \eqref{eq:d_matrix} for the set of elements related to the last refinement of the Voronoi mesh \ref{mesh:C}. Different colors represent different polygon classes. Each column represents a different order $k=1,2,3$ from left to right. Function type \ref{function:A}.}
    \label{fig:test2:stab_bound_ft1}
\end{figure}

\begin{figure}[!ht]
    \centering
    \begin{subfigure}{0.3\textwidth}
        \includegraphics[width=1\textwidth]{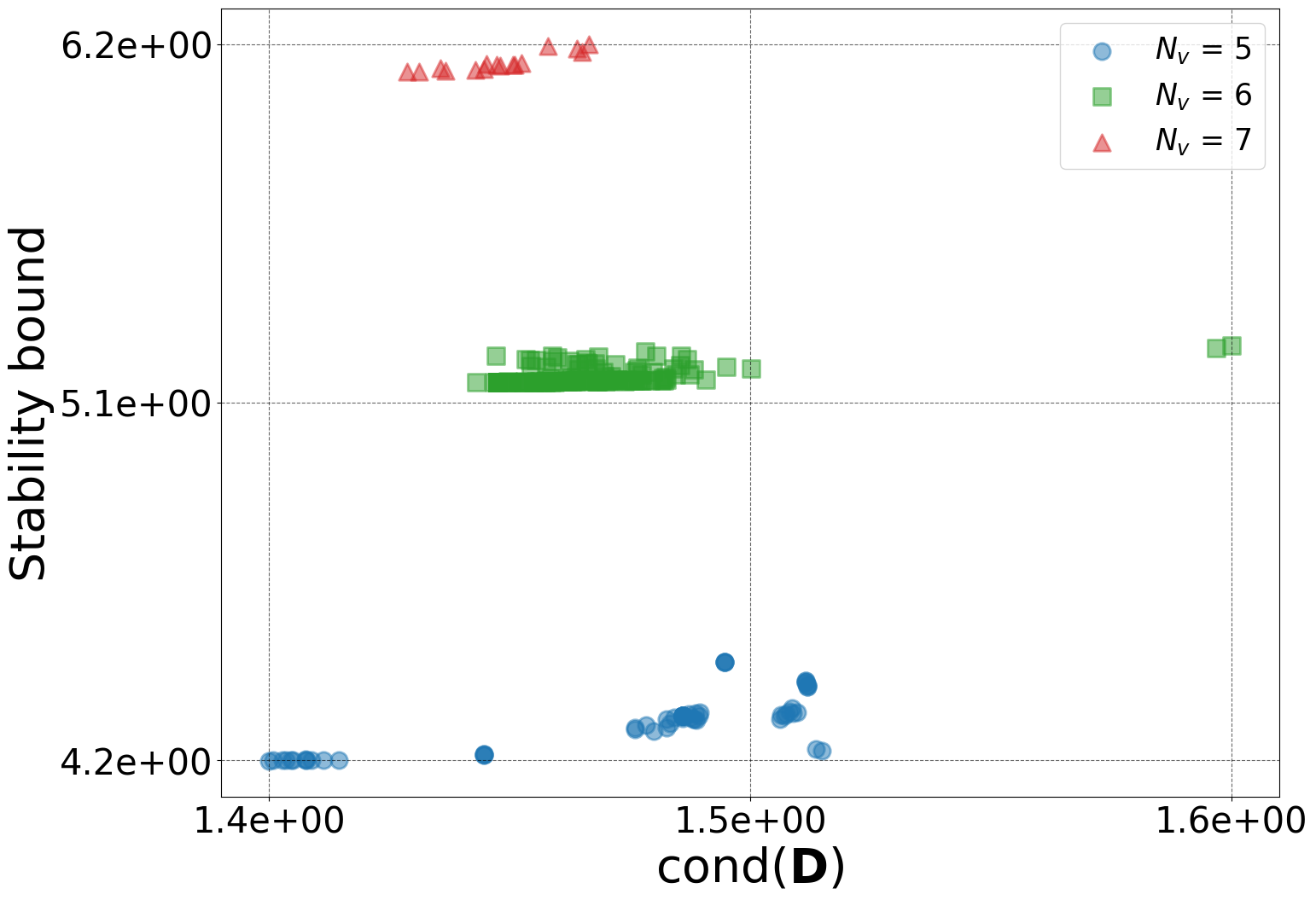}
        \caption{}
    \end{subfigure}
    \begin{subfigure}{0.3\textwidth}
        \includegraphics[width=1\textwidth]{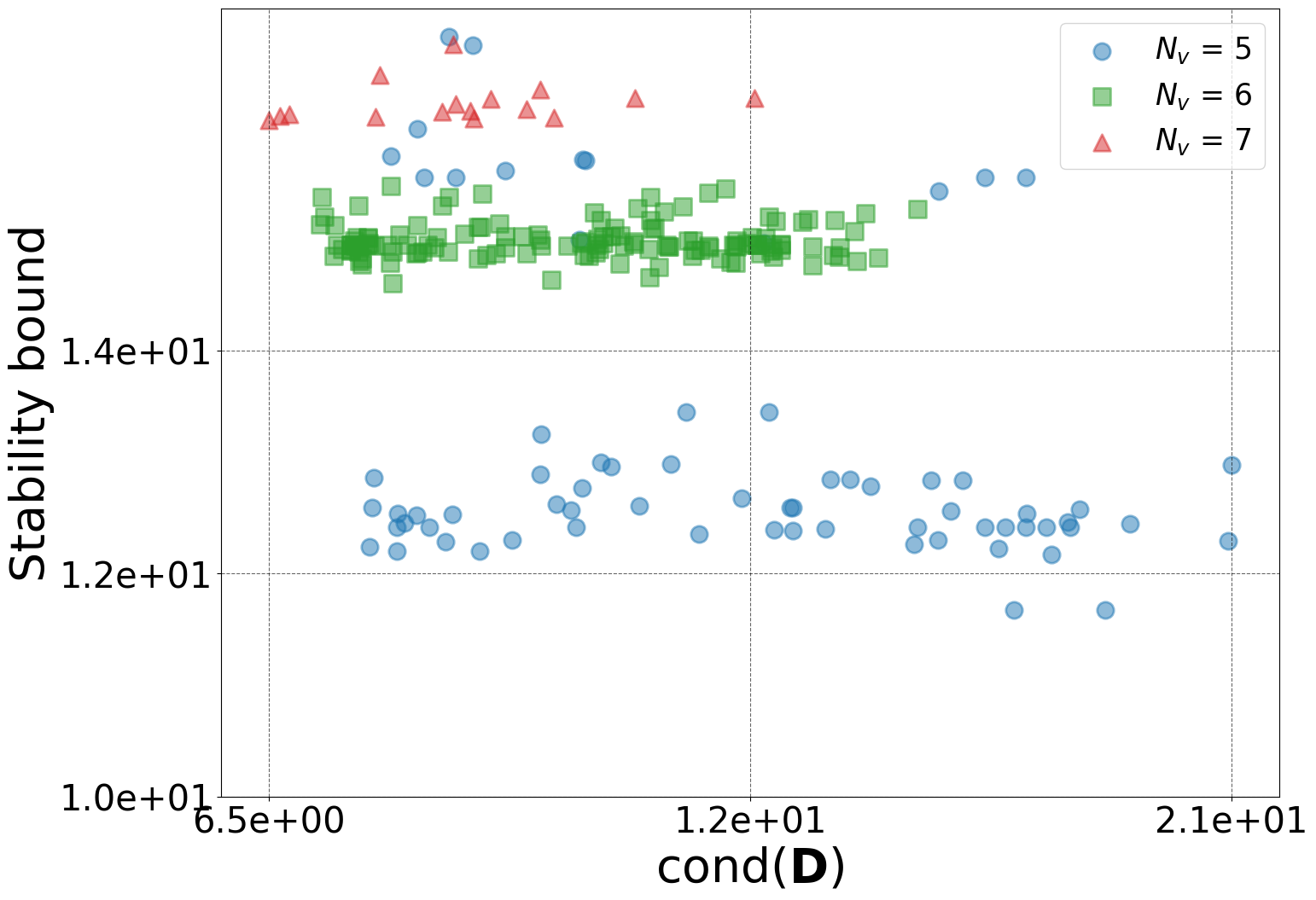}
        \caption{}
    \end{subfigure}
    \begin{subfigure}{0.3\textwidth}
        \includegraphics[width=1\textwidth]{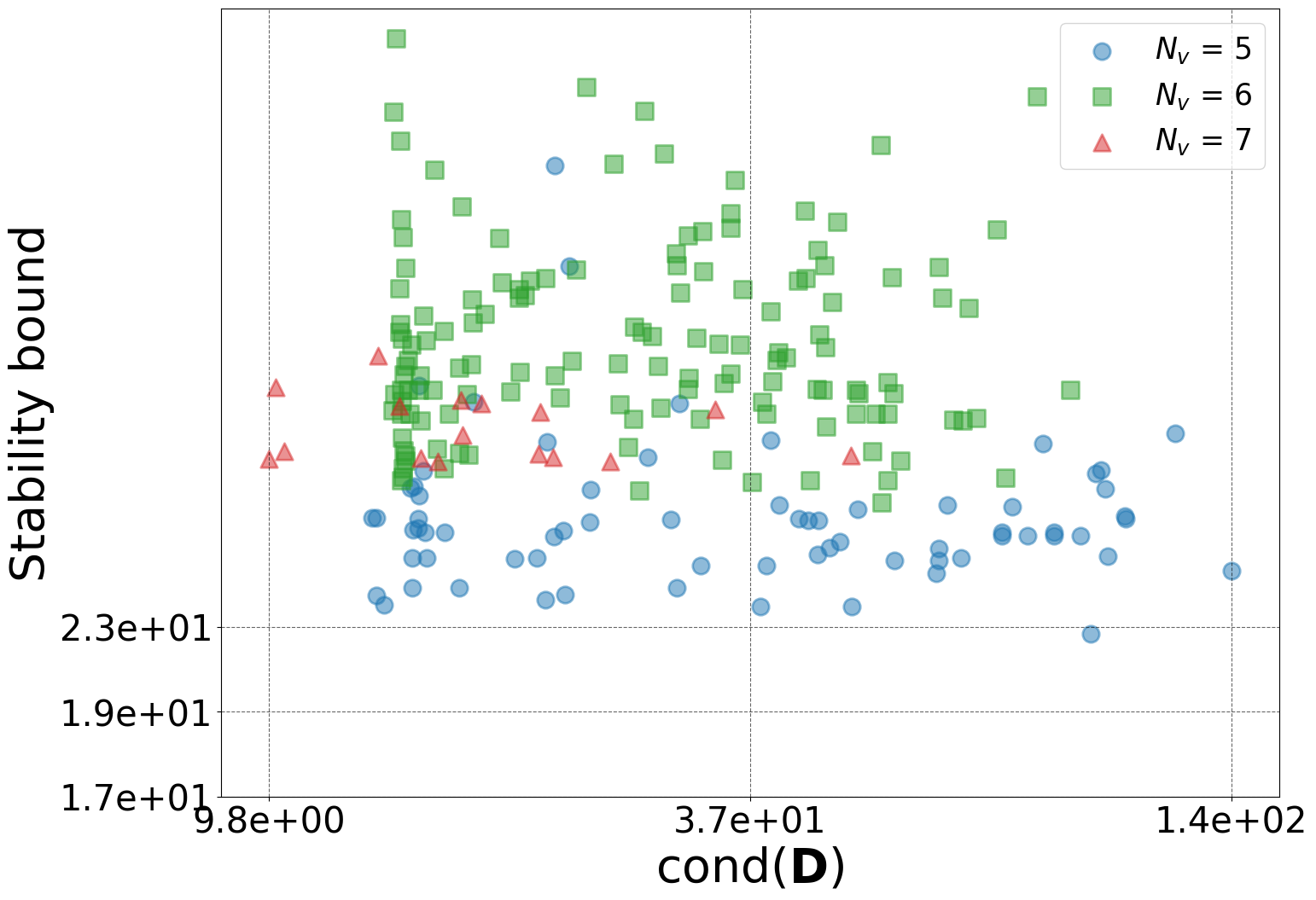}
        \caption{}
    \end{subfigure}
    \begin{subfigure}{0.3\textwidth}
        \includegraphics[width=1\textwidth]{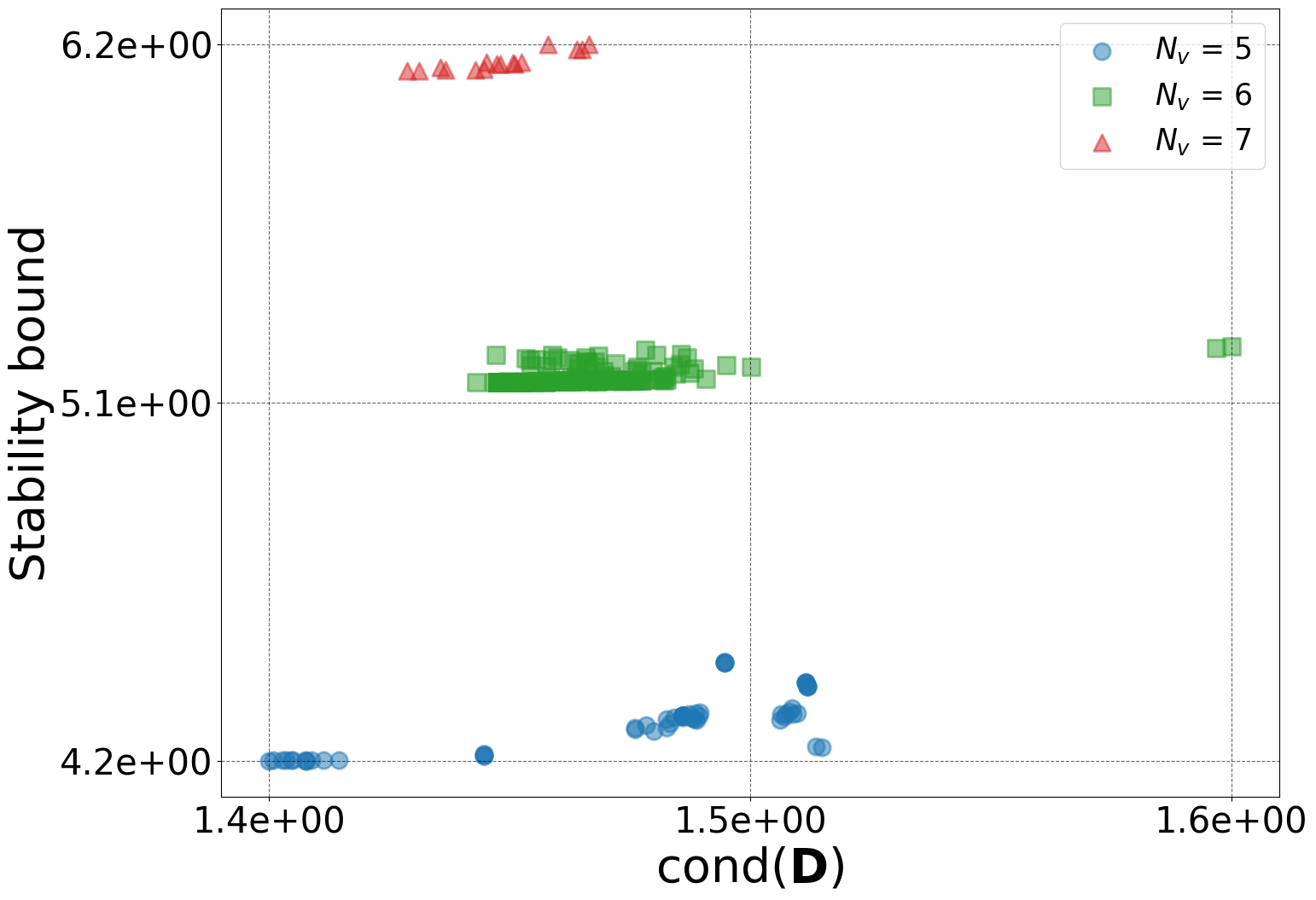}
        \caption{}
    \end{subfigure}
    \begin{subfigure}{0.3\textwidth}
        \includegraphics[width=1\textwidth]{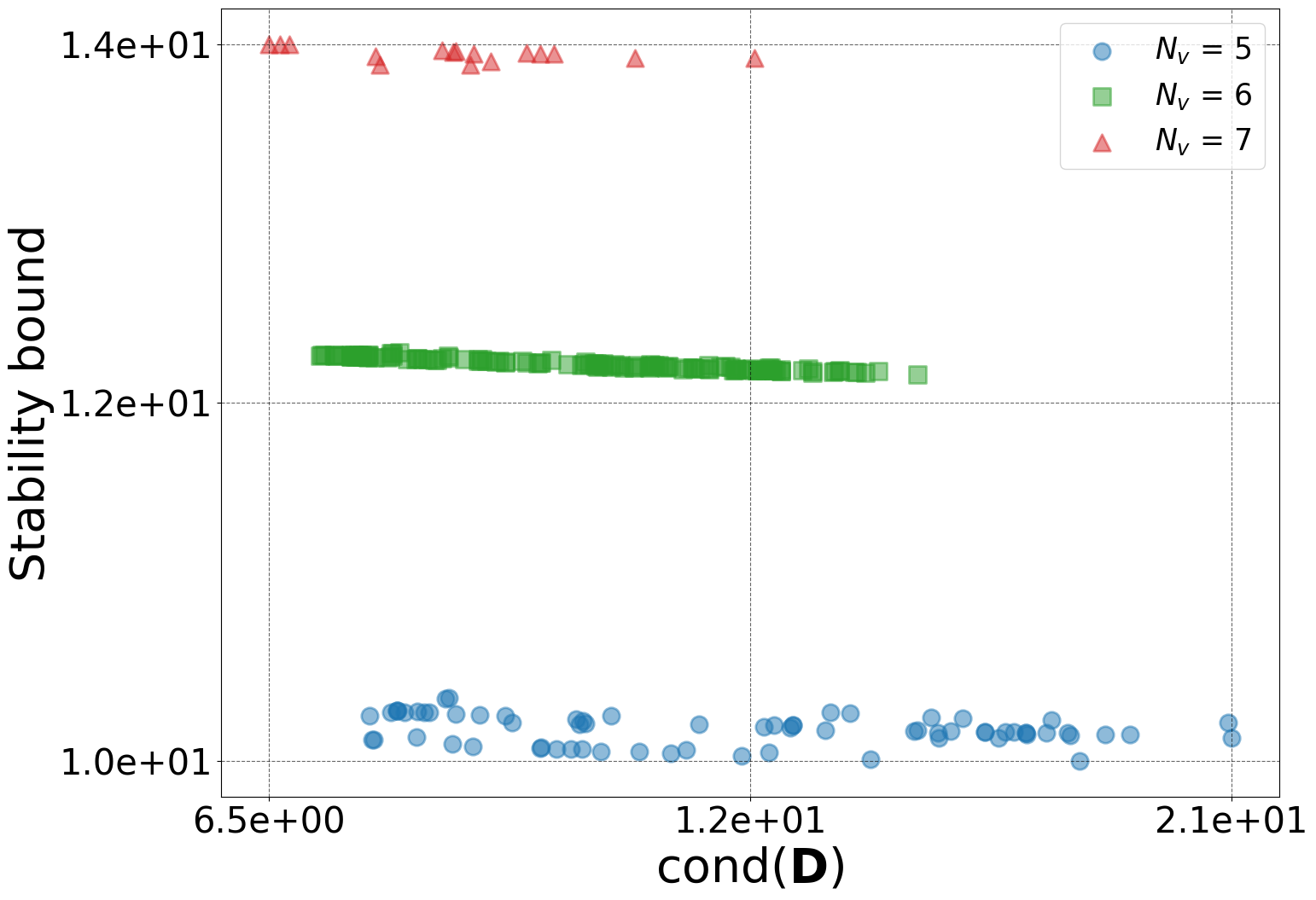}
        \caption{}
    \end{subfigure}
    \begin{subfigure}{0.3\textwidth}
        \includegraphics[width=1\textwidth]{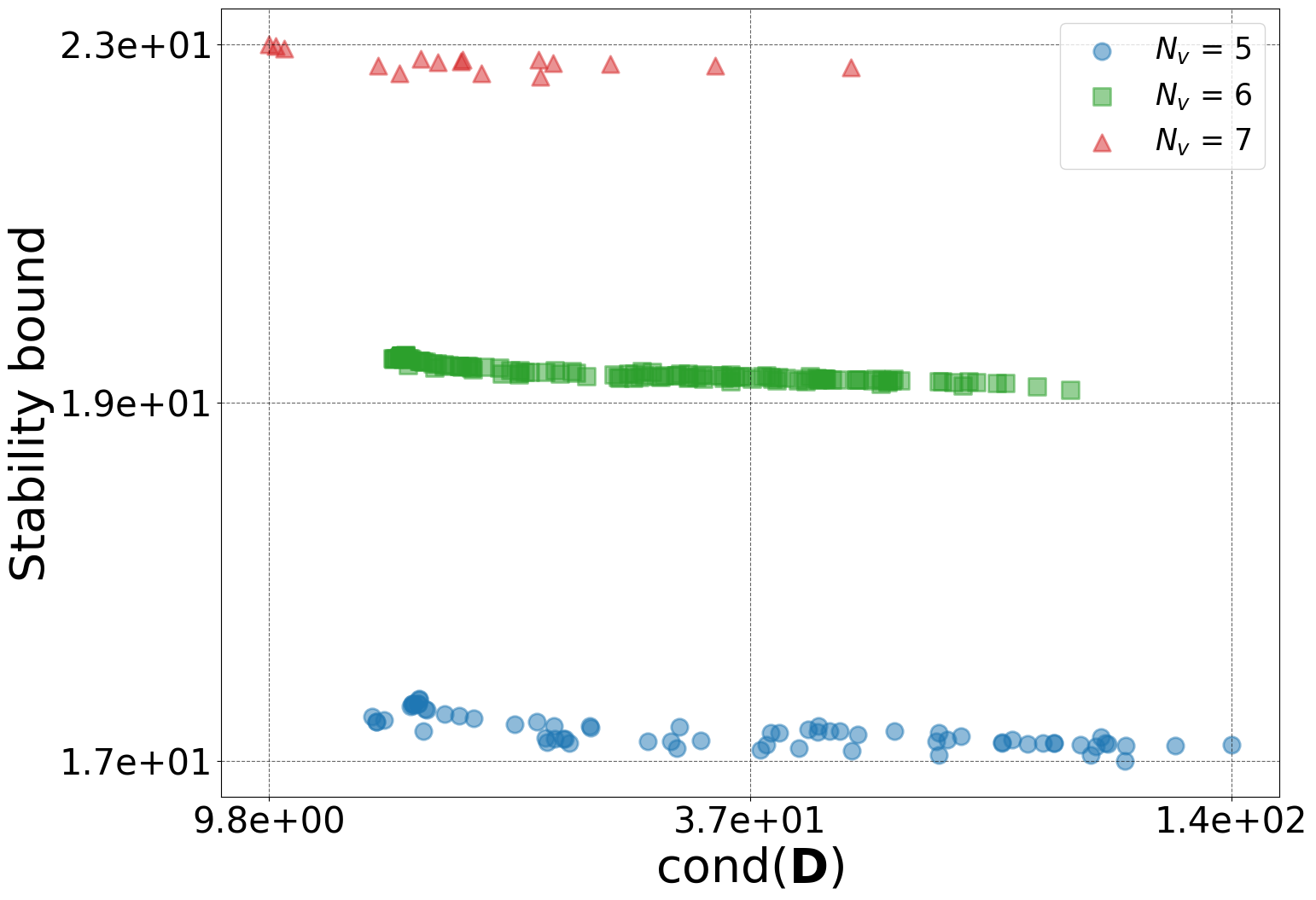}
        \caption{}
    \end{subfigure}
    \caption{Test 2: Behavior of the stability bounds \eqref{eq:function_regularity} (first row) and \eqref{eq:function_regularity_post} (second row) with respect to the condition number of the matrix $\mathbf{D}$ defined in \eqref{eq:d_matrix} for the set of elements related to the last refinement of the Voronoi mesh \ref{mesh:C}. Different colors represent different polygon classes. Each column represents a different order $k=1,2,3$ from left to right. Function type \ref{function:B}.}
    \label{fig:test2:stab_bound_ft2}
\end{figure}

\begin{figure}[!ht]
    \centering
    \begin{subfigure}{0.3\textwidth}
        \includegraphics[width=1\textwidth]{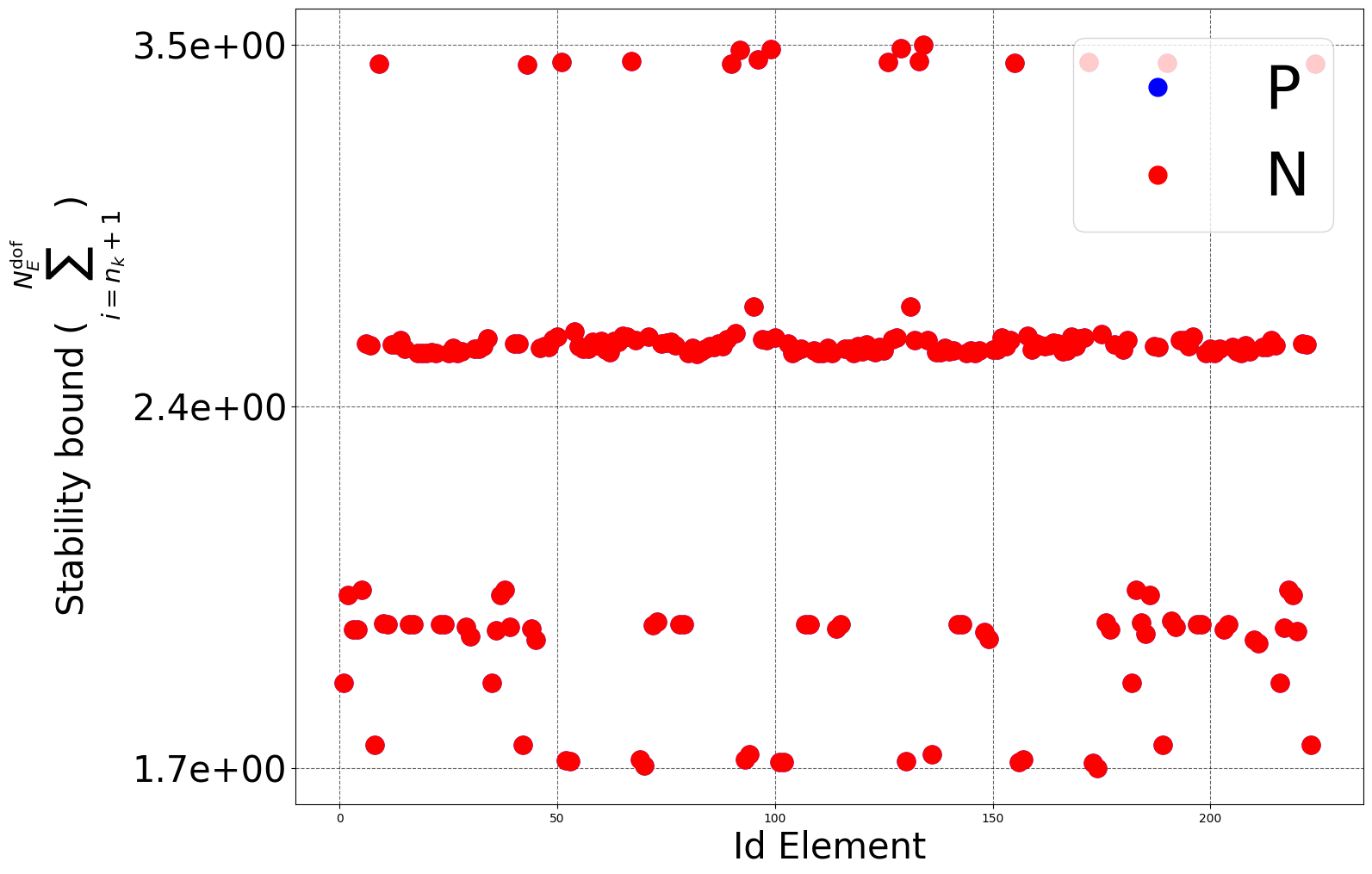}
        \caption{}
    \end{subfigure}
    \begin{subfigure}{0.3\textwidth}
        \includegraphics[width=1\textwidth]{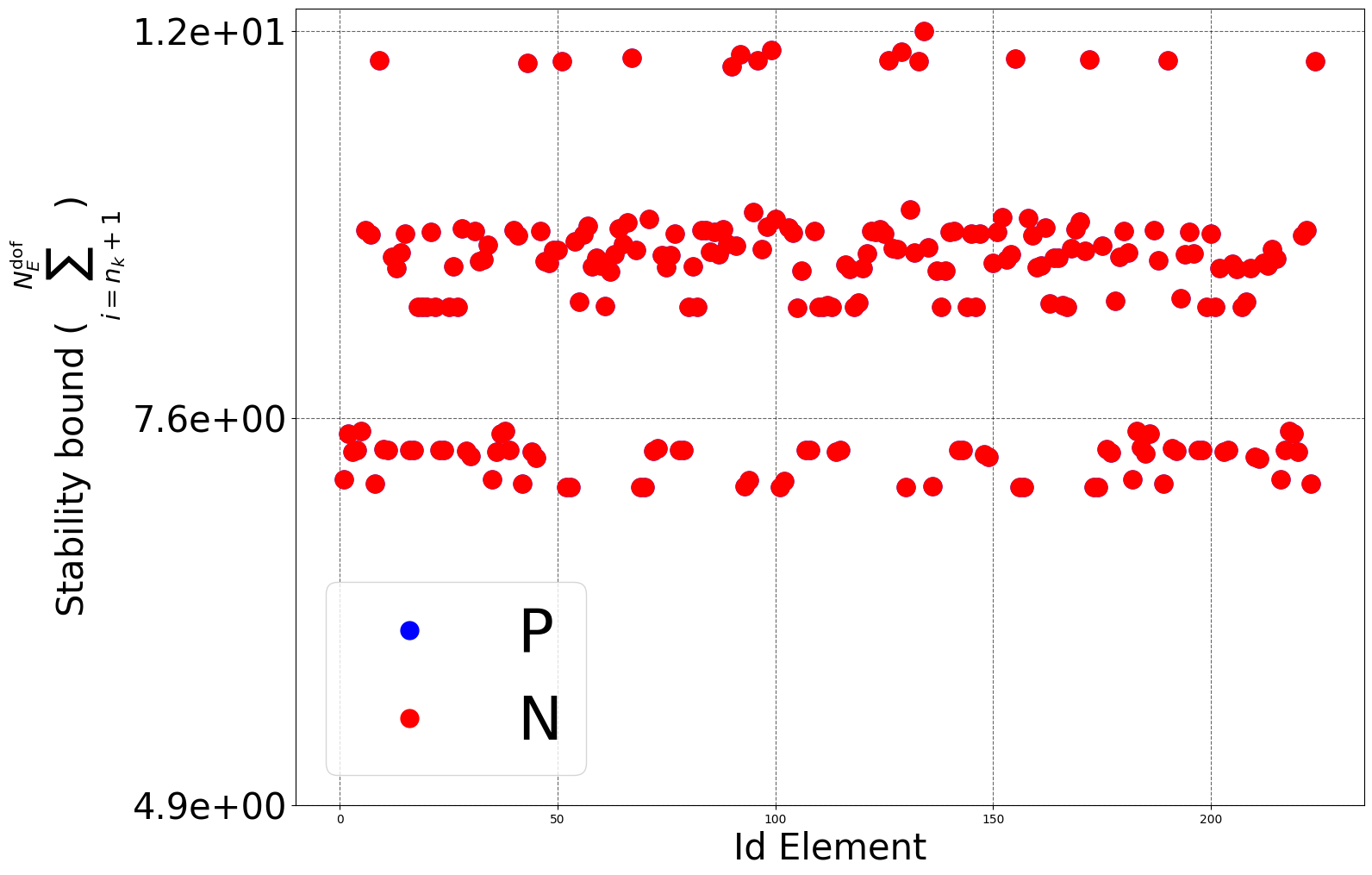}
        \caption{}
    \end{subfigure}
    \begin{subfigure}{0.3\textwidth}
        \includegraphics[width=1\textwidth]{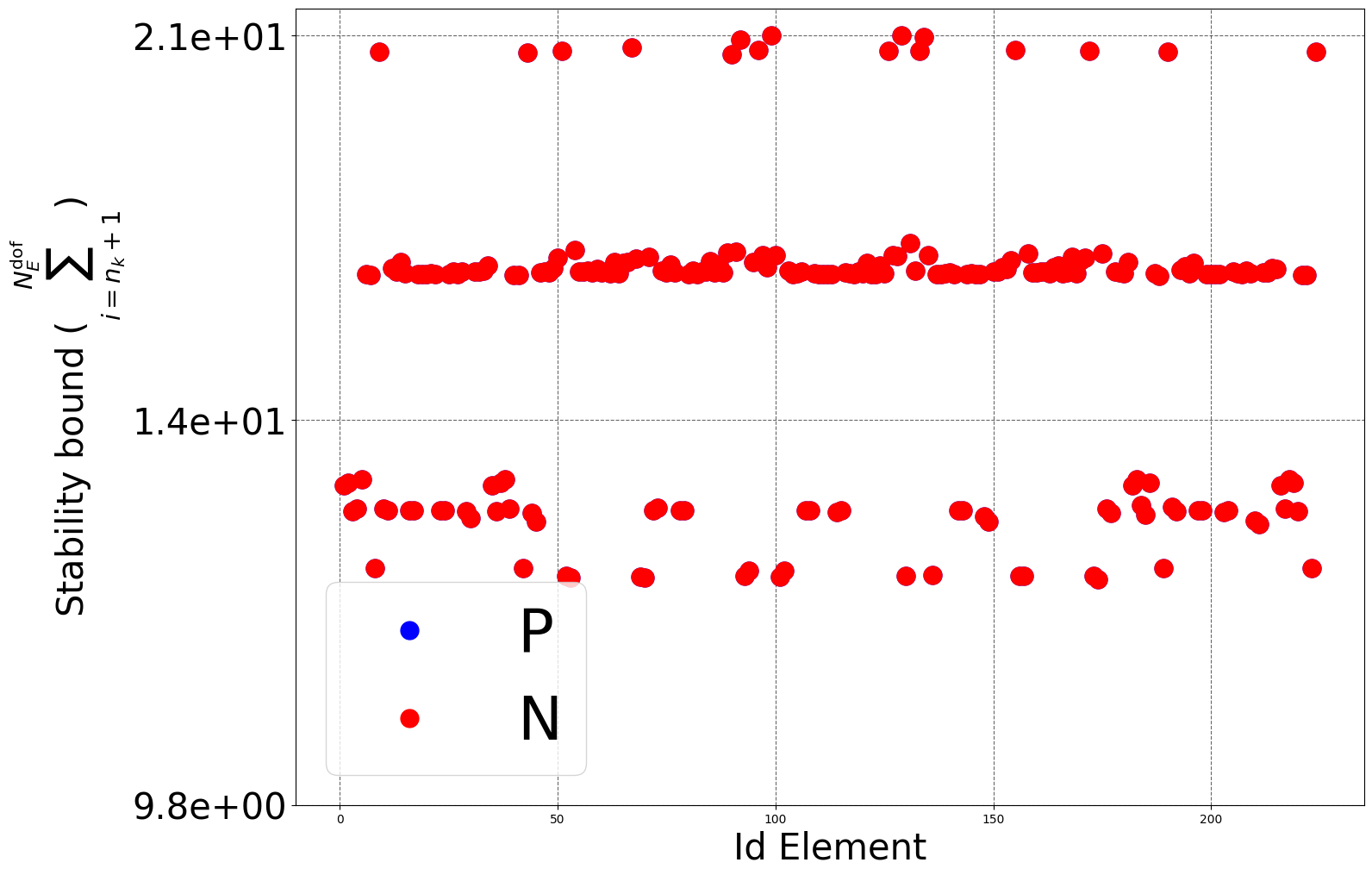}
        \caption{}
    \end{subfigure}
    \begin{subfigure}{0.3\textwidth}
        \includegraphics[width=1\textwidth]{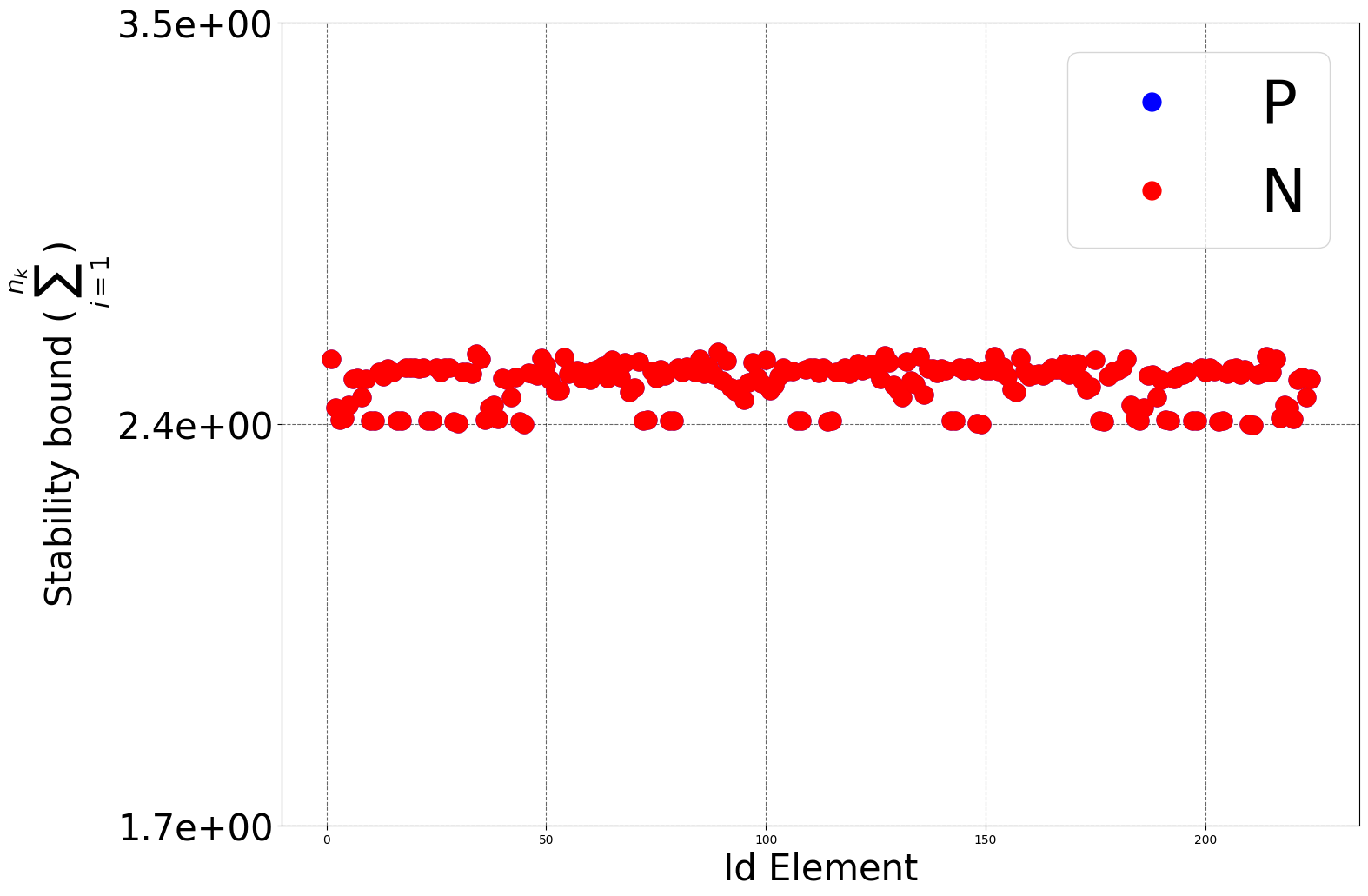}
        \caption{}
    \end{subfigure}
    \begin{subfigure}{0.3\textwidth}
        \includegraphics[width=1\textwidth]{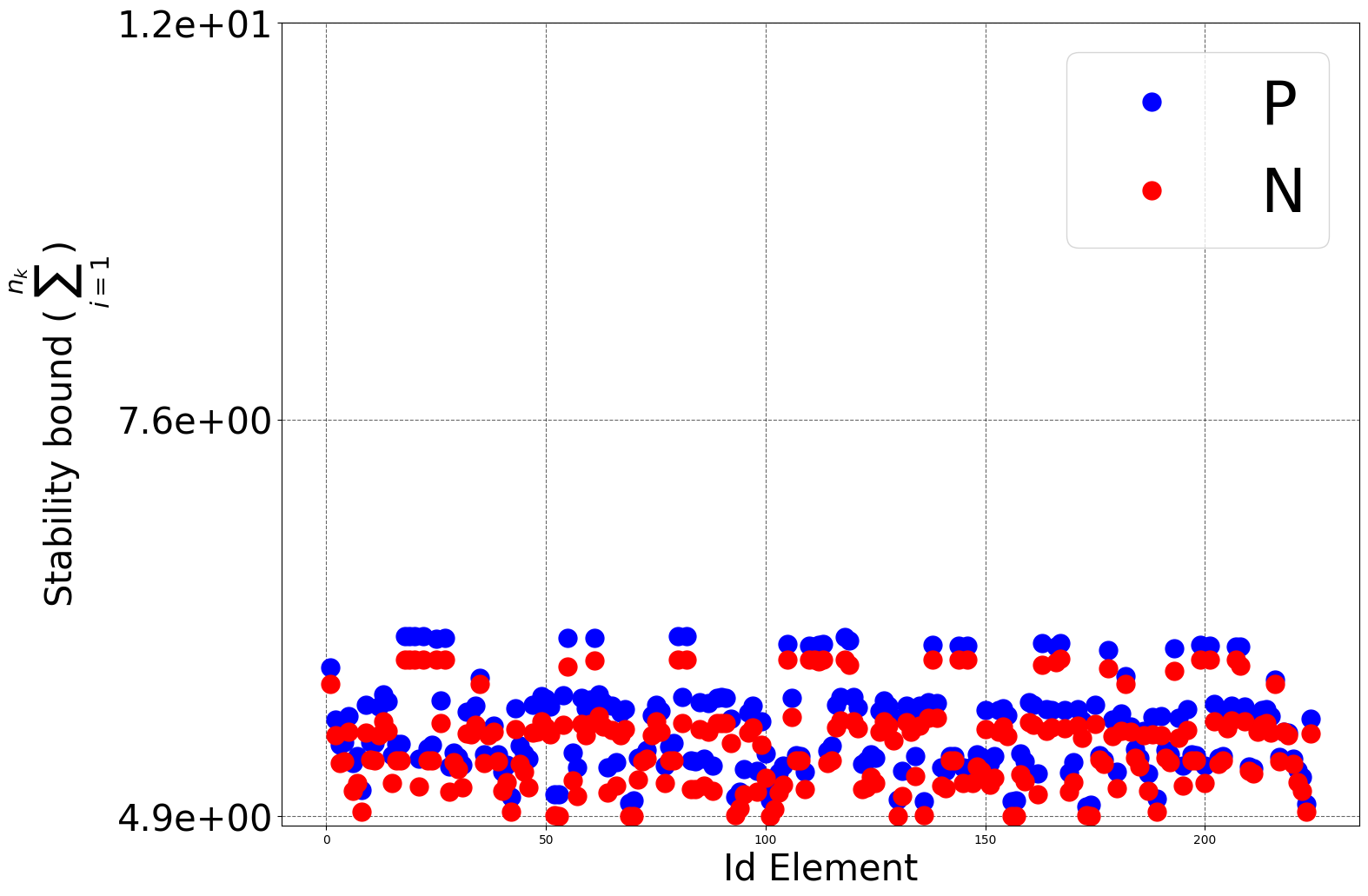}
        \caption{}
    \end{subfigure}
    \begin{subfigure}{0.3\textwidth}
        \includegraphics[width=1\textwidth]{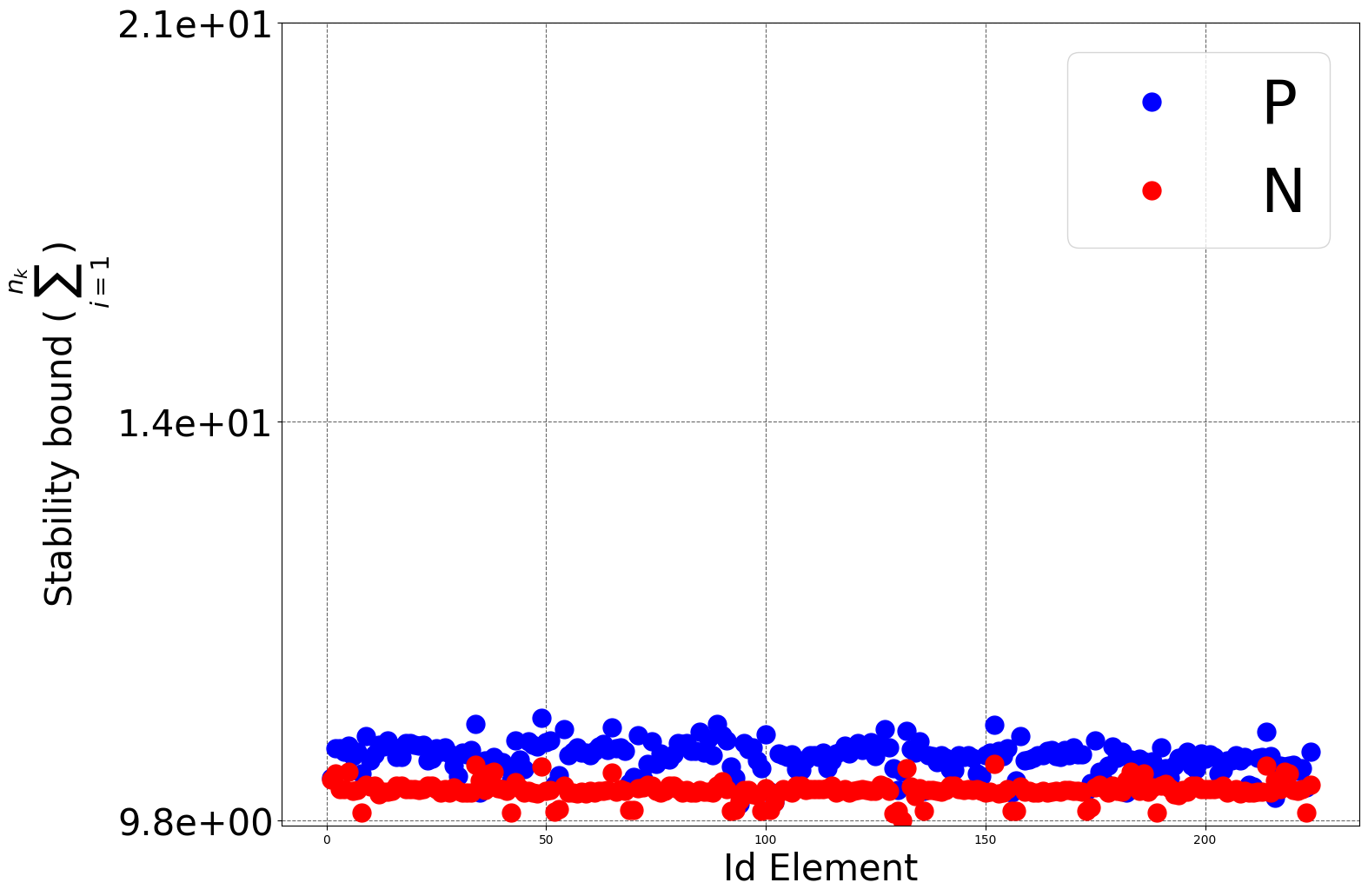}
        \caption{}
    \end{subfigure}
    \caption{Test 2: Behavior of the partial stability bounds \eqref{eq:partial_stability_last} (first row) and \eqref{eq:partial_stability_first} (second row) for the elements in the last refinement of the Voronoi mesh \ref{mesh:C}. The bounds are associated with the original functions (red) and with the post-processed functions (blue). Each column represents a different order $k=1,2,3$ from left to right. Function type \ref{function:A}.}
    \label{fig:test2:stab_bound_ft1_split}
\end{figure}

\begin{figure}[!ht]
    \centering
    \begin{subfigure}{0.3\textwidth}
        \includegraphics[width=1\textwidth]{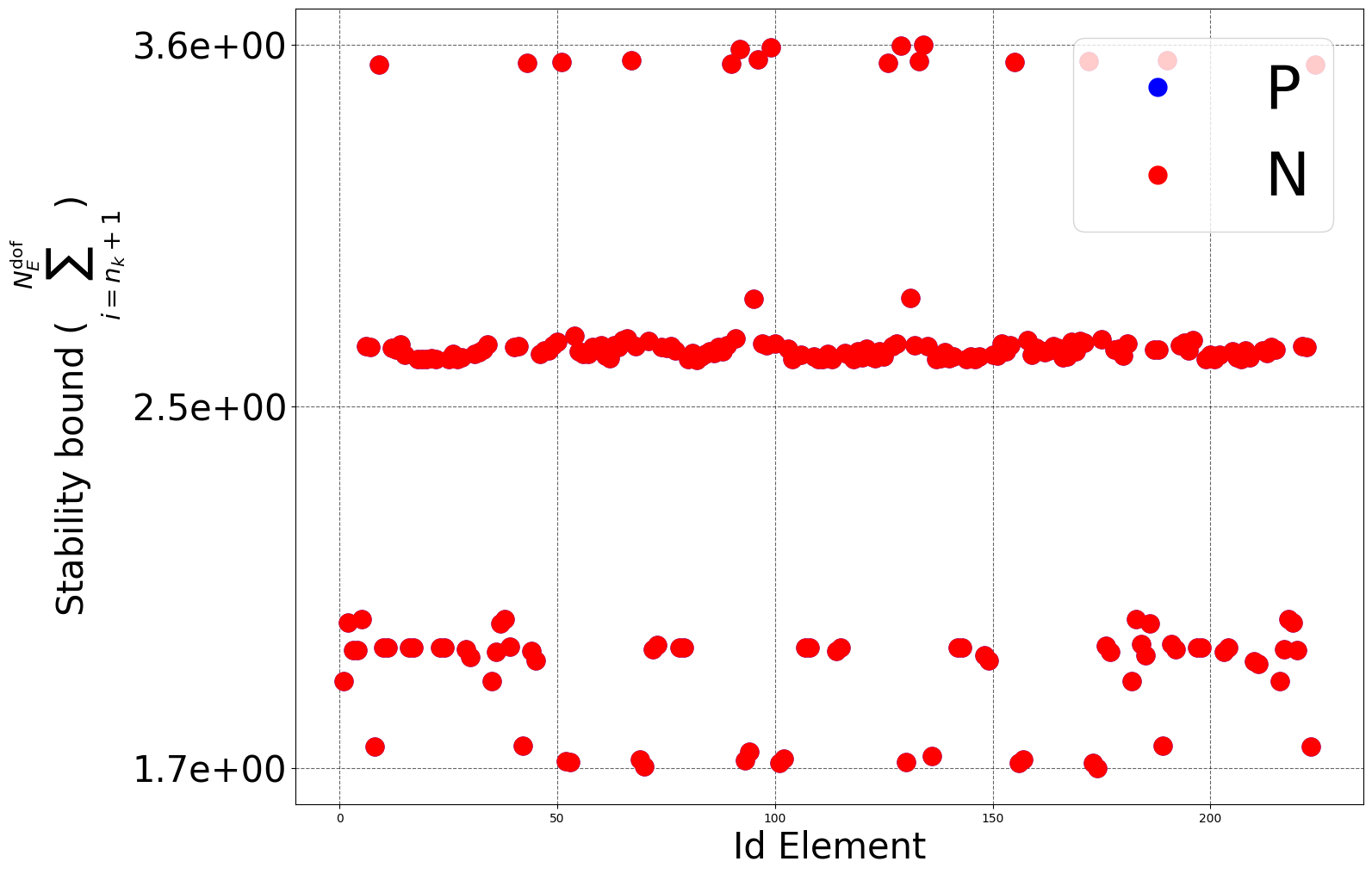}
        \caption{}
    \end{subfigure}
    \begin{subfigure}{0.3\textwidth}
        \includegraphics[width=1\textwidth]{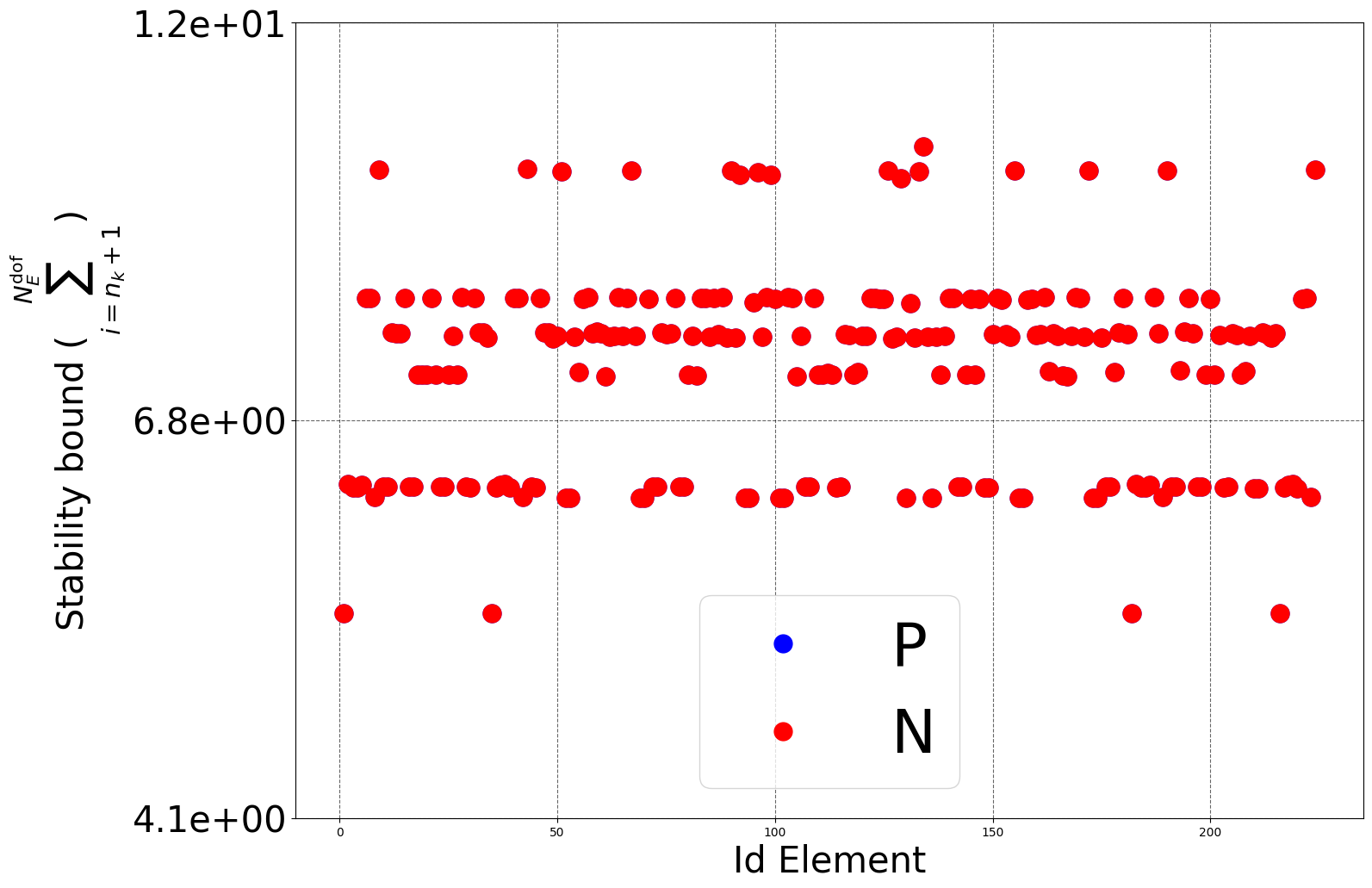}
        \caption{}
    \end{subfigure}
    \begin{subfigure}{0.3\textwidth}
        \includegraphics[width=1\textwidth]{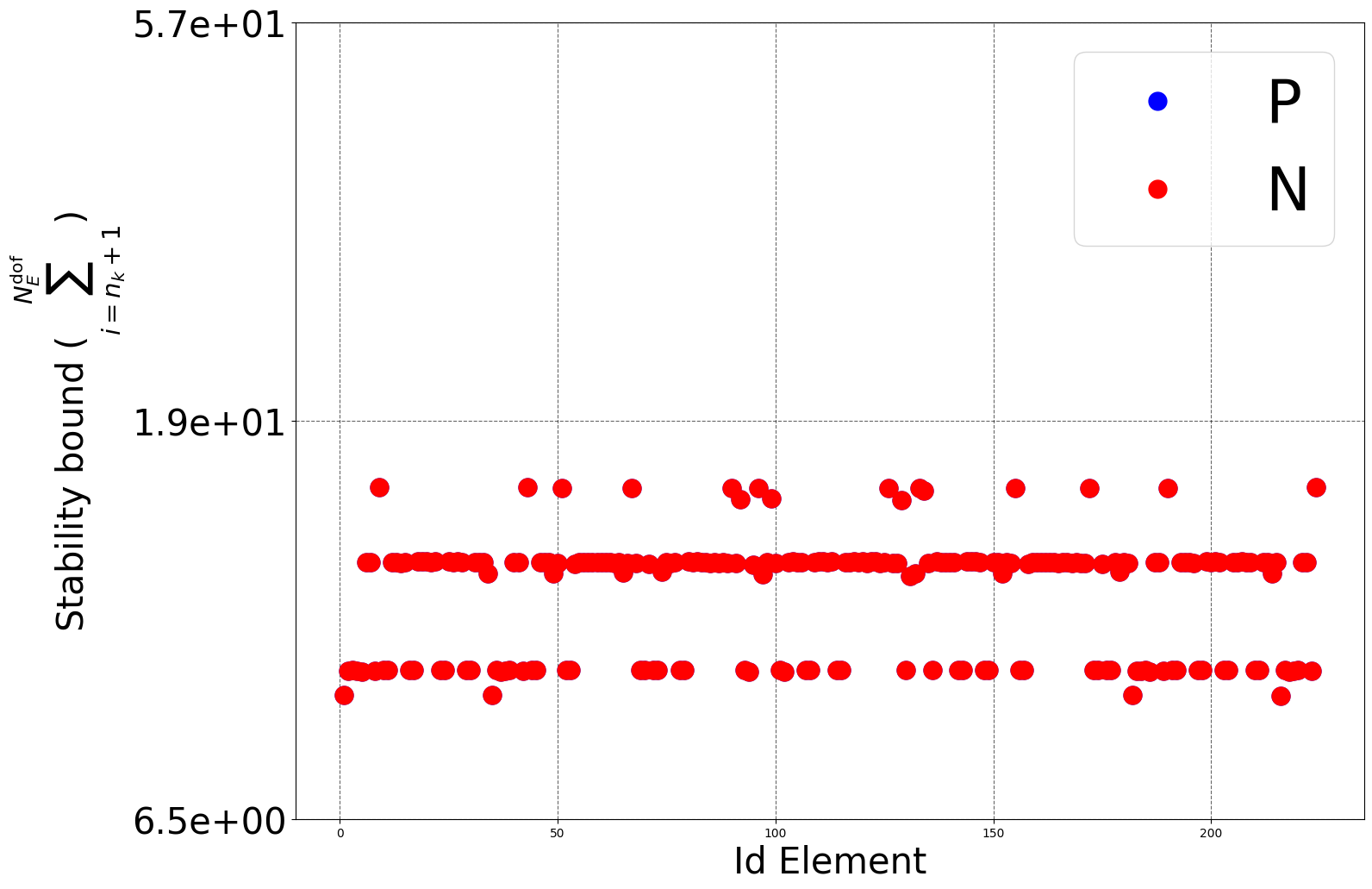}
        \caption{}
    \end{subfigure}
    \begin{subfigure}{0.3\textwidth}
        \includegraphics[width=1\textwidth]{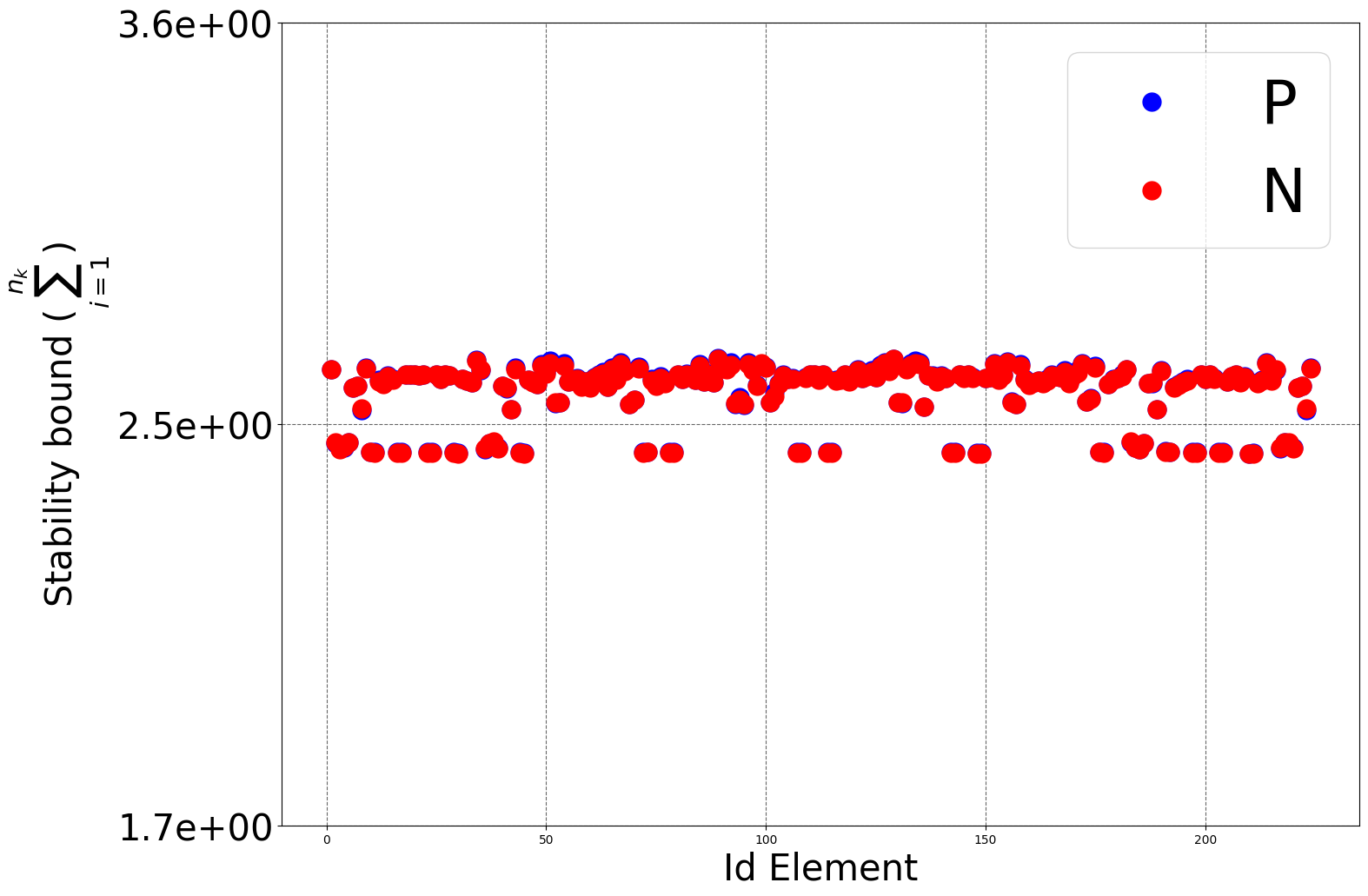}
        \caption{}
    \end{subfigure}
    \begin{subfigure}{0.3\textwidth}
        \includegraphics[width=1\textwidth]{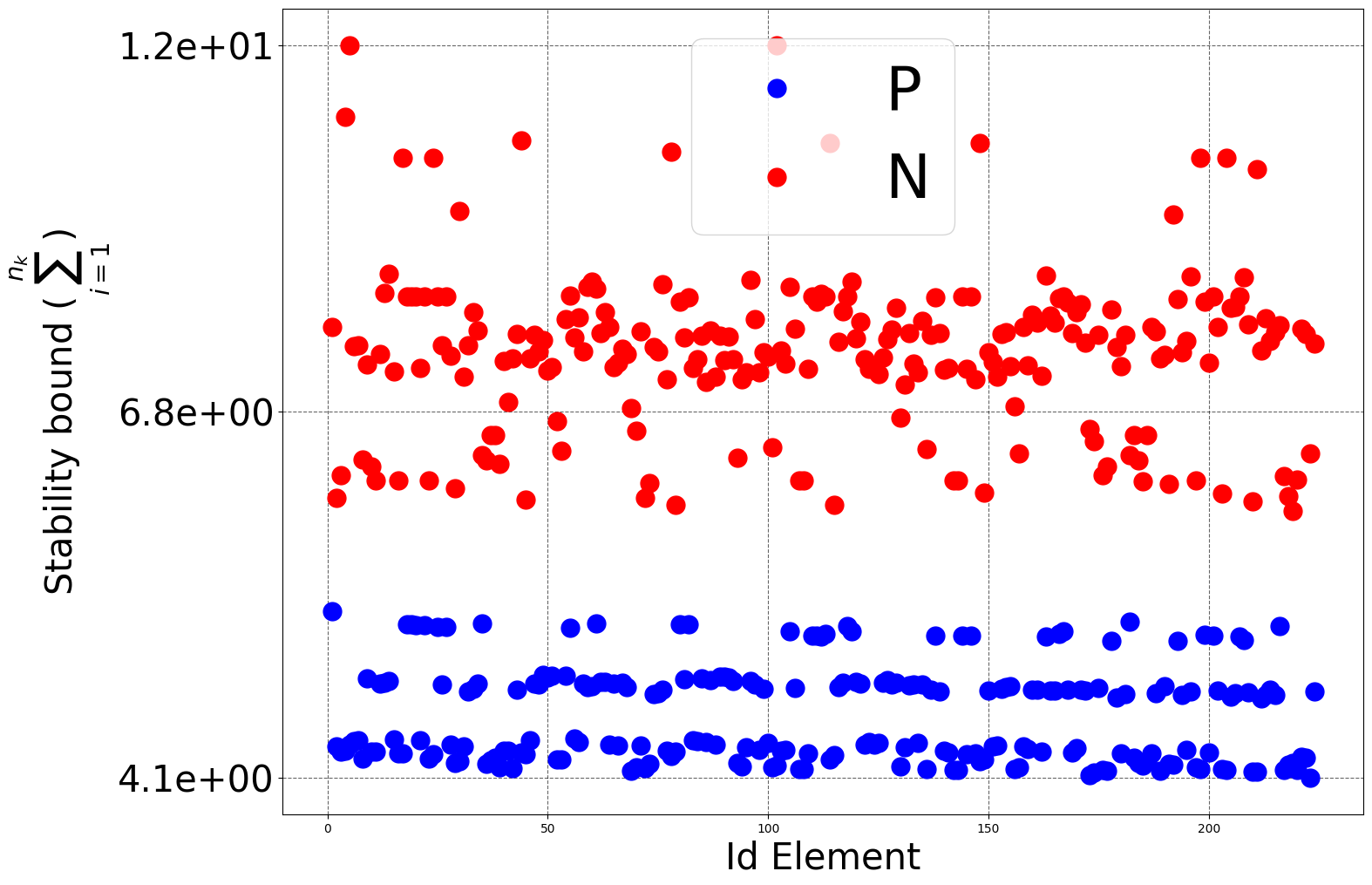}
        \caption{}
    \end{subfigure}
    \begin{subfigure}{0.3\textwidth}
        \includegraphics[width=1\textwidth]{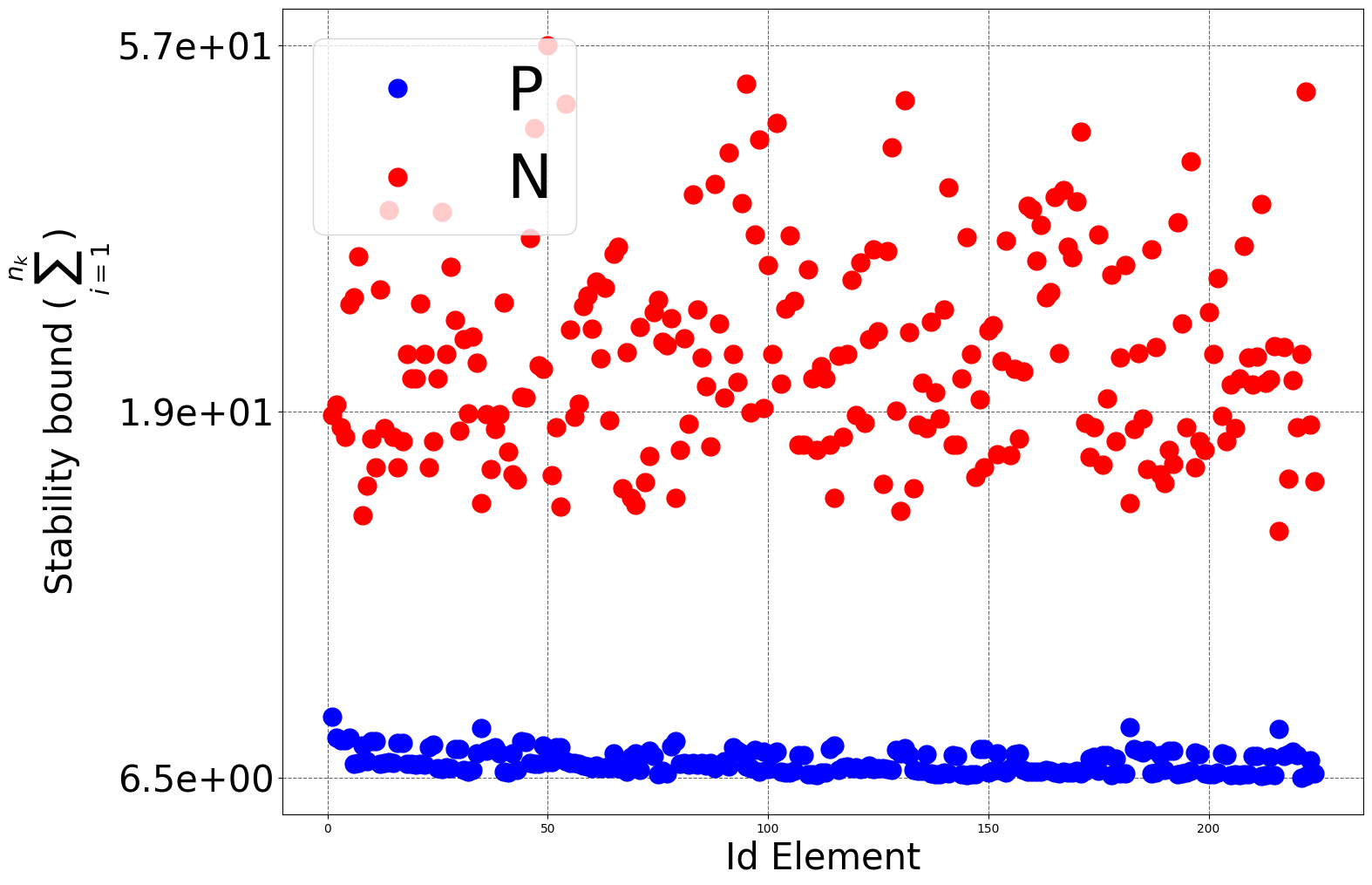}
        \caption{}
    \end{subfigure}
\caption{Test 2: Behavior of the partial stability bounds \eqref{eq:partial_stability_last} (first row) and \eqref{eq:partial_stability_first} (second row) for the elements in the last refinement of the Voronoi mesh \ref{mesh:C}. The bounds are associated with the original functions (red) and with the post-processed functions (blue). Each column represents a different order $k=1,2,3$ from left to right. Function type \ref{function:B}.}
    \label{fig:test2:stab_bound_ft2_split}
\end{figure}

To conclude, we briefly analyze the behavior of the stability bounds \eqref{eq:function_regularity} and \eqref{eq:function_regularity_post}, i.e., before and after applying {\postprocess}, with respect to the condition number of the matrix $\mathbf{D}$ defined in \eqref{eq:d_matrix}, for the set of elements in the last refinement of the Voronoi mesh \ref{mesh:C}. Figures~\ref{fig:test2:stab_bound_ft1} and \ref{fig:test2:stab_bound_ft2} report such behavior for function types \ref{function:A} and \ref{function:B}, resp., and for the polynomial degrees $k=1,2,3$, by clustering elements in accordance with their number of edges. In this case, we avoid reporting the behavior of function type \ref{function:C}, since this is very similar to \ref{function:B}.

Observe that such a Voronoi mesh only contains pentagons, hexagons and heptagons, and three clusters are easily recognizable in almost all the plots. Moreover, due to the higher number of vertices (and DOFs), the stability bounds for the heptagons are naturally the largest and the ones for the pentagons are generally the smallest. Nonetheless, all the stability bounds are of the same order of magnitude for both function types \ref{function:A} and \ref{function:B}, with and without post-processing. It is interesting to note that the post-processing acts as a stabilizer for the functions of type \ref{function:B}, because the related clusters are better divided and the stability bounds are smaller after the post-processing. This phenomenon is even more evident for higher-order basis functions.

To better investigate such a phenomenon, in Figures \ref{fig:test2:stab_bound_ft1_split} and \ref{fig:test2:stab_bound_ft2_split} for function types \ref{function:A} and \ref{function:B}, resp., we report the stability bound restricted to the sum of the functions modified by such a procedure (the first $n_k$ functions, according to the ordering in Section \ref{sec:pol_reproducibility}) and the unchanged ones (the remaining $\Ndof[E]-n_k$ functions). 

More precisely, in both figures, in the first row we show the stability bounds restricted to the last $\Ndof[E]-n_k$ functions, i.e.
\begin{equation}
\sum_{j=n_k+1}^{\Ndof[E]} \norm[\leb{2}{E}]{\nabla \bvarphi_{j,E}}, \qquad\text{and}\qquad \sum_{j=n_k+1}^{\Ndof[E]} \norm[\leb{2}{E}]{\nabla \pvarphi_{j,E}}.
\label{eq:partial_stability_last}
\end{equation}
Since these functions are such that $\bvarphi_{j,E}=\pvarphi_{j,E}$, for all $j=n_k+1,\dots,\Ndof[E]$, and all elements, all the blue dots are perfectly overlapped by their red counterpart. Moreover, the three clusters highlighted in the previous figures are always clearly recognizable in all considered cases.

When considering the bounds restricted to the first $n_k$ functions (second row in each figure), i.e. the quantities 
\begin{equation}
\sum_{j=1}^{n_k} \norm[\leb{2}{E}]{\nabla \bvarphi_{j,E}}, \qquad\text{and}\qquad \sum_{j=1}^{n_k} \norm[\leb{2}{E}]{\nabla \pvarphi_{j,E}},
\label{eq:partial_stability_first}
\end{equation}
the points associated with the original functions are not visibly clustered in three groups, whereas post-processed functions are better separated. In this case, the post-processed functions differ from the initial ones, so the blue dots and the red ones do not coincide (except for functions \ref{function:A} and $k=1$ because the corresponding local space already contains all polynomials of order 1). This effect is less visible even when \ref{function:B} and $k=1$ are considered, because the polynomials of order 1 can be accurately approximated by functions in the B-NAVEM space. 

Moreover, we remark that the two stability bounds are very close when the functions of type \ref{function:A} are used, while the post-processing visibly reduces the noise and lessens such bounds for functions of type \ref{function:B}. We can thus conclude that the QR algorithm selects the functions with more scattered stability bounds, and the post-processing generates functions with smaller oscillations when applied to functions of type \ref{function:B}.

\section{Conclusion}\label{sec:conclusion}
In this paper, we propose an algebraic algorithm to construct an elemental function space on a generic polygonal element $E$ that includes the polynomial space $\Poly{k}{E}$ of a given order $k$.
The algorithm takes as input a set of basis functions that span a space of dimension at least equal to $\dim \Poly{k}{E}$, and constructs a new set of functions that preserve some structural properties of the original set, such as continuity and regularity, while spanning a space that includes all the polynomials in $\Poly{k}{E}$.

The algorithm acts only at the elemental level, and no information about the global mesh is required. Nevertheless, the proposed strategy restores optimal polynomial convergence rates in both $\leb{2}{}$ and $\sob{1}{}$ on the entire domain. The effectiveness of the algorithm is demonstrated by using it to recover optimal convergence rates in the higher-order extension of the B-NAVEM method, introduced here for the first time. Theoretical analysis, showcasing the accuracy of the B-NAVEM method before and after the application of the proposed algorithm, proves how the proposed algorithm removes the limitation imposed by neural network accuracy on the convergence errors.

The theoretical results are paired with a set of numerical experiments confirming the expected properties and investigating the effect of the chosen quadrature rule. The numerical experiments confirm that, when all the operators are accurately integrated, the theoretical optimal polynomial convergence rate is attained numerically. 

Finally, we remark that the proposed procedure is independent of the spatial dimension of the problem, and that the theoretical framework extends naturally to the three-dimensional case. However, since three-dimensional simulations typically introduce considerable additional complexity, the numerical investigation of the 3D setting is left for future works.

\section*{Acknowledgments}

The authors are members of INdAM-GNCS.

The author S.B. kindly acknowledges partial financial support provided by European Union through project Next Generation EU, M4C2, PRIN 2022 PNRR project P2022BH5CB\_001 ``Polyhedral Galerkin methods for engineering applications to improve disaster risk forecast and management: stabilization-free operator-preserving methods and optimal stabilization methods''. The author G.T. kindly acknowledges partial financial support provided by INdAM-GNCS Project ``Metodi numerici politopali stabilization-free e neural-based per problemi accoppiati e non lineari'' (CUP: E53C25002010001).

\appendix
\renewcommand{\thesection}{Appendix \Alph{section}}

\section{A more efficient and stable factorization for the post-processing phase}\label{app:qr}

In Section~\ref{sec:pol_reproducibility}, we adopt the LU factorization to produce a polynomial-complete basis. This choice was natural: when performing the computations from scratch, computations are exactly the same as performed in the LU factorization. Nevertheless, before computing the LU factorization, the QR-factorization with column pivoting is needed to properly order the degrees of freedom and obtain a final block $\mathbf{D}(:, 1:n_k)$ that is non-singular. This factorization can be used in the post-processing in place of LU to avoid computing another factorization. By proceeding in this way, we factorize the matrix $\mathbf{D} \in \R^{n_k \times \Ndof[E]}$ defined by following a natural ordering for the degrees of freedom as
\begin{equation}
    \mathbf{D} \mathbf{P} = \mathbf{Q} \begin{bmatrix}
        \mathbf{R} & \overline{\mathbf{R}},
    \end{bmatrix}
\end{equation}
where $\mathbf{P}$ represent the matrix of column pivoting, $\mathbf{Q}\in \R^{n_k \times n_k}$ is an orthonormal matrix, i.e. $\mathbf{Q}^T \mathbf{Q} = \mathbf{I}$, and the first block of the matrix $\begin{bmatrix}
    \mathbf{R} & \overline{\mathbf{R}},
\end{bmatrix} \in \R^{n_k \times \Ndof[E]}$, i.e. $\mathbf{R}  \in \R^{n_k \times n_k}$, is a upper triangular matrix. The matrix $\mathbf{D} \mathbf{P}$ can be split as
\begin{equation*}
    \mathbf{D} \mathbf{P} = \begin{bmatrix}
        \mathbf{A} & \mathbf{B}
    \end{bmatrix},
\end{equation*}
where $\mathbf{A} = \mathbf{Q}\mathbf{R}$ coincides with the coefficient matrix used in \eqref{eq:postproc_initial_system}, whereas $\mathbf{B} = \mathbf{Q} \overline{\mathbf{R}}$ collects the remaining columns. The final system \eqref{eq:postproc_initial_system}, i.e.
\begin{equation*}
    \mathbf{A} \begin{bmatrix}
        \pvarphi_{1,E}(\xx)\\
        \vdots\\
        \pvarphi_{n_k,E}(\xx)
    \end{bmatrix} = \begin{bmatrix}
         m_1(\xx) \\
        \vdots \\
         m_{n_k}(\xx) \\
    \end{bmatrix} - \mathbf{B}\begin{bmatrix}
          \bvarphi_{n_k + 1,E}(\xx)\\
        \vdots \\
        \bvarphi_{\Ndof[E],E}(\xx)\\
    \end{bmatrix}
\end{equation*}
can thus be rewritten as
\begin{equation}
    \mathbf{R} \begin{bmatrix}
        \pvarphi_{1,E}(\xx)\\
        \vdots\\
        \pvarphi_{n_k,E}(\xx)
    \end{bmatrix} = \mathbf{Q}^T\begin{bmatrix}
         m_1(\xx) \\
        \vdots \\
         m_{n_k}(\xx) \\
    \end{bmatrix} - \overline{\mathbf{R}}\begin{bmatrix}
          \bvarphi_{n_k + 1,E}(\xx)\\
        \vdots \\
        \bvarphi_{\Ndof[E],E}(\xx)\\
    \end{bmatrix}.
    \label{eq:qr_stystem}
\end{equation}
Let us prove that this factorization, besides being more efficient and stable, allows us to preserve the Lagrange property of the basis functions and thus their continuity, as LU does.
Let us write \eqref{eq:qr_stystem} component-wise as
\begin{equation}
    \sum_{i=1}^{n_k} R_{\alpha i} \pvarphi_{i,E} = \sum_{\beta=1}^{n_k} (Q^T)_{\alpha \beta} m_{\beta} - \sum_{j=n_k +1}^{\Ndof[E]} \overline{R}_{\alpha j} \bvarphi_{j,E}, \quad \forall \alpha=1,\dots,n_k,
\end{equation}
where the DOFs ordering is now defined by QR.
By assuming that $\dof_{j}(\bvarphi_{i,E}) = \delta_{ij}$ for all $i,j = n_k + 1,\dots,\Ndof[E]$, we have that
\begin{equation*}
    \dof_{\gamma}(\pvarphi_{i,E}) = \delta_{\gamma i} \quad \forall \gamma = 1,\dots,\Ndof[E]\quad \forall i = n_k + 1,\dots, \Ndof[E],
\end{equation*}
since $\pvarphi_{i,E} \coloneq \bvarphi_{i,E}$. Moreover,
we can prove the following.
\begin{itemize}
    \item For all $\gamma = n_k + 1, \dots, \Ndof[E]$, 
    \begin{align*}
        \sum_{i=1}^{n_k} R_{\alpha i} \dof_{\gamma}(\pvarphi_{i,E}) &= \sum_{\beta=1}^{n_k} (Q^T)_{\alpha \beta} \dof_{\gamma}(m_{\beta}) - \sum_{j=n_k +1}^{\Ndof[E]} \overline{R}_{\alpha j} \dof_{\gamma}(\bvarphi_{j,E}) \\
        &=\sum_{\beta=1}^{n_k} (Q^T)_{\alpha \beta} D_{\beta \gamma} - \sum_{j=n_k +1}^{\Ndof[E]} \overline{R}_{\alpha j} \delta_{j \gamma}\\
        &= \overline{R}_{\alpha \gamma} - \overline{R}_{\alpha \gamma} = 0.
    \end{align*}
    We observe that, in this case, we are actually solving a problem in the form
    \begin{equation}
    \mathbf{R} \begin{bmatrix}
        \dof_{\gamma}(\pvarphi_{1,E})\\
        \vdots\\
        \dof_{\gamma}(\pvarphi_{n_k,E})
    \end{bmatrix} = \begin{bmatrix}
         0 \\
        \vdots \\
         0 \\
    \end{bmatrix}.
\end{equation}
    By the non-singularity of the coefficient matrix $\mathbf{R}$, we obtain that
    \begin{equation*}
        \dof_{\gamma}(\pvarphi_{i,E}) = 0 \quad \forall \gamma = n_k + 1 ,\dots, \Ndof[E]\quad \forall i =1,\dots,n_k. 
    \end{equation*}
    \item For $\gamma = 1,\dots,n_k$,
    \begin{align*}
        \sum_{i=1}^{n_k} R_{\alpha i} \dof_{\gamma}(\pvarphi_{i,E}) &= \sum_{\beta=1}^{n_k} (Q^T)_{\alpha \beta} \dof_{\gamma}(m_{\beta}) - \sum_{j=n_k +1}^{\Ndof[E]} \overline{R}_{\alpha j} \dof_{\gamma}(\bvarphi_{j,E}) \\
        &=\sum_{\beta=1}^{n_k} (Q^T)_{\alpha \beta} D_{\beta \gamma} - \sum_{j=n_k +1}^{\Ndof[E]} \overline{R}_{\alpha j} \delta_{j \gamma} = R_{\alpha \gamma} .
    \end{align*}
    We observe that, in this case, we are actually solving a problem in the form
    \begin{equation}
    \mathbf{R} \begin{bmatrix}
        \dof_{1}(\pvarphi_{1,E}) & \dots &  \dof_{n_k}(\pvarphi_{1,E})\\
        \vdots & \ddots & \vdots\\
        \dof_{1}(\pvarphi_{n_k,E}) & \dots & \dof_{n_k}(\pvarphi_{n_k,E})
    \end{bmatrix} = \mathbf{R}.
    \end{equation}
    By the non-singularity of $\mathbf{R}$, we obtain that the solution is represented by the identity matrix, i.e.
    \begin{equation*}
        \dof_{\gamma}(\pvarphi_{i,E}) = \delta_{\gamma i} \quad \forall \gamma, i =1,\dots,n_k. 
    \end{equation*}
\end{itemize}

Thus, the basis functions satisfy the Lagrange property, and, as a consequence, they are continuous across the domain (see Proposition~\ref{prop:continuity}).

\bibliographystyle{IEEEtranDOI}
\bibliography{biblio.bib}

\end{document}